\documentclass{article}
\usepackage{arxiv}
\usepackage{amsmath,amssymb,mathtools,wasysym}    % Math packages
\usepackage{hyperref}
\usepackage{xcolor}
\usepackage[capitalise]{cleveref}
\usepackage{regexpatch}
\usepackage{algorithm}
\usepackage{algorithmic}
\usepackage{ulem} % for striking out with \sout command
\usepackage{tikz}
\usepackage{fontawesome5} % \faGithubSquare
\usepackage{parskip}
\usepackage{booktabs}
\usepackage{calc,intcalc}
\usepackage[list=true, font=large, labelfont=bf,labelformat=brace, position=top]{subcaption}
\usepackage[numbers,sort&compress]{natbib}
\usepackage{empheq}
\usepackage{ifthen}
\usepackage{pdflscape}
\usepackage{siunitx}      %units \num{1e-10} \SI{1e-10}{\meter\per\second}
\usepackage{blkarray,bigstrut} % for block array in lgt_data
\usepackage{wrapfig}
\usepackage{grffile}        % unknown graphic extension error
\usepackage{svg}            % \includesvg{figure1.svg}
\usetikzlibrary{calc}
\usepackage{import}
\usetikzlibrary[pgfplots.groupplots]
\hypersetup{
    colorlinks,
    linkcolor={black},
    citecolor={blue!30!black},
    urlcolor={black}
}
\usetikzlibrary{shapes.misc}

\usepackage{pgfplots}
  \pgfplotsset{compat=newest}
  \usetikzlibrary{plotmarks}
  \usetikzlibrary{arrows.meta}
  \usepgfplotslibrary{patchplots}
  \usepackage{grffile}
  \usepackage{amsmath}

  \pgfplotsset{plot coordinates/math parser=false}
  \newlength\figureheight
  \newlength\figurewidth

\pgfplotsset{every axis/.append style={
                    label style={font=\scriptsize},
                    tick label style={font=\scriptsize},
                    legend style={font=\scriptsize}
                    }}

\usetikzlibrary{plotmarks}
\usepackage{scalerel}
\usetikzlibrary{svg.path}
\definecolor{orcidlogocol}{HTML}{A6CE39}
\tikzset{
  orcidlogo/.pic={
    \fill[orcidlogocol] svg{M256,128c0,70.7-57.3,128-128,128C57.3,256,0,198.7,0,128C0,57.3,57.3,0,128,0C198.7,0,256,57.3,256,128z};
    \fill[white] svg{M86.3,186.2H70.9V79.1h15.4v48.4V186.2z}
                 svg{M108.9,79.1h41.6c39.6,0,57,28.3,57,53.6c0,27.5-21.5,53.6-56.8,53.6h-41.8V79.1z M124.3,172.4h24.5c34.9,0,42.9-26.5,42.9-39.7c0-21.5-13.7-39.7-43.7-39.7h-23.7V172.4z}
                 svg{M88.7,56.8c0,5.5-4.5,10.1-10.1,10.1c-5.6,0-10.1-4.6-10.1-10.1c0-5.6,4.5-10.1,10.1-10.1C84.2,46.7,88.7,51.3,88.7,56.8z};
  }
}
\newcommand{\orcid}[1]{\href{https://orcid.org/#1}{\mbox{\scalerel*{
\begin{tikzpicture}[yscale=-1,transform shape]
\pic{orcidlogo};
\end{tikzpicture}
}{|}}}}

\input{commands.tex}
\title{Front Transport Reduction for Complex Moving Fronts}
\def\subtitle{Nonlinear Model Reduction for an Advection-Reaction-Diffusion Equation with a Kolmogorov--Petrovsky--Piskunov Reaction Term}
\def\shorttitle{Front Transport Reduction for Complex Moving Fronts}
\begin{document}


\author{
Philipp Krah \orcid{0000-0001-8982-4230},\\
Technische Universität Berlin,\\
Institute of Mathematics,\\
Straße des 17. Juni 136,\\
10623 Berlin, Germany\\
\email{krah@math.tu-berlin.de}
\And
Steffen Büchholz,\\
Technische Universität Berlin,\\
Institute of Fluid Mechanics and Engineering Acoustics,\\
Müller-Breslau-Str. 15,\\
10623 Berlin, Germany
\And
Matthias Häringer,\\
Technische Universität München,\\
Institute of Thermofluiddynamics,\\
Boltzmannstr. 15,\\
85747 Garching, Germany
\And
Julius Reiss \orcid{0000-0003-3692-5390},\\
Technische Universität Berlin,\\
Institute of Fluid Mechanics and Engineering Acoustics,\\
Müller-Breslau-Str. 15,\\
10623 Berlin, Germany
}

\date{Received: date / Accepted: date}

\maketitle

\begin{abstract}

This work addresses model order reduction for complex moving fronts, which are transported by advection or through a reaction-diffusion process. Such systems are especially challenging for model order reduction since the transport cannot be captured by linear reduction methods. Moreover, topological changes, such as splitting or merging of fronts pose difficulties for many nonlinear reduction methods and the small non-vanishing support of the underlying partial differential equations dynamics makes most nonlinear hyper-reduction methods infeasible.
We propose a new decomposition method together with a hyper-reduction scheme that addresses these shortcomings. The decomposition uses a level-set function to parameterize the transport and a nonlinear activation function that captures the structure of the front. This approach is similar to autoencoder artificial neural networks, but additionally provides insights into the system, which can be used for efficient reduced order models. We make use of this property and are thus able to solve the advection equation with the same complexity as the POD-Galerkin approach while obtaining errors of less than one percent for representative examples. Furthermore, we outline a special hyper-reduction method for more complicated advection-reaction-diffusion systems. 
The capability of the approach is illustrated by various numerical examples in one and two spatial dimensions, including real life applications to a two-dimensional Bunsen flame.
\end{abstract}

\newpage
%\tableofcontents

% !TEX root = ../wPOD.tex

\section{Introduction}

This article addresses model order reduction for reactive flows.
These flows often exhibit sharp fronts, like flames, which makes their simulation computational expensive.
This suggests applying model reduction for reducing simulation costs. 
However, classical model order reduction methods fail \cite{HuangDuraisamyMerkle2018} due to the sharp, moving fronts, that pose challenges for reducing and predicting new system states. This manuscript addresses these issues by presenting a new decomposition method together with efficient strategies to evaluate the dynamics of the reduced system.
For our study we use advection-reaction-diffusion systems (ARD) with a nonlinear Kolmogorov--Petrovsky--Piskunov (KPP) reaction term, as the complex front dynamics with topology changes of such systems feature essential difficulties for MOR, while the analysis is simplified because the reacting quantity is scalar and bounded. 

% This article addresses model order reduction for complex moving fronts, as they often occur in advection-reaction-diffusion systems (ARD). Analyzing such systems is of importance for various fields in science. For example in the simulation and optimization of combustion processes \cite{GrayLemkeReissPaschereitSesterhennMoeck2017}, neutron diffusion theory \cite{Kenneth1953}, studying spread of epidemic fronts \cite{Kendal1965}, forest fires \cite{Perry1998} or financial risk analysis \cite{BinderJadhavMehrmann2021}.
Model order reduction (MOR) has been studied for various ARD systems \cite{Grepl2012,BlackSchulzeUnger2021,KrahSrokaReiss2020,SwischukKramerHuangWillcox2020,BinderJadhavMehrmann2021}. In this study we focus on systems that exhibit locally one-dimensional traveling fronts.
The compact support of the moving fronts is challenging for linear reduced basis methods, such as the proper orthogonal decomposition (POD).
The POD approximates a set of snapshots $q(\vec{x},t_i), i=1,\dots,\Ntime$ by separation of variables
\begin{equation}
    \label{eq:POD-ansatz}
    q(\vec{x},t) \approx  \sum_{k=1}^r \hat{a}_k(t) \hat{\psi}_k(\vec{x})\,,
\end{equation}
with help of time amplitudes $\hat{a}_k(t)$ and spatial modes $\hat{\psi}_k(\vec{x})$, computed by a singular value decomposition (SVD). Unfortunately, moving fronts with sharp gradients significantly slow down the convergence of \cref{eq:POD-ansatz}. This has been numerically investigated for reactive flows \cite{HuangDuraisamyMerkle2018} and is theoretically quantified with help of the Kolmogorov
$n$-width in \cite{OhlbergerRave2015,GreifUrban2019}.

The convergence can be improved by compensating the transport, for which many authors use one-to-one mappings \cite{ReissSchulzeMehrmann2018,DecomposeBlackSchulzeUnger2021,Reiss2021,FedeleAbessiRoberts2015,RowleyKevrekidisMardsen2003} to align the front onto a reference frame in which the moving front is stationary. This allows to efficiently decompose the temporal variation of the front shape into  few basis functions.  On the downside, however, it assumes that
the transport dependent movement is known \cite{ReissSchulzeMehrmann2018,Reiss2021} or at least a sufficiently smooth function in time \cite{DecomposeBlackSchulzeUnger2021}, which is easy to parametrize and itself independent of $\vec{x}$. Unfortunately, for complicated transports, where fronts may split or merge, this approach does not work, because no smooth one-to-one mapping exists.

In this work, we therefore follow a more direct approach, in which we make use of an auxiliary field $\phi(\vec{x},t)$, which parameterizes the transport efficiently, together with a shape function $f$ to retain the front shape:
\begin{equation}
    \label{eq:ftr-ansatz0}
    q(\vec{x},t) \approx  f(\phi(\vec{x},t)) \quad \text{s.~t.}\quad \phi(\vec{x},t) = \sum_{k=1}^r a_k(t) \psi_k(\vec{x})\,,
\end{equation}
The auxiliary field $\phi\colon\mathbb{R}^d\times[0,T]\to \mathbb{R}$ allows to embed the local one dimensional front movement into a $d$-dimensional transport. Since the transport is only parameterized locally, changes in the topology of the front surface can be captured. A similar approach was introduced in \cite{KrahSrokaReiss2020}, where $\phi$ was constructed with help of a signed distance function and the front function $f$ was determined from a fit to the reacting front. 
While improving the approximation, this was found to be not optimal, since the obtained signed distance function does not have to be of low rank. Here, we follow a similar approach, but we formulate an optimization problem to compute $\phi$. The resulting description \cref{eq:ftr-ansatz0} is called \textit{Front Transport Reduction} (FTR) in the remainder of this manuscript.
% \begin{problem}{\textit{Front Transport Reduction}}
%   \label{FTR-opt}
%   For a given snapshot matrix $\matr{Q}\in\mathbb{R}^{M\times N_t}$  with $\matr{Q}_{ij}=q(\vec{x}_i,t_j)\in[0,1]$ and nonlinear smooth monotone increasing function $f\colon \mathbb{R}\to [0,1]$, find a rank $r$ matrix $\matr{\Phi}\in\mathbb{R}^{M\times N_t}$, such that the error $\norm{\matr{Q}-\tilde{\matr{Q}}}_\mathrm{F}^2$ for $\tilde{\matr{Q}}_{ij}=f(\matr{\Phi}_{ij})$ is minimized.
% \end{problem}
 Due to the nonlinearly activated linear space created by the span of $\{\psi_k(\vec{x})\}_{k=1,\dots,r}$, this approach shows many parallels to artificial neural networks. It can be seen as the decoder part of a shallow autoencoder structure. While shallow autoencoders have been used in previous studies \cite{KimChoiWidemannZohdi2021}, we are explicitly incorporating the underlying physical assumptions and thereby obtain interpretable results of reduced variables. 

The second part of this manuscript addresses dynamical ROM predictions of ARD systems using the FTR ansatz \cref{eq:ftr-ansatz0}. Here, many different methods exist in the literature, which can be categorized into intrusive or non-intrusive reduced order models. \textit{Intrusive} refers to data models where the resulting predictions are based on the initial ARD model, whereas a purely data-driven, \textit{non-intrusive} model is based on additional assumptions such as smoothness in the reduced parameter space. Intrusive models project the original equation system on the reduced manifold, which is nonlinear in our approach. These so-called \textit{manifold Galerkin methods} have been used in combination with neural networks in \cite{CarlbergFarhatCortial2013} and with dynamical transformed modes in \cite{BlackSchulzeUnger2020}. Unfortunately, manifold Galerkin methods require special hyper-reduction schemes to gain speedups in the resulting ROM. Examples of these methods are the extended-ECSW scheme proposed by \cite{JainTiso2019}, the gappy-POD based GNAT procedure \cite{CarlbergFarhatCortial2013} first introduced for nonlinear manifolds in \cite{KimChoiWidemannZohdi2021} or the shifted DEIM algorithm in \cite{BlackSchulzeUnger2021}. The idea of all of these methods is to evaluate the nonlinear dynamics of the underlying system for a small number of points to determine the evolution of the parameters in the reduced space.
 Unfortunately, the extended-ECSW scheme \cite{JainTiso2019} and the GNAT procedure \cite{KimChoiWidemannZohdi2021,JainTiso2019}, cannot be used for ARD systems with sharp advected fronts, since they preselect a fixed set of points, but the dynamics are localized only near the moving front.
We discuss this problem and state a practical solution using a special hyper-reduction scheme, based on the reduced integration domain (RID) method \cite{Ryckelynck2005}. Furthermore, we examine the ability of the FTR mapping to predict new system states with the help of non-intrusive methods.

\subsubsection*{Structure of the Article}

The remainder of the article is structured into three parts. The first part, \cref{sec:dim-red-moving-fronts}, is dedicated to the FTR decomposition, where we motivate the decomposition and introduce two algorithms to solve the corresponding high dimensional optimization problem. While the first algorithm is based on iterative thresholding of singular values (see \cref{subsec:FTR}), the other algorithm uses artificial neural networks (see \cref{subsec:ML}). The algorithms are applied and compared for two synthetic examples in \cref{subsec:synthetic-examples} and one example of a two-dimensional (2D) ARD system with topology change in \cref{subsec:offline-ARD}.
In the second part, \cref{sec:non-lin-roms-moving-fronts}, we use the low dimensional description of the FTR to predict new system states via non-intrusive MOR (\cref{subsec:data-driven}) and intrusive MOR (\cref{subsec:manifoldgalerkin}). For the latter we propose a special hyper-reduction method. The resulting ROM is tested for 1D and 2D ARD systems in \cref{subsec:hyper-ftr-numerical-examples}.
Finally, we summarize our results in the last part, \cref{sec:conclusion_and_outlook}.
%, including advection \cref{subsec:MovingDisc} and advection with topology changes \cref{subsubsec:advection_with_topo_change}. Furthermore 

\subsubsection*{Nomenclature}
Matrices are denoted in capital letters with straight, bold font $\matr{A}\in\mathbb{R}^{M\times M}$ and vectors are denoted by $\vec{x}\in\mathbb{R}^M$.
Whenever a scalar function $f\colon \mathbb{R}\to\mathbb{R}$ is applied on a vector valued quantity $\vec{x}$, we assume pointwise operation on the entries of $\vec{x}=(x_1,\cdots,x_M)$, if not stated otherwise,  and write $f(\vec{x})$ instead of $(f(x_1),\cdots,f(x_M))$, similarly for matrices $f(\matr{A})$.
Furthermore, if $q(\vec{x},t,\mu)\in\mathbb{R}$ is the solution of a scalar PDE, its (discrete space) ODE counterpart is denoted by a vector $\vec{q}(t,\mu)\in \mathbb{R}^M$ containing the spatial values of $q$ in its components. Correspondingly, the snapshot matrix $\matr{Q}$ contains all time and parameter snapshots in its columns: $\matr{Q}=[\vec{q}(t_1,\mu_1), \vec{q}(t_2,\mu_1), \dots, \vec{q}(t_{\Ntime},\mu_P)]$. Partial derivatives in space and time are denoted by $\partial_x=\frac{\partial}{\partial x}, \partial_t=\frac{\partial}{\partial t}, \partial_{xx}=\frac{\partial^2}{\partial x^2}$ and $\dot{q}=\frac{\partial{q}}{\partial t}$.
% !TEX root = FTR.tex
\section{Dimension Reduction Methods for Complex Moving Reaction Fronts}
\label{sec:dim-red-moving-fronts}
In this section, we motivate why special nonlinear reduction methods are advantageous when decomposing reactive flows and we introduce the \textit{Front Transport Reduction} as an iterative thresholding algorithm in \cref{subsec:FTR} and as an \textit{autoencoder network} with one decoder layer in \cref{subsec:ML}.

\subsection{The need for a nonlinear decomposition approach for moving fronts.}
\label{subsec:need_for_non_linear_decomposition}
To motivate our decomposition approach, we consider advection-reaction-diffusion systems of the form:
\begin{equation}
\label{eq:react-diff-advect}
    \frac{\partial q}{\partial t}  = \vec{u}\cdot\nabla q + \diffconst\Delta q + R(q,\reactconst)\,.
\end{equation}
These systems describe how a quantity or reactant $q(\vec{x},t)$ spreads in space $\vec{x}\in \domain\subset\mathbb{R}^d,\, d>0$ over time $t\in[0,T]$. This spread can be caused by the advection with velocity $\vec{u}\in\mathbb{R}^d$ or an interplay between diffusion $\Delta q$ and reaction processes $R(q,\reactconst)$. For the sake of simplicity, we focus on the reaction-diffusion described by a nonlinear Kolmogorov--Petrovsky--Piskunov (KPP) reaction term $R(q,\reactconst)=\reactconst q^\alpha(q-1)$ with $\alpha>0, \reactconst>0$.
These systems exhibit traveling or pulsating fronts \cite{HadelerRothe1975,BerestyckiHamelNadirashvili2004,BerestyckiHamelRoques2004,Fischer1937} and without loss of generality one can assume that $d=1$ near the front, since the moving structures are locally one-dimensional%
%%%%%%%%footnote%%%%%%%
\footnote{Formally, this work makes extensive use of the physically justified assumption that the spatial variable $\vec{x}=(x_1,x_2\dots,x_d)$ of the reactant $q$ can be transformed to $\vec{x}'=(\phi',x_2',\dots,x_d')$, where $x_2',\cdots,x_d'$ are on a hyperplane tangential to the front of the traveling wave. On this hyperplane, all gradients in the equation vanish relative to the terms that are normal to the traveling wave. Therefore, the flow can be described by a one-dimensional equation in the variable $\phi'$ and we simply can rewrite $q(\vec{x}',t)=q(\phi',t)$ (see \citep[p.87]{PoinsotVeynante2005} for details). }
%%%%%%%%footnote%%%%%%%
\cite{PoinsotVeynante2005,Williams1985,Peters2001}.
Therefore, the solution of \cref{eq:react-diff-advect} can be simply transformed into a co-moving frame
\begin{equation}
\label{eq:f_of_phi}
q(\vec x,t) = f(\phi(\vec x,t))\,,
\end{equation}
where the \textit{front-profile} of the traveling wave is described by $f$ and $\phi(\vec x,t) = (\vec x - \vec\Delta(\vec{x},t))\cdot \vec{e}_{\vec{v}}$, the location of the front with respect to the direction $\vec{e}_{\vec{v}}=\vec{v}/\norm{\vec{v}}$ of the wave speed $\vec{v}$. For a one dimensional traveling wave this is illustrated in \cref{subfig:frontmovement}.
 The profile of the wave $f$ can be analytically computed with help of perturbation theory after transforming \cref{eq:react-diff-advect} into the co-moving frame (see for example \cite{TropkinaShchepakina2021}) or by fitting the front profile \cite{KrahSrokaReiss2020}.
 However, the wave speed $\vec{v}$ is the most complex part in typical applications, since it is coupled to an outer transport/velocity field $\vec{u}$ in \cref{eq:react-diff-advect} and an additional constant propagation speed $c^*$ of the reacting wave which depends on $R(q,\reactconst)$ (i.e. minimal propagation speed $c^*\ge2\sqrt{\diffconst R'(0,\reactconst)}$ for KPP nonlinearities \cite{HadelerRothe1975,BerestyckiHamelNadirashvili2004,BerestyckiHamelRoques2004}).

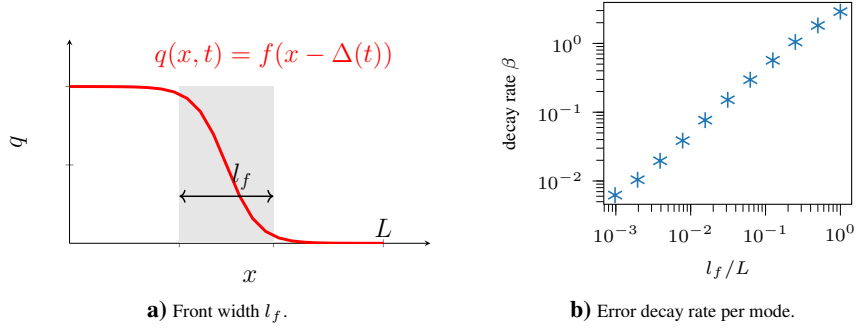
\begin{figure}
    \centering
    \begin{subfigure}[t]{0.4\textwidth}
        \centering
        \begin{tikzpicture}[thick,scale=0.7]
            \begin{axis}[
                unit vector ratio*=1 1 1,
                xmin=0, xmax=2.3,
                ymin=0, ymax=1.3,
                axis lines=center,
                %axis on top=true,
                domain=0:2,
                xtick={0.7,1.3,2},
                xticklabels={ },
                yticklabels={},
                x label style={at={(axis description cs:0.5,-0.1)},anchor=north},
                y label style={at={(axis description cs:-0.1,.5)},rotate=90,anchor=south},
                ylabel={\Large$q$},
                xlabel={\Large $x$},
                ]
                \fill[color=black!10] (axis cs:0.7,0) rectangle (axis cs:1.3,1);
                \addplot [mark=none,draw=red,ultra thick] {0.5-0.5*tanh((\x-1)*5)};
                \node [right, red] at (axis cs: 0.5,1.2) {\Large$q(x,t)= f(x-\Delta(t))$};

                %% Add the asymptotes
                \draw [thick,<->] (axis cs:0.7,0.3) -- (axis cs:1.3,0.3);
                \node [above] at (axis cs:1.1,0.3) {\Large $l_f$};
                \node [above] at (axis cs:2,0) {\Large $L$};
                %\draw [blue, dotted, thick] (axis cs:2,+1)-- (axis cs:0,+1);
            \end{axis}
        \end{tikzpicture}
        \caption{\scriptsize Front width $l_f$.}
        \label{subfig:front-graph}
    \end{subfigure}
    \begin{subfigure}[t]{0.4\textwidth}
        \centering
        \setlength{\figureheight}{0.7\linewidth}
        \setlength{\figurewidth}{0.8\linewidth}
        % This file was created by tikzplotlib v0.9.8.
\begin{tikzpicture}

\definecolor{color0}{rgb}{0.12156862745098,0.466666666666667,0.705882352941177}

\begin{axis}[
height=\figureheight,
log basis x={10},
log basis y={10},
tick align=outside,
tick pos=left,
width=\figurewidth,
x grid style={white!69.0196078431373!black},
xlabel={\(\displaystyle l_f/L\)},
xmin=0.000690533966002488, xmax=1.4142135623731,
xmode=log,
xtick style={color=black},
y grid style={white!69.0196078431373!black},
ylabel={decay rate \(\displaystyle \beta\)},
ymin=0.0045853111609863, ymax=3.94302679627954,
ymode=log,
ytick style={color=black}
]
\addplot [semithick, color0, mark=asterisk, mark size=3, mark options={solid}, only marks]
table {%
1 2.90031481762917
0.5 1.82569781833937
0.25 1.05043944769858
0.125 0.568241527543457
0.0625 0.2963262499751
0.03125 0.151466567685262
0.015625 0.0766309984637047
0.0078125 0.0386300068496256
0.00390625 0.0196561275302464
0.001953125 0.0104297085827299
0.0009765625 0.00623380767741204
};
\end{axis}

\end{tikzpicture}
        \caption{\scriptsize Error decay rate per mode.}
        \label{subfig:front-error}
    \end{subfigure}
    \caption{Relationship between characteristic length scale $l_f$ and approximation error of the POD. \Cref{subfig:front-graph} shows a wave front $f(x)=\sig(x l_f)$ with front width $l_f$ traveling along a distance $L$ with constant speed $c=L/T$ and shift $\Delta(t)=c t,\, t\in[0,T[$. \Cref{subfig:front-error} shows the decay rate $\norm{q-\tilde{q}_n}\sim e^{-\beta n}$ as a function of the relative front width $l_f/L$, when approximating $q$ with $n$ POD-modes.}
    \label{fig:mov-front1D}
\end{figure}

Dimensional analysis yields a definition of the thickness of the propagating front in terms of the fraction of the diffusion and propagation speed:
\begin{equation}
\label{eq:characteristic-length}
    l_f=\diffconst/c^*\le 0.5\sqrt{\diffconst/R'(0,\reactconst)}\,.
\end{equation}
This characteristic length scale  of the system is shown in \cref{fig:mov-front1D}.
In a linear projection-based MOR approach, $l_f$ plays an essential role, because its length is directly related to the success of the approximation. As already pointed out by \cite{GreifUrban2019,OhlbergerRave2015} for transport systems with vanishing front width ($l_f\to 0$), every front position is linear independent of the others and therefore equally important when defining a projection basis. Therefore, the typical exponential decay $\norm{q-\tilde{q}_n}\sim e^{-\beta n}$ of the approximation error is reduced to $\sim n^{-\frac{1}{2}}$, when increasing the dimension of the ROM-basis $n$. Since the authors \cite{GreifUrban2019,OhlbergerRave2015} give no general results for $l_f>0$, we quantify the decay of the approximation errors numerically in \cref{fig:mov-front1D}. It can be seen from \cref{subfig:front-error} that the error decay rate per POD mode diminishes if the traveling distance $L$ becomes large relative to the front width $l_f$, making a linear MOR approach impractical.
In order to compensate the transport, many studies use one-to-one mappings \cite{ReissSchulzeMehrmann2018,DecomposeBlackSchulzeUnger2021,Reiss2021,FedeleAbessiRoberts2015,RowleyKevrekidisMardsen2003}, which can not be used here, since reacting fronts may split or merge.

Here, we thus make explicit use of the underlying physical structure of ARD systems \cref{eq:f_of_phi}. For given snapshot data $\matr{Q}\in \mathbb{R}^{M\times N_t}$ with $\matr{Q}_{ij}=q(\vec{x}_i,t_j)\in[0,1]$ and front function $f\colon \mathbb{R}\to [0,1]$, the approach decomposes the data with help of the nonlinear mapping
\begin{equation}
    \label{eq:ftr-ansatz}
    q(\vec{x},t) \approx \tilde{q}(\vec{x},t) =  f(\phi(\vec{x},t)) \quad \text{s.~t.}\quad \phi(\vec{x},t) = \sum_{k=1}^r a_k(t) \psi_k(\vec{x})\,, r\ll N_t
\end{equation}
and a low rank field $\matr{\Phi}_{ij}=\phi(\vec{x}_i,t_j)$, that allows to embed the local one dimensional front movement into a $d$ dimensional transport. The idea is visualized in \cref{fig:FTR-1D}. Since the transport is only parameterized locally, changes in the topology of the front surface can be captured.
The decomposition goal is formulated as an optimization problem.
\begin{problem}{\textit{Front Transport Reduction}}
  \label{FTR-opt}
  For a given snapshot matrix $\matr{Q}\in\mathbb{R}^{M\times N_t}$  with $\matr{Q}_{ij}=q(\vec{x}_i,t_j)\in[0,1]$ and nonlinear smooth monotone increasing function $f\colon \mathbb{R}\to [0,1]$, find a rank $r$ matrix $\matr{\Phi}\in\mathbb{R}^{M\times N_t}$, such that the error $\norm{\matr{Q}-\tilde{\matr{Q}}}_\mathrm{F}^2$ for $\tilde{\matr{Q}}_{ij}=f(\matr{\Phi}_{ij})$ is minimized.
\end{problem}

\begin{figure}[htp!]
\centering
\begin{subfigure}[t]{0.45\textwidth}
        \centering
        \setlength{\figureheight}{0.5\linewidth}
        \setlength{\figurewidth}{1\linewidth}
        \input{levelset1Dv2.tex}
        \caption{Transported quantities}
        \label{subfig:frontmovement}
    \end{subfigure}
    \begin{subfigure}[t]{0.45\textwidth}
        \centering
        \setlength{\figureheight}{0.9\linewidth}
        \setlength{\figurewidth}{1\linewidth}
        % This file was created by matplotlib2tikz v0.7.4.
\begin{tikzpicture}[scale=0.8]

\definecolor{color0}{rgb}{1,0.647058823529412,0}

\begin{axis}[
height=\figureheight,
legend cell align={left},
legend style={at={(0.5,0.6)}, anchor=north, draw=none},
log basis y={10},
tick align=outside,
tick pos=left,
width=\figurewidth,
x grid style={white!69.01960784313725!black},
xlabel={rank \(\displaystyle r\)},
xmin=0, xmax=20,
xtick style={color=black},
y grid style={white!69.01960784313725!black},
ylabel={relative error},
ymin=4.54213436950621e-18, ymax=6.69645204306777,
ymode=log,
ytick style={color=black}
]
\addplot [blue, mark=*, mark size=2, mark options={solid}]
table {%
0 1
1 0.333307703774789
2 0.199954027121491
3 0.142791811095647
4 0.111027012985494
5 0.0908068272042107
6 0.0768032340419105
7 0.0665299577867899
8 0.0586707576941898
9 0.0524636522024691
10 0.0474369849968113
11 0.0432831983625682
12 0.0397931956823475
13 0.0368198763784882
14 0.0342567600498823
15 0.0320248853257144
16 0.0300644729534269
17 0.0283294474291522
18 0.0267837355609547
19 0.0253987040912931
20 0.0241513474298762
21 0.0230229812714542
22 0.0219982847067642
23 0.0210645870127596
24 0.020211329204153
25 0.0194296523606727
26 0.0187120792278299
27 0.0180522653322048
28 0.0174448025186014
29 0.0168850624490217
30 0.0163690708692308
31 0.0158934057818367
32 0.0154551143525037
33 0.0150516446111154
34 0.0146807889238899
35 0.014340636896266
36 0.0140295358836925
37 0.013746057681911
38 0.0134889702748358
39 0.0132572135707566
40 0.0130498795780921
41 0.0128661931241879
42 0.0127054977392726
43 0.0125672391273757
44 0.0124514240610197
45 0.0123563603073238
46 0.0122830307689976
47 0.0122308092701341
48 0.0121995376914025
49 0.00557113083247043
};
\addlegendentry{$\Vert \matr{Q} - \tilde{\matr{Q}}^\text{POD} \Vert_\text{F}^2$}
\addplot[color0, dashed, mark=triangle*, mark size=2, mark options={solid,rotate=90}]
table {%
0 1
1 0.983468842746643
2 2.98409190362933e-16
3 2.39087578108502e-16
4 1.71640175526394e-16
5 1.06455300209217e-16
6 9.95076345691556e-17
7 9.95076345691556e-17
8 9.95076345691556e-17
9 9.95076345691556e-17
10 9.95076345691556e-17
11 9.95076345691556e-17
12 9.95076345691556e-17
13 9.95076345691556e-17
14 9.95076345691556e-17
15 9.95076345691556e-17
16 9.95076345691556e-17
17 9.95076345691556e-17
18 9.95076345691556e-17
19 9.95076345691556e-17
20 9.95076345691556e-17
21 9.95076345691556e-17
22 9.95076345691556e-17
23 9.95076345691556e-17
24 9.95076345691556e-17
25 9.95076345691556e-17
26 9.95076345691556e-17
27 9.95076345691556e-17
28 9.95076345691556e-17
29 9.95076345691556e-17
30 9.95076345691556e-17
31 9.95076345691556e-17
32 9.95076345691556e-17
33 9.95076345691556e-17
34 9.95076345691556e-17
35 9.95076345691556e-17
36 9.95076345691556e-17
37 9.95076345691556e-17
38 9.95076345691556e-17
39 9.95076345691556e-17
40 9.95076345691556e-17
41 9.95076345691556e-17
42 9.95076345691556e-17
43 9.95076345691556e-17
44 9.95076345691556e-17
45 9.95076345691556e-17
46 9.95076345691556e-17
47 9.95076345691556e-17
48 9.95076345691556e-17
49 3.0416184978568e-17
};
\addlegendentry{$\Vert \matr{\Phi}- \tilde{\matr{\Phi}}^\text{POD} \Vert_\text{F}^2$}
\end{axis}

\end{tikzpicture}
        \caption{Relative POD approximation errors}
        \label{subfig:levelset1Dsigma}
    \end{subfigure}
    \caption{Illustration of the basic idea of the front transport reduction method. \Cref{subfig:frontmovement}: The FTR replaces the sharp traveling front structure $q$ (blue curves), by a level set function $\phi$ (orange lines) and a nonlinear mapping $f$ (indicated by the red arrow). Both quantities share locally the same transport. However, the level set field $\phi(x,t)=x-\Delta(t)$ is of low rank and can be therefore parameterized with only a few POD basis functions (here: $\{x,1\}$). \Cref{subfig:levelset1Dsigma}: The generated snapshot data $\matr{\Phi}_{ij}=\phi(x_i,t_j)$ can be approximated efficiently with the POD, compared to  $\matr{Q}_{ij}=q(x_i,t_j)$. }
    \label{fig:FTR-1D}
\end{figure}

Two possible algorithms that solve the optimization problem \ref{FTR-opt} are provided in the following sections.

\subsection{Front Transport Reduction via Iterative Thresholding of Singular Values}
\label{subsec:FTR}

A simple iterative algorithm to determine the auxiliary field $\matr{\Phi}\in \mathbb{R}^{M\times N_t}$ of the front transport reduction Problem \ref{FTR-opt} is stated in \cref{alg:iterFTR}.

Our algorithm is constructed by combining a gradient descent step (line \ref{alg:iterFTR-GD}) to minimize  $\norm{\matr{Q}-\tilde{\matr{Q}}}_\mathrm{F}^2$,
together with a rank-$r$ projection step of $\matr{\Phi}$ (line \ref{alg:iterFTR-thresholding}).
In the gradient descent step, the FTR residual
\begin{equation}
\mathcal{L}_\mathrm{FTR}(\matr{\Phi}) =\frac{1}{2} \norm{\matr{Q}-\tilde{\matr{Q}}}_\mathrm{F}^2  \quad \text{ with }\quad \tilde{\matr{Q}} = f(\matr{\Phi})
\end{equation}
is minimized in direction of the gradient $D_\matr{\Phi} \mathcal{L}_\mathrm{FTR}(\matr{\Phi})=f'(\matr{\Phi})\odot(f(\matr{\Phi})-\matr{Q})$. Here, $f'(\matr{\Phi}),f(\matr{\Phi})$ are element-wise operations of $f,f'$ on $\matr{\Phi}$.
Since $f$ is monotonically increasing, it is sufficient to replace $D_\matr{\Phi} \mathcal{L}_\mathrm{FTR}$ by $\matr{R}=f(\matr{\Phi})-\matr{Q}$ in line \ref{alg:iterFTR-GD}. Neglecting $f'(\matr{\Phi})$ in the gradient prevents a dying gradient for points where $f'(\Phi)\to 0$, i.e. $\abs{\matr{\Phi}_{ij}}\gg 0$.
Note that replacing the simple gradient descent step by a
quasi Newton method or a line search would not affect the convergence rate, since it is followed by a projection step (line \ref{alg:iterFTR-thresholding}), which is likely to
destroy the possible larger step of a more sophisticated method.

\begin{algorithm}[H]
\caption{FTR as iterativ thresholding}
\label{alg:iterFTR}
\begin{algorithmic}[1]
\REQUIRE $\matr{Q}\in\mathbb{R}^{M\times N_t}$ data $\matr{Q}_{ij}=q(\vec{x}_i,t_j)$, $\tau$ step size, $r$ rank
\STATE init $\matr{\Phi}^k = 0$
\WHILE{ not converged }
%\STATE {\color{black!60!white}($Y_{ij}^{n} = Y_{ij}^n-%\frac{1}{N_x}\sum_{k=1}^{N_x}Y_{kj}^n $)}
    \STATE residual $\matr{R} = f(\matr{\Phi}^k) -\matr{Q}$
    \STATE $\matr{\Phi}^{k+1/2} = \matr{\Phi}^{k} - \tau \matr{R}$ \label{alg:iterFTR-GD}
    \STATE decompose and truncate\\
                $\matr{\Phi}^{k+1} = \mathrm{svd}(\matr{\Phi}^{k+1/2},r)$ \label{alg:iterFTR-thresholding}
    \STATE $k \leftarrow k + 1$
\ENDWHILE
\RETURN $\matr{\Phi}^k$
\end{algorithmic}
\end{algorithm}

% A similar algorithm exists in the literature, known as iterative thresholding algorithms for the robust principal component analysis proposed in \cite{WrightGaneshRaoPengMa2009}, where the authors use a soft thresholding operator to accurately recover the low rank part of a possible low rank data field, which is distorted by noise.
% Therefore, one could exchange the hard rank $r$-threshold by a soft threshold
% in \cref{alg:iterFTR} line \ref{alg:iterFTR-thresholding}, which has not been pursued in this work.

The computational costs of \cref{alg:iterFTR} scale with the complexity of the singular value decomposition (SVD). For large systems it can be advantageous to use randomized- or  wavelet-techniques \cite{HalkoMartinssonTropp2011,KrahEngelsSchneiderReiss2020} to compute the SVD.

\newcommand{\enc}{g_\text{enc}}
\newcommand{\dec}{g_\text{dec}}
\subsection{Front Transport Reduction via Neural Autoencoder Networks}
\label{subsec:ML}
Another way to solve the optimization problem \ref{FTR-opt} is with the help of neural autoencoder networks, which are commonly used in dimensionality reduction \cite{VanDerMaatenPostmaHerik2009}. For a general introduction to neural autoencoder networks, we refer to \cite{Goodfellow-et-al-2016}.
Here, we briefly explain the concept and the specifications of our network.

An autoencoder tries to reproduce the input data, while squeezing it through an informational bottleneck. It consists of two parts, the
\begin{description}
    \item[\textbf{Encoder}] $\enc \colon \mathbb{R}^{M}\to\mathbb{R}^r,\, \vec{q}\mapsto \vec{a} = \enc(\vec{q}) $, mapping the input data $\vec{q}$ onto points $\vec{a}$ in a learned lower dimensional latent space and the
    \item[\textbf{Decoder}]
    $\dec \colon\mathbb{R}^{r}\to\mathbb{R}^{M},\vec{a} \mapsto \dec(\vec{a}) =  \tilde{\vec{q}}$, mapping the latent representation back to the input space.
\end{description}
The composition of the two parts
\begin{align*}
    \tilde{\vec{q}} = \dec(\enc(\vec{q}))
\end{align*}
defines the autoencoder.
 The task of the optimization procedure is now to determine $\dec, \enc$, such that the reconstruction error over the training data $\matr{Q}=[ \vec{q}_1, \dots, \vec{q}_{\Ntime}]$:
\begin{align*}
    \mathcal{L}_{\mathrm{FTR}} = \sum_{i=1}^{\Ntime}\Vert \vec{q}_i - \tilde{\vec{q}}_i\Vert^2_\mathrm{F} = \sum_{i=1}^{\Ntime}\Vert \vec{q}_i - \dec(\enc(\vec{q}_i))\Vert^2_\mathrm{F}
\end{align*}
is minimized.
After the network has been trained, the reduction is achieved as the dimension $r\ll M$ of the latent variables $\vec{a}_i=\enc(\vec{q}_i)\in\mathbb{R}^r$ is much smaller than the input dimension $M$. Therefore, the decoder $\vec{q}_i\approx\dec(\vec{a}_i)$
represents a reduced map of the high dimensional data contained in the columns of $\matr{Q}$.

In the training procedure, the functions $\enc,\dec$ are determined by trainable parameters of the network, called weights and biases. The networks are constructed by a composition of \define{layers}, $\enc=L_1\circ L_2 \circ \dots \circ L_N$. Usually, the layers of the network $L_n\colon \mathbb{R}^i \to \mathbb{R}^o$ are given by an affine linear mapping $\vec{x}\mapsto h_n(\matr{W}_n \vec{x}+\vec{b}_n)$, with weights $\matr{W}_n\in\mathbb{R}^{o,i}$ and biases $\vec{b}_n \in\mathbb{R}^o$ together with a predefined nonlinear function $h_n$. The choice of the input and output dimension $i,o$ in each layer, the activation function and the number of layers is called \define{architecture of the network}.

As the FTR-autoencoder network (FTR-NN) should implement the structure motivated in Problem \ref{FTR-opt}, we choose a special architecture. It consists of a single layer decoder, without bias
 \begin{align*}
     \tilde{\vec{q}}=\dec(\vec{a}) = f(\FTRmodes \vec{a})\,, \quad \FTRmodes \in \mathbb{R}^{M\times r}\,,
 \end{align*}
which is activated by the physics dependent front function $f$.
Here, the images of the linear part $\boldsymbol{\phi}_i= \FTRmodes\vec{a}_i$, with respect to $\vec{a}_i=\enc(\vec{q}_i)$  correspond to the columns of the discrete transport field $\matr{\Phi}=[\boldsymbol{\phi}_1,\dots, \boldsymbol{\phi}_{\Ntime}]$. Since the image of the linear part is represented by $\FTRmodes\in \mathbb{R}^{M \times r},\,r\ll M$ the resulting matrix is at most of rank $r$.

The encoder network consists of four convolutional layers, each followed by an exponential linear unit (ELU) and a batch normalization layer \cite{IoffeSzegedy2015}. After flattening the output, the convolutional layers are followed by two linear layers, where the first one is again followed by an ELU activation and a batch normalization layer. We apply a stride of two in all convolutional layers after the first, to downsample the spatial resolution of the input data. Further details of the architecture and training procedure can be found in \cref{appx:AutoencoderArchitecture}.

For the training of the FTR-NN, an additional smoothness constraint is added to the optimization goal $\mathcal{L}_\mathrm{FTR}$, which penalizes the non-smoothness of the columns $\vec{\psi}_n$ of $\FTRmodes\in \mathbb{R}^{M \times r}$
\begin{align}
    \label{eq:NN-reg}
   \mathcal{L}_\mathrm{smooth} = \smoothness\sum_{n=1}^{r}\frac{\Vert \matr{D} \vec{\psi}_n\Vert^2_\mathrm{F}}{\Vert \vec{\psi}_n \Vert^2_\mathrm{F}}.
\end{align}
Here, $\matr{D}\in\mathbb{R}^{M\times M}$ denotes the coefficient matrix of a forward finite difference, which is implemented as a convolution operation over the columns of $\FTRmodes$. For the examples in this manuscript $\smoothness=10^{-7}$, was found to be optimal.  The additional smoothness constraint allows for faster convergence of the network in the validation phase. The constraint is reasonable since the columns represent the transport field $\matr{\Phi}_{ij}=\phi(\vec{x}_i,t_j)$, which is assumed to be smooth.

\subsection{Synthetic Examples}
\label{subsec:synthetic-examples}
In this subsection, we provide two synthetic examples. The first example illustrates the application of the FTR to linear advection and compares the two decomposition methods outlined above with previous results \cite{KrahSrokaReiss2020}. The second example addresses topology changing fronts.

%%%%%%%%%%%%%%%%%%%%%%%%%%%%%%%%%%%%%%%%%%%%
\subsubsection{Linear Advection of a Disk}
\label{subsec:MovingDisc}
%%%%%%%%%%%%%%%%%%%%%%%%%%%%%%%%%%%%%%%%%%%%

\begin{figure}[htp!]
		  \centering
							 \setlength\figureheight{0.25\linewidth}
							 \setlength\figurewidth{0.25\linewidth}
			  % This file was created with tikzplotlib v0.9.15.
\begin{tikzpicture}

\begin{groupplot}[group style={group size=2 by 1, horizontal sep=0.3cm}]
\nextgroupplot[
colorbar horizontal,
colorbar style={xtick={0,0.5,1},xticklabels={\(\displaystyle {0.0}\),\(\displaystyle {0.5}\),\(\displaystyle {1.0}\)},minor xtick={}},
colormap/viridis,
colorbar/width=2mm,
height=\figureheight,
minor xtick={},
minor ytick={},
point meta max=1,
point meta min=0,
tick pos=left,
title={data \(\displaystyle q(\boldsymbol{x},t)\)},
width=\figurewidth,
xmin=0, xmax=129,
xtick=\empty,
ymin=0, ymax=129,
ytick=\empty
]
\addplot graphics [includegraphics cmd=\pgfimage,xmin=0, xmax=129, ymin=0, ymax=129] {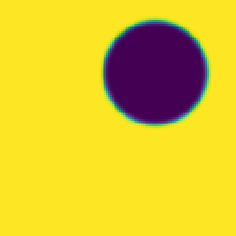};

\nextgroupplot[
colorbar horizontal,
colorbar style={xtick={0.5,1},xticklabels={\(\displaystyle {0.5}\),\(\displaystyle {1.0}\)},minor xtick={}},
colormap/viridis,
colorbar/width=2mm,
height=\figureheight,
minor xtick={},
minor ytick={},
point meta max=1.21717631816864,
point meta min=0.0989818349480629,
tick pos=left,
title={POD \(\displaystyle \tilde{q}(\boldsymbol{x},t)\)},
width=\figurewidth,
xmin=0, xmax=129,
xtick=\empty,
ymin=0, ymax=129,
ytick=\empty
]
\addplot graphics [includegraphics cmd=\pgfimage,xmin=0, xmax=129, ymin=0, ymax=129] {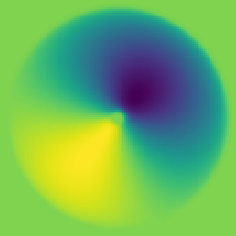};
\end{groupplot}
\end{tikzpicture}
				\caption{Data at time $t=0.11$ and its approximation with the Proper Orthogonal Decomposition (POD) using $r=3$ modes.}
	  			\label{fig:mov-disc-POD}
	\end{figure}

\begin{figure}[htp!]
			\centering
				\begin{subfigure}[t]{0.7\textwidth}
					\centering
				\includegraphics[width=1\linewidth,trim={4cm 1.3cm 3cm 1.5cm},clip]{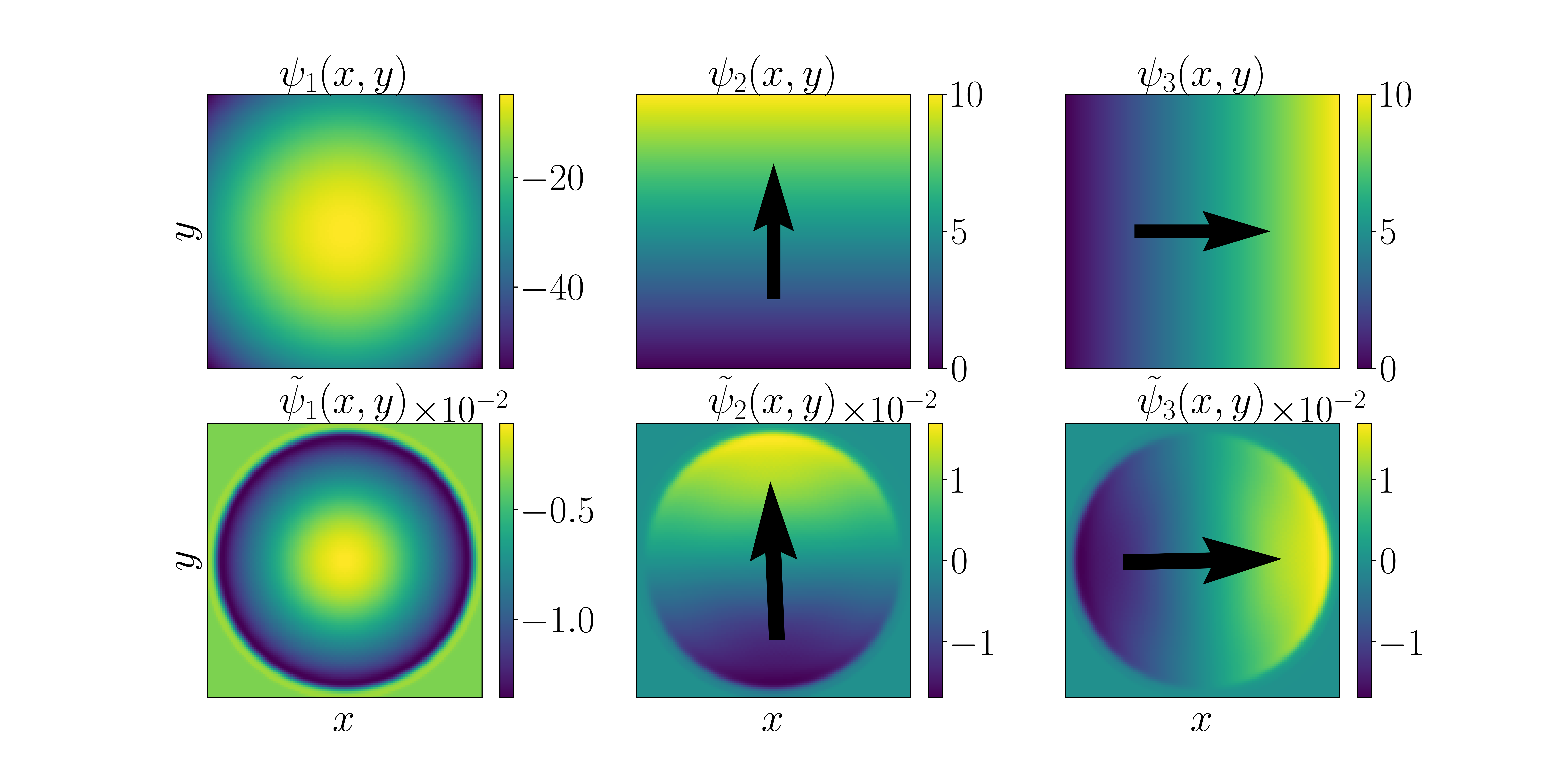}%
							\caption{spatial modes}
							\label{fig:transportfield-modes}
				\end{subfigure}%
				\begin{subfigure}[t]{0.3\textwidth}
					\centering
						\setlength\figureheight{0.9\linewidth}%
						\setlength\figurewidth{0.9\linewidth}%
									 % This file was created with tikzplotlib v0.9.17.
\begin{tikzpicture}

\definecolor{color0}{rgb}{0.12156862745098,0.466666666666667,0.705882352941177}
\definecolor{color1}{rgb}{1,0.498039215686275,0.0549019607843137}
\definecolor{color2}{rgb}{0.172549019607843,0.627450980392157,0.172549019607843}

\begin{axis}[
height=\figureheight,
legend cell align={left},
legend style={
  fill opacity=0.8,
  draw opacity=1,
  text opacity=1,
  at={(0.97,1.03)},
  anchor=south east,
  draw=white!80!black
},
tick align=outside,
tick pos=left,
width=\figurewidth,
x grid style={white!69.0196078431373!black},
xlabel={time \(\displaystyle t\)},
xmin=-0.039, xmax=1.039,
xtick style={color=black},
xtick={-0.5,0,0.5,1,1.5},
xticklabels={
  \(\displaystyle {\ensuremath{-}0.5}\),
  \(\displaystyle {0.0}\),
  \(\displaystyle {0.5}\),
  \(\displaystyle {1.0}\),
  \(\displaystyle {1.5}\)
},
y grid style={white!69.0196078431373!black},
ylabel={\(\displaystyle \tilde{a}_i(t)\)},
ymin=-306.796581368581, ymax=261.121549643845,
ytick style={color=black},
ytick={-400,-200,0,200,400},
yticklabels={
  \(\displaystyle {\ensuremath{-}400}\),
  \(\displaystyle {\ensuremath{-}200}\),
  \(\displaystyle {0}\),
  \(\displaystyle {200}\),
  \(\displaystyle {400}\)
}
]
\addplot [semithick, color0]
table {%
0.01 -280.975240315565
0.03 -280.936442874551
0.05 -280.904046561962
0.07 -280.893569302584
0.09 -280.888517074398
0.11 -280.889839938761
0.13 -280.888355194699
0.15 -280.888926389136
0.17 -280.891098201216
0.19 -280.897046990654
0.21 -280.917625310647
0.23 -280.958121847229
0.25 -280.982120868016
0.27 -280.958121847155
0.29 -280.917625310502
0.31 -280.897046990443
0.33 -280.891098200998
0.35 -280.88892638886
0.37 -280.88835519442
0.39 -280.889839938476
0.41 -280.888517074148
0.43 -280.893569302383
0.45 -280.904046561804
0.47 -280.936442874459
0.49 -280.975240315576
0.51 -280.975240315609
0.53 -280.93644287456
0.55 -280.904046561965
0.57 -280.893569302581
0.59 -280.88851707438
0.61 -280.889839938723
0.63 -280.888355194671
0.65 -280.888926389088
0.67 -280.891098201154
0.69 -280.897046990554
0.71 -280.917625310559
0.73 -280.958121847145
0.75 -280.982120867944
0.77 -280.958121847062
0.79 -280.917625310412
0.81 -280.89704699034
0.83 -280.891098200869
0.85 -280.888926388762
0.87 -280.888355194316
0.89 -280.8898399384
0.91 -280.888517074069
0.93 -280.893569302296
0.95 -280.904046561758
0.97 -280.93644287441
0.99 -280.975240315534
};
\addlegendentry{$i=1$}
\addplot [semithick, color1]
table {%
0.01 -234.837496474915
0.03 -231.102993375539
0.05 -223.730335180678
0.07 -212.847341777692
0.09 -198.611847908445
0.11 -181.248566285802
0.13 -161.024721991547
0.15 -138.262342230182
0.17 -113.319788599688
0.19 -86.5909946448963
0.21 -58.4992289215098
0.23 -29.4851242186596
0.25 -0.000208005251156441
0.27 29.4847115228499
0.29 58.4988260706268
0.31 86.5906079618368
0.33 113.319424161984
0.35 138.262005780372
0.37 161.024418834074
0.39 181.248301200073
0.41 198.611625077944
0.43 212.847164712877
0.45 223.730206672519
0.47 231.102915445578
0.49 234.837470358358
0.51 234.837496474912
0.53 231.102993375515
0.55 223.730335180659
0.57 212.847341777677
0.59 198.611847908427
0.61 181.248566285778
0.63 161.024721991535
0.65 138.262342230165
0.67 113.319788599671
0.69 86.5909946448683
0.71 58.4992289214875
0.73 29.485124218641
0.75 0.000208005232930895
0.77 -29.4847115228668
0.79 -58.4988260706445
0.81 -86.5906079618492
0.83 -113.319424161984
0.85 -138.262005780386
0.87 -161.024418834081
0.89 -181.248301200092
0.91 -198.611625077956
0.93 -212.847164712878
0.95 -223.730206672545
0.97 -231.102915445594
0.99 -234.837470358371
};
\addlegendentry{$i=2$}
\addplot [semithick, color2]
table {%
0.01 -14.7720328555589
0.03 -44.0791068369317
0.05 -72.6875348055381
0.07 -100.152535212683
0.09 -126.038884354048
0.11 -149.939627359829
0.13 -171.474091003592
0.15 -190.305188090131
0.17 -206.135959348975
0.19 -218.718576525058
0.21 -227.863581324871
0.23 -233.432178776181
0.25 -235.307089143228
0.27 -233.432230903863
0.29 -227.863684747876
0.31 -218.718729612748
0.33 -206.136159691739
0.35 -190.305432529961
0.37 -171.474375686099
0.39 -149.939947797017
0.41 -126.039235488624
0.43 -100.152911514892
0.45 -72.6879303483397
0.47 -44.0795154142443
0.49 -14.7724480353247
0.51 14.7720328556044
0.53 44.0791068369746
0.55 72.6875348055839
0.57 100.152535212729
0.59 126.038884354092
0.61 149.939627359868
0.63 171.474091003642
0.65 190.305188090176
0.67 206.135959349014
0.69 218.718576525075
0.71 227.863581324902
0.73 233.43217877622
0.75 235.307089143281
0.77 233.432230903902
0.79 227.863684747918
0.81 218.718729612781
0.83 206.136159691753
0.85 190.305432529996
0.87 171.474375686127
0.89 149.939947797062
0.91 126.039235488666
0.93 100.152911514928
0.95 72.6879303483866
0.97 44.0795154142885
0.99 14.7724480353691
};
\addlegendentry{$i=3$}
\end{axis}

\end{tikzpicture}
									 \caption{time coefficients}
									 \label{fig:transportfield-time}
				\end{subfigure}
				\caption{Visualization of the FTR transport field $\phi(x,t)\approx\sum_{i=0}^3 a_i(t) \psi_i(x,y)$ for the disk moving in a circle (see \cref{eq:movDisc,eq:movDisc1}). Displayed are the expected spatial modes $\psi_i$ (\cref{fig:transportfield-modes}), their temporal amplitudes $a_i$ (\cref{fig:transportfield-time}) and FTR approximation $\tilde{\psi}_i,\tilde{a}_i,i=1,2,3$. The arrows in \cref{fig:transportfield-modes} indicate the direction of the shift. They are computed from the spatial mean of $\nabla\psi_2(x,y)$, $\nabla\psi_3(x,y)$. The corresponding amplitudes $a_2,a_3$ parameterize the circular movement in time.}
				\label{fig:transportfield}
	\end{figure}
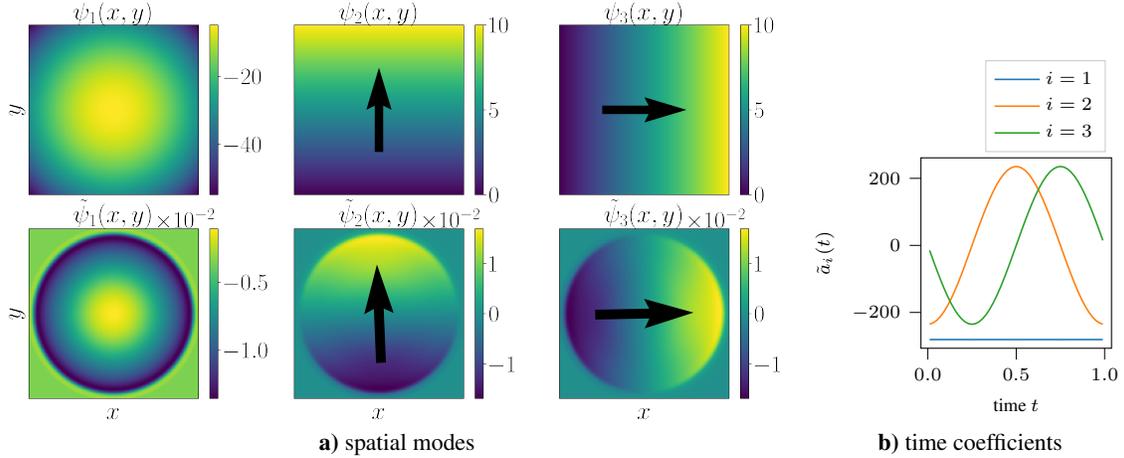

The first synthetic example is taken from \cite{KrahSrokaReiss2020}. It illustrates the idea of the FTR in the pure advection case, without any topological change. The example parameterizes a disk of radius $R=0.15L$, which is moving in a circle:
	\begin{align}
	  q(\vec{x},t)& = f(\phi(\vec{x},t))
	        \quad\text{and}
	\quad \phi(\vec{x},t)=\frac{1}{2R}(\norm{\vec{x}-\vec{x}_0(t)}_2^2-R^2) \label{eq:movDisc1}\\
		& \text{where } \vec{x}_0(t) =L
		\begin{pmatrix}
			0.5 + 1/4 \cos(2 \pi t)\\
			0.5 + 1/4 \sin(2\pi t)
		\end{pmatrix}\,.
		\label{eq:movDisc}
	\end{align}
	The snapshot data $q(\vec{x},t)$ is generate from a level set field $\phi(\vec{x},t)$, which is zero at the outer radius of the disk, i.e. location of the front. The front is generated with help of the function $f(x) = (\tanh(x/\lambda)+1)/2,\, \lambda = 0.1$. One representative snapshot of the data is shown together with its approximation using the POD in \cref{fig:mov-disc-POD}. 
	As the authors of \cite{KrahSrokaReiss2020} have already pointed out, for this example $\phi(\vec{x},t)=\sum_{k=1}^3 \psi_k(\vec{x})a_k(t)$ can be parameterized by only three functions and is therefore of low rank, even if the field $q(\vec{x},t)$ is not. The basis functions $\psi_1,\psi_2,\,\psi_3$ are shown in the top row of \cref{fig:transportfield-modes}. They can be interpreted as a quadratic basis function $\psi_1(x,y)=(x-0.5L)^2+(y-0.5L)^2+R^2+L^2/4$ that represents the initial shape of the contour line with constant time amplitude $a_1(t)=a_1\in\mathbb{R}$, and the linear transport functions $\psi_2(x,y)=x,\psi_3(x,y)=y$ for the shift in $x$/$y$-direction with $a_2(t) \sim \cos(2\pi t),a_3(t) \sim \sin(2\pi t)$.  Note, that the arrows in \cref{fig:transportfield-modes} indicate $\nabla\psi_2(x,y)$, $\nabla\psi_3(x,y)$ the direction of the shift. 
	
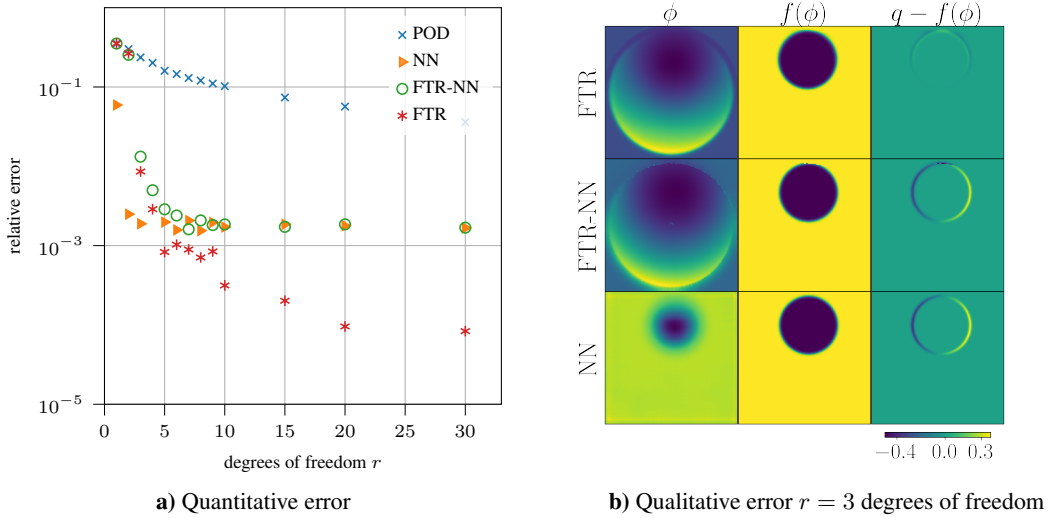
\begin{figure}[htp!]
	      \centering
				\begin{subfigure}[t]{0.5\textwidth}
					\centering
							 \setlength\figureheight{0.9\linewidth}
							 \setlength\figurewidth{0.9\linewidth}
			  % This file was created with tikzplotlib v0.9.17.
\begin{tikzpicture}

\definecolor{color0}{rgb}{0.12156862745098,0.466666666666667,0.705882352941177}
\definecolor{color1}{rgb}{1,0.498039215686275,0.0549019607843137}
\definecolor{color2}{rgb}{0.172549019607843,0.627450980392157,0.172549019607843}
\definecolor{color3}{rgb}{0.83921568627451,0.152941176470588,0.156862745098039}

\begin{axis}[
height=\figureheight,
legend cell align={left},
legend style={fill opacity=0.8, draw opacity=1, text opacity=1, draw=none},
log basis y={10},
tick align=outside,
tick pos=left,
width=\figurewidth,
x grid style={white!69.0196078431373!black},
xlabel={degrees of freedom \(\displaystyle r\)},
xmajorgrids,
xmin=0, xmax=33,
xminorgrids,
xtick style={color=black},
xtick={0,5,10,15,20,25,30,35,40},
y grid style={white!69.0196078431373!black},
ylabel={relative error},
ymajorgrids,
ymin=1e-05, ymax=1,
yminorgrids,
ymode=log,
ytick style={color=black},
ytick={1e-07,1e-05,0.001,0.1,10,1000},
yticklabels={
  \(\displaystyle {10^{-7}}\),
  \(\displaystyle {10^{-5}}\),
  \(\displaystyle {10^{-3}}\),
  \(\displaystyle {10^{-1}}\),
  \(\displaystyle {10^{1}}\),
  \(\displaystyle {10^{3}}\)
}
]
\addplot [semithick, color0, mark=x,  mark options={solid}, only marks]
table {%
1 0.35294216245474
2 0.300317731687437
3 0.236248407176739
4 0.201673335860509
5 0.159784102894094
6 0.145569784672175
7 0.129808184955096
8 0.120579092728741
9 0.11058242315232
10 0.102379531451645
15 0.0737289605395835
20 0.0566677474829191
30 0.0361078041322912
};
\addlegendentry{POD}
\addplot [semithick, color1, mark=triangle*,  mark options={solid,rotate=270}, only marks]
table {%
1 0.0592303201556206
2 0.00248891534283757
3 0.00188650877680629
4 100
5 0.00197975640185177
6 0.00157511851284653
7 0.00207153242081404
8 0.00154433329589665
9 0.0019627925939858
10 0.00171413295902312
15 0.00185864674858749
20 0.00178126059472561
30 0.00167297967709601
};
\addlegendentry{NN}
\addplot [semithick, color2, mark=o,  mark options={solid,fill opacity=0}, only marks]
table {%
1 0.352942109107971
2 0.25429704785347
3 0.0132349142804742
4 0.00499370275065303
5 0.00286845955997705
6 0.00239710952155292
7 0.00160562468226999
8 0.00207886449061334
9 0.00181597576010972
10 0.00185103830881417
15 0.00172364071477205
20 0.00185600691474974
30 0.00168378627859056
};
\addlegendentry{FTR-NN}
\addplot [semithick, color3, mark=asterisk,  mark options={solid}, only marks]
table {%
1 0.352942599676485
2 0.26503815331951
3 0.00859242614360061
4 0.00287922328732646
5 0.000831141118207966
6 0.00102821970879097
7 0.000893421107770546
8 0.000708422918700336
9 0.000844104235225542
10 0.000314822058078974
15 0.000201817116890624
20 9.56191307577227e-05
30 8.35009392370084e-05
};
\addlegendentry{FTR}
\end{axis}

\end{tikzpicture}
							 \caption{Quantitative error}
							 \label{fig:mov-disc-compare-rel-error}
				\end{subfigure}%
	      \begin{subfigure}[t]{0.5\textwidth}
	  			\includegraphics[width=0.9\linewidth, trim=00 00 00 20,clip]{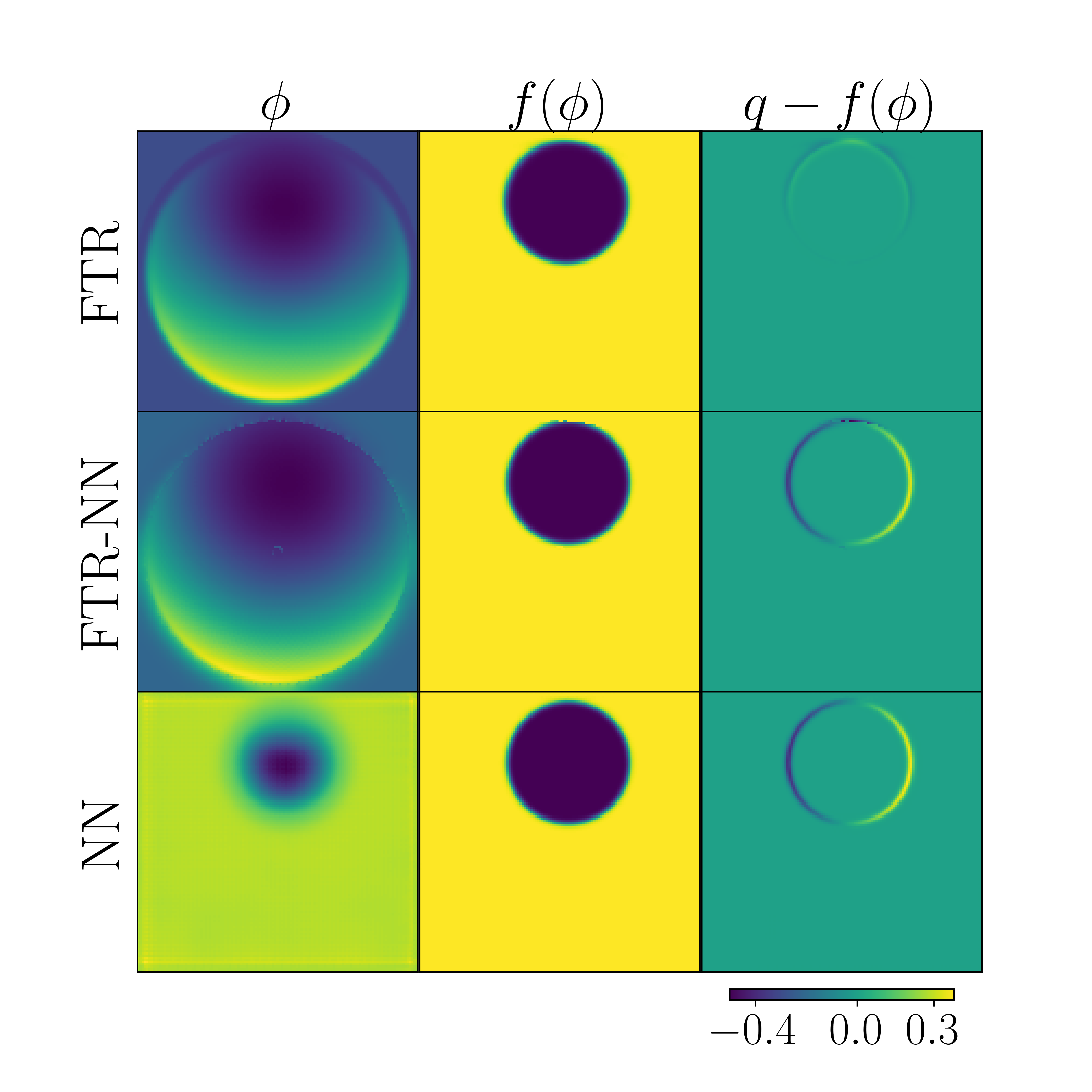}
					\caption{Qualitative error $r=3$ degrees of freedom}
					\label{fig:mov-disc-compare-diff}
	    	\end{subfigure}%
	  			\caption{Comparison of POD, FTR, FTR-NN and the symmetrical autoencoder structure labeled with NN.  \cref{fig:mov-disc-compare-rel-error} compares the relative errors in the Frobenius norm for different degrees of freedom. \cref{fig:mov-disc-compare-diff} visualizes the level-set field $\phi$ together with the approximation of the data $\tilde{q}=f(\phi)$ and the deviation from the exact data $q-\tilde{q}$ for one selected snapshot.}
	  			\label{fig:mov-disc-compare}
	  \end{figure}

	To show that the singular value thresholding \cref{alg:iterFTR} (FTR) and the neural network approach \cref{subsec:ML} (FTR-NN) can find a similar basis set, we generate 200 equally spaced snapshots from $q$ in the time interval $0<t\le 1$, discretized with $129\times129$ grid points in the rectangular domain $[0,L]^2$. The data was split into a train and test set, where every second sample is a test sample. After training the neural network on the training samples, it is compared to the results of the POD and the thresholding \cref{alg:iterFTR} using the test samples. 
    The results are visualized in \cref{fig:transportfield,fig:mov-disc-compare}.
	In \cref{fig:mov-disc-compare} we compare the results of both FTR-algorithms (FTR, FTR-NN)
  and a simple symmetrical autoencoder structure, labeled with NN (for details see \cref{appx:AutoencoderArchitecture}). The NN decoder attempts to implement the encoder in an inverse manner (see details \cref{appx:AutoencoderArchitecture}), which is a common practice in dimension reduction.
   The relative errors in the Frobenius norm are
	shown in \cref{fig:mov-disc-compare-rel-error}. The quantitative errors of FTR-NN and the FTR show a significant drop using $r=3$ degrees of freedom (FTR basis functions/latent space dimension), which is in accordance with the proposed level-set field. In contrast, the POD is showing a much slower convergence of the relative error. 
	Comparing the two networks NN and FTR-NN regarding the quantitative errors shows that the additional depth of the NN-decoder compared to the one layer decoder in FTR-NN does not influence the minimal relative error. This leads us to conclude that additional depth is not needed for a better representation. However, note that the NN needs fewer degrees of freedom to converge to its minimal relative error, which is due to the higher expressivity of a deeper network. 
    Furthermore, it is important to note that the FTR-thresholding algorithm outperforms both networks, when increasing the number of degrees of freedom, for this special example.   
	For qualitative comparison, \cref{fig:mov-disc-compare-diff} shows the approximation of one snapshot before and after activation (first and second column), together with the difference in the last column. Comparing the POD in \cref{fig:mov-disc-POD} to the FTR results shows, that the typical stair casing behavior (which becomes a blurring of the sharp structures for many snapshots as used here) of the POD can be overcome with the FTR ansatz that recaptures the sharp front.
	We observe that both qualitative and quantitative errors of the FTR-NN and iterative thresholding approach yield similar results. In this study, we use $\smoothness=10^{-7}$ for regularizing the smoothness of $\phi$ at the output of the FTR-NN and NN decoder. As visualized in \cref{appx:AutoencoderArchitecture}  \cref{fig-apx:disc-smoothness}, for larger smoothness parameter $\smoothness>10^{-7}$ the transport field of the FTR-NN is smoothly continued at areas of no information (no transport), but the additional constraint \cref{eq:NN-reg} can cause a larger overall approximation error. %\footnote{The regularization was necessary for convergence of the autoencoder structure.}. 
	However, the level-set fields of the iterative thresholding approach and NN are almost identical inside areas where fronts have been transported.
	This is due to the special choice of the encoder.  

	 \Cref{fig:transportfield} compares the fields $\psi_1,\psi_2,\,\psi_3$, obtained by contemplation, to the first three modes of $\phi$.
	Similar to the proposed functions, the auxiliary field can be split into a mode ($\tilde{\psi}_1$) responsible for the shape of the disk and two modes that parameterize the transport ($\tilde{\psi}_2,\tilde{\psi}_3$).
	As expected for this special case, $\tilde{a}_1$ is constant and
	$\tilde{a}_2,\tilde{a}_3\sim \cos(2\pi t + \theta)$, with $\theta \in \mathbb{R}$ depends on the alignment (indicated as arrows in \cref{fig:transportfield-modes}) of the two shifting functions $\tilde{\psi}_2,\tilde{\psi}_3$. 
	The modes $\tilde{\psi}_2,\tilde{\psi}_3$ only have meaningful values along the trajectories of the front because the algorithm can not in-paint $\phi$ in areas of no transport. This explains that the modes in \cref{fig:transportfield} are zero outside the circle.

%%%%%%%%%%%%%%%%%%%%%%%%%%%%%%%%%%%%%%%%%%%%%%%%%%%%%%%%%%%%%%%%%%%%%%%%%%%%%%%%%%%%%%
% Topo Change
%%%%%%%%%%%%%%%%%%%%%%%%%%%%%%%%%%%%%%%%%%%%%%%%%%%%%%%%%%%%%%%%%%%%%%%%%%%%%%%%%%%%%%
\subsubsection{Advection with Topology Change}
\label{subsubsec:advection_with_topo_change}

\begin{figure}[htp!]
    \centering
%\vspace{-0.2cm}
\begin{subfigure}[t]{0.4\textwidth}
	\centering
		\includegraphics[width=1\linewidth]{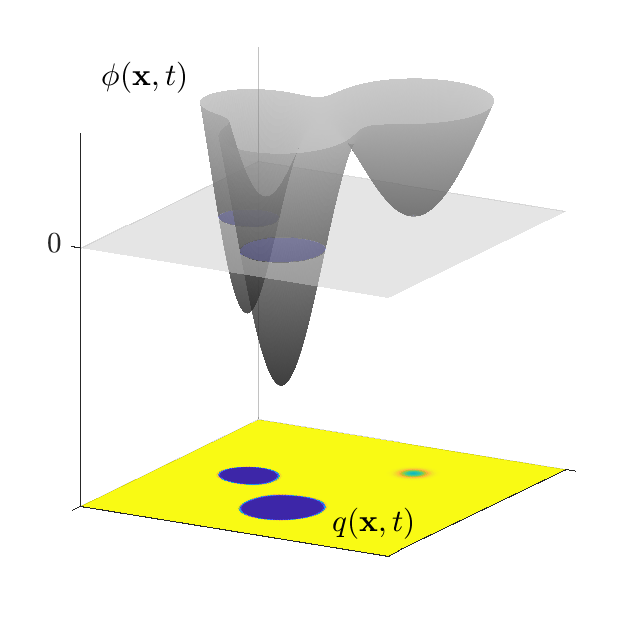}%
	    \caption{data $t=0$}
	    \label{subfig:lvl_topo0}
	  \end{subfigure}
		\begin{subfigure}[t]{0.4\textwidth}
			\centering
		\includegraphics[width=1\linewidth]{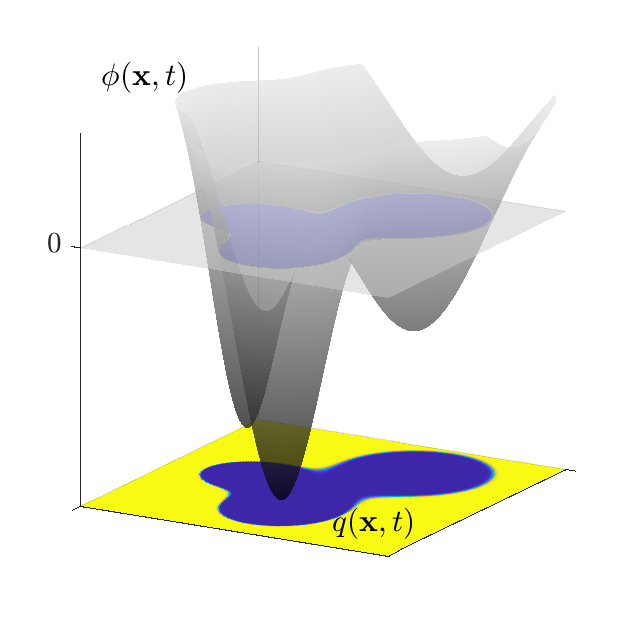}
    \caption{data $t=0.4$}
    \label{subfig:lvl_topo40}
  \end{subfigure}
	\begin{subfigure}[t]{0.4\textwidth}
		\centering
	\includegraphics[width=1\linewidth]{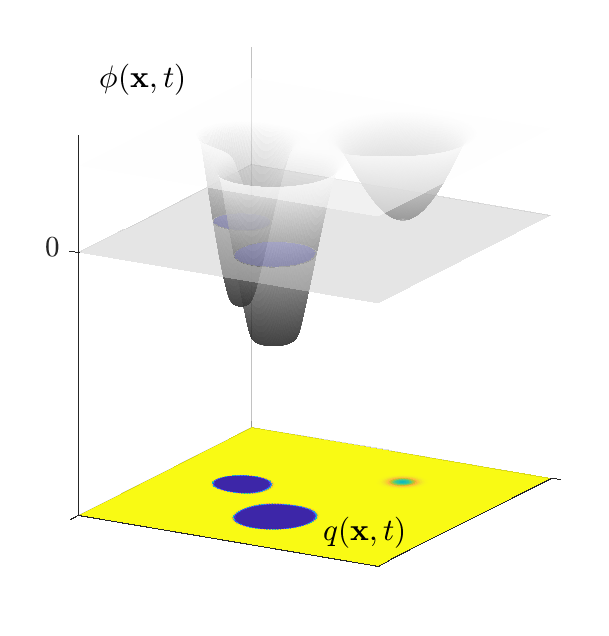}
	\caption{FTR $t=0$, rank $r=2$}
	\label{subfig:lvl_topo0_ftr}
\end{subfigure}
	\begin{subfigure}[t]{0.4\textwidth}
		\centering
	\includegraphics[width=1\linewidth]{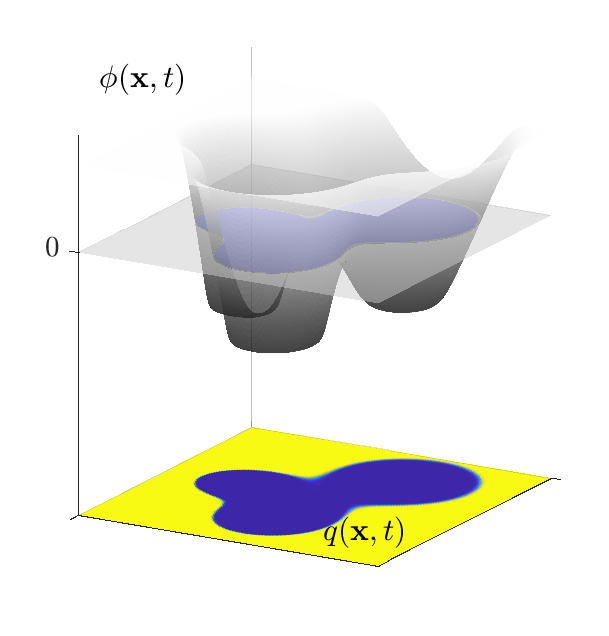}
	\caption{FTR $t=0.4$,  rank $r=2$}
	\label{subfig:lvl_topo40_ftr}
\end{subfigure}
  \caption{Graph of the auxiliary field $\phi$ (\cref{eq:levl_topo} in a)-b) and its FTR approximation in c)-d). The resulting snapshots $q=f(\phi)$ are shown as color plot in the $xy$ plane. The intersection with the zero level is visualized.}
  \label{fig:adv_levl_topo}
\end{figure}

In this example we show that our approach is capable of handling
transport with topological changes. Therefore, we introduce the synthetic snapshot data $q(\vec{x},t)=f(\phi(\vec{x},t))$ build from the
level-set field
\begin{equation}
	\label{eq:levl_topo}
	\varphi(\vec{x},t) = \sum_{k=1}^3 -A_k e^{-\sigma_k r_k}-t,\quad r_k = \norm{\vec{x}-\vec{x}_i}_2\,,
\end{equation}
which we try to approximate. The front $f$ is chosen as above. The level-set field is sampled equidistantly using $256\times 265$ grid points
in $[0,10]^2$, with $(A_1,A_2,A_3)=(1,1.4,1.2)$, $(\sigma_1,\sigma_2,\sigma_3)=(0.1,0.3,0.5)$ and $\vec{x}_1=(7.5,3.5),\vec{x}_2=(2.5,5.0),\, \vec{x}_3=(5.0,7.6)$.
Furthermore, 101 equally spaced snapshots with $0\le t\le 0.5$ are constructed from \cref{eq:levl_topo}.
As above, we split the samples in a test and train set, where every second sample is used for testing the autoencoders. After training the networks, they are compared to the reconstruction errors of the POD and FTR using the test samples. The level-set fields for $t=0$ and $t=0.4$  are visualized as a surface plot in \cref{fig:adv_levl_topo}, together with the resulting snapshots of $q$ as a color plot.
%%%%%%%%%%%%%%%%%%%%%%%%%%%%%%%%%%%%%%%%
% figure
%%%%%%%%%%%%%%%%%%%%%%%%%%%%%%%%%%%%%%%%
%%%%%%%%%%%%%%%%%%%%%%%%%%%%%%%%%%%%%%%%
%%%%%%%%%%%%%%%%%%%%%%%%%%%%%%%%%%%%%%%%
 The intersection of $\phi$ with the zero plane parameterizes the surface of the front. For increasing $t$, the level-set function is shifted vertically and produces an expanding surface of the front, which is merged from three independent into one single front contour. The merging of the fronts allows no smooth bijective mapping between the contour lines of the front at time $t=0$ to $t=0.4$. This property makes it difficult for most dimension reduction methods, which can handle transports because these rely on one-to-one mappings between different time or parameter instances.
 
\begin{figure}[htp!]
		\centering
	\setlength\figureheight{0.45\linewidth}
		\setlength\figurewidth{0.45\linewidth}
			% This file was created with tikzplotlib v0.9.17.
\begin{tikzpicture}

\definecolor{color0}{rgb}{0.12156862745098,0.466666666666667,0.705882352941177}
\definecolor{color1}{rgb}{1,0.498039215686275,0.0549019607843137}
\definecolor{color2}{rgb}{0.172549019607843,0.627450980392157,0.172549019607843}
\definecolor{color3}{rgb}{0.83921568627451,0.152941176470588,0.156862745098039}

\begin{axis}[
height=\figureheight,
legend cell align={left},
legend style={fill opacity=0.8, draw opacity=1, text opacity=1, draw=white!80!black},
log basis y={10},
tick align=outside,
tick pos=left,
width=\figurewidth,
x grid style={white!69.0196078431373!black},
xlabel={degrees of freedom \(\displaystyle r\)},
xmin=0, xmax=9.5,
xtick style={color=black},
xtick={0,1,2,3,4,5,6,7,8,9},
xticklabels={
  \(\displaystyle {0}\),
  \(\displaystyle {1}\),
  \(\displaystyle {2}\),
  \(\displaystyle {3}\),
  \(\displaystyle {4}\),
  \(\displaystyle {5}\),
  \(\displaystyle {6}\),
  \(\displaystyle {7}\),
  \(\displaystyle {8}\),
  \(\displaystyle {9}\)
},
y grid style={white!69.0196078431373!black},
ylabel={relativ error},
ymin=9.75642319088971e-05, ymax=0.329476887594134,
ymode=log,
ytick style={color=black},
ytick={1e-06,1e-05,0.0001,0.001,0.01,0.1,1,10},
yticklabels={
  \(\displaystyle {10^{-6}}\),
  \(\displaystyle {10^{-5}}\),
  \(\displaystyle {10^{-4}}\),
  \(\displaystyle {10^{-3}}\),
  \(\displaystyle {10^{-2}}\),
  \(\displaystyle {10^{-1}}\),
  \(\displaystyle {10^{0}}\),
  \(\displaystyle {10^{1}}\)
}
]
\addplot [semithick, color0, mark=x, mark options={solid}, only marks]
table {%
1 0.227315157651901
2 0.0972356274724007
3 0.047885786741972
4 0.024815620854497
5 0.012823055498302
6 0.00670261355116963
7 0.0035079806111753
8 0.00184363266453147
9 0.000972695706877857
};
\addlegendentry{POD}
\addplot [semithick, color1, mark=triangle*,  mark options={solid,rotate=270}, only marks]
table {%
1 0.00191282760351896
2 0.00314735434949398
3 0.00209194677881896
4 0.00276960548944771
5 0.00227608112618327
6 0.00219403998926282
7 0.00221412861719728
8 0.00177926721516997
9 0.00211092340759933
};
\addlegendentry{NN}
\addplot [semithick, color2, mark=o, mark options={solid,fill opacity=0}, only marks]
table {%
1 0.202121630311012
2 0.00167121866252273
3 0.00201942399144173
4 0.00119881948921829
5 0.00144371145870537
6 0.000679648364894092
7 0.00171721505466849
8 0.00132444023620337
9 0.00128549023065716
};
\addlegendentry{FTR-NN}
\addplot [semithick, color3, mark=asterisk, mark options={solid}, only marks]
table {%
1 0.22773880155583
2 0.000970599598433606
3 0.000345817797625974
4 0.000201052933032172
5 0.000155342124547687
6 0.000145920528549174
7 0.0001427262279033
8 0.000141554929470171
9 0.000141149243125245
};
\addlegendentry{FTR}
\end{axis}

\end{tikzpicture}
			\caption{Comparison of the relative errors for the advection example with topology change using the proper orthogonal decomposition (POD), the FTR iterative thresholding algorithm (FTR), the FTR autoencoder structure (FTR-NN) and a standard autoencoder (NN).}
			\label{fig:adv_lvl_topo_error}
\end{figure}
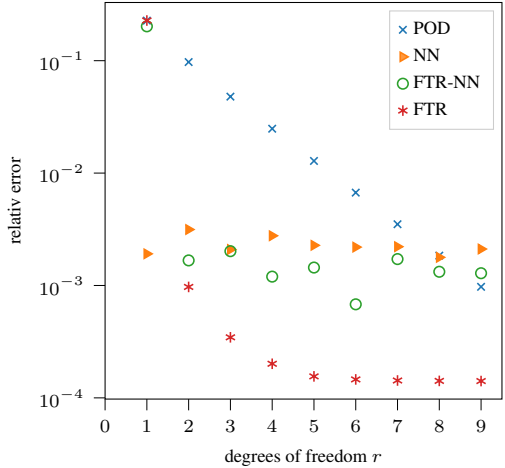

As presented in \cref{fig:adv_lvl_topo_error}, the FTR approximates the dynamics within two dyadic pairs with an error smaller than 0.2$\%$, which is expected from the two-term dyadic structure in \cref{eq:levl_topo}. The networks approximation errors behave as in the case of the moving disk. The FTR-NN gives similar results as the FTR, but with larger minimal relative error. Due to the additional depth, the NN only needs one degree of freedom to converge towards its minimal error.  
Topology changes of the zero level-set are nicely recovered as is illustrated in
\cref{subfig:lvl_topo0_ftr,subfig:lvl_topo40_ftr}, since the FTR approach can recover the initial auxiliary field $\phi$ in the regions of transport.

\subsection{Application to Advection-reaction-diffusion Systems}
\label{subsec:offline-ARD}
To motivate the FTR approach for more complex examples,
we introduce the advection-diffusion-reaction PDE with a KKP reaction term
		\begin{equation}
		    \label{eq:react-diff-system}
			\begin{cases}
				\partial_t q &= - \vec{u}\cdot \vec{\nabla} q + \diffconst\Delta q - \reactconst q^2(q-1)\\
				q(\vec{x},0) &= q_0(\vec{x})
			\end{cases}\,,
		\end{equation}
		on a square, two dimensional domain $\Omega=[0,L]^2$ with periodic boundary conditions and time interval $[0,T]$.
		The PDE is discretized in space using 6th order central finite differences, and in time with an explicit Runge-Kutta method of 5th(4th) order \cite{DormandPrince1980}. In the following, we refer to the discretized system as the full order model (FOM). All simulation parameters are listed in \cref{tabl:advec-diff-react}.
		\begin{table}
			\centering
			\begin{tabular}{l l c}
				\toprule
				Name & Value  \\
				\midrule
				\textbf{FOM - parameters}\\
				Simulation time & $T=3$ \\
				Grid resolution & $M=512 \times 512$\\
				Domain size &  $L=1$\\
				Diffusion constant & $\diffconst=10^{-3}$\\
				Reaction constant & $\reactconst=10$ \\
				Advection of vortex-pair &$c=10$\\
				\midrule
				\textbf{ROM - parameters}\\
				Number of snapshots &$\Ntime=100$ \\
				Front function & $f(x) = \sigmoid (x)$\\
				\bottomrule
			\end{tabular}
			\caption{Parameters of the 2D ARD simulation of \cref{subsec:offline-ARD}}
			\label{tabl:advec-diff-react}
		\end{table}

\begin{figure}[htp!]
			\centering
					\includegraphics[width=0.4\linewidth,trim={0 0 3cm 0},clip]{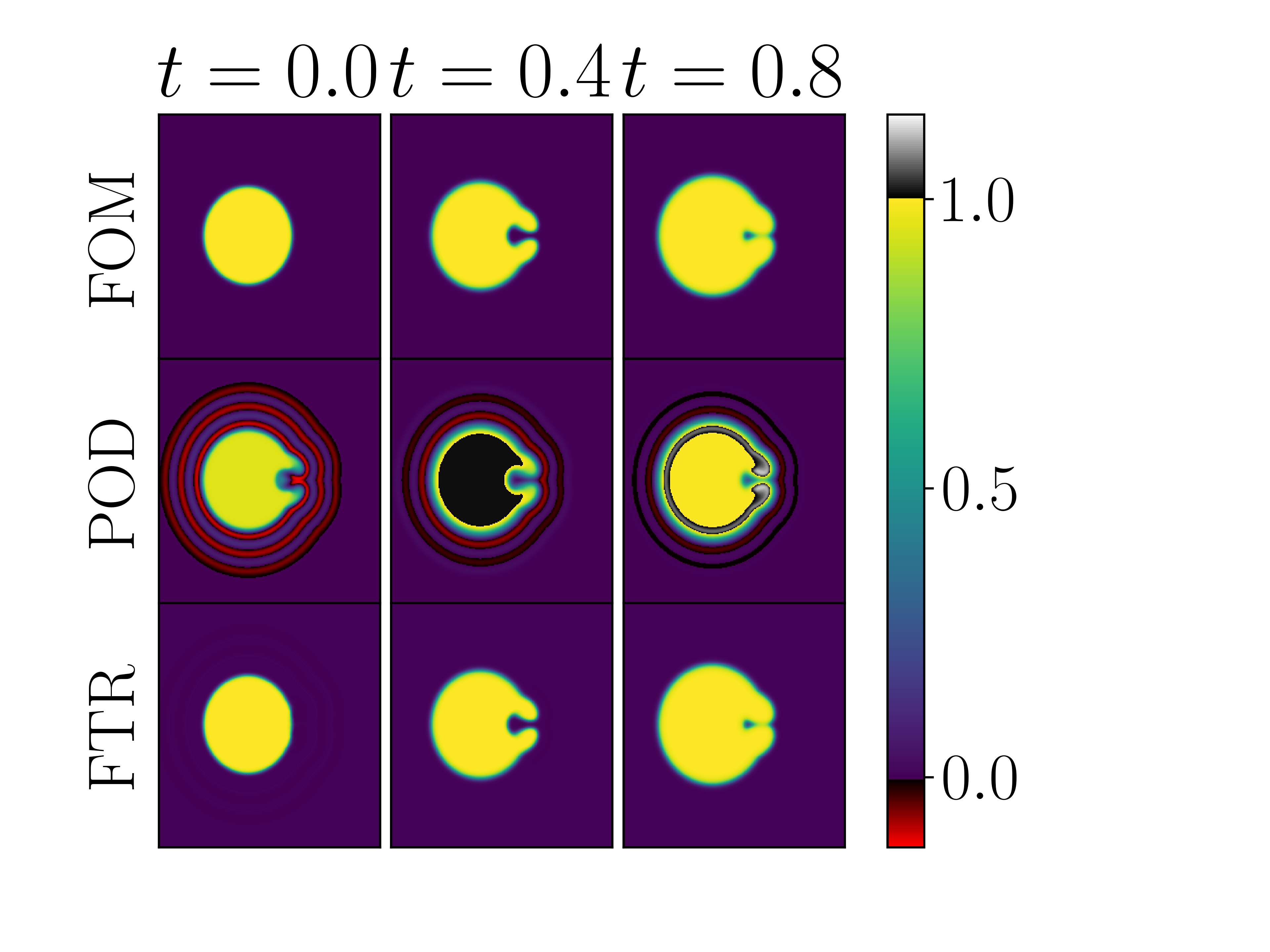}
		  \caption{Qualitative comparison of the reconstruction errors of the 2D ARD system \cref{eq:react-diff-system} at three different time instances $t=0.0, 0.4, 0.8$ (respectively left, middle, right column). The plot shows the FOM data (top row) and its reconstructions using the POD (middle row) and FTR (bottom row). For the FTR and POD, $r=6$ degrees of freedom are used. The colorbar is chosen such that values outside the initial range of values $0\le q \le 1$ are highlighted in black or red.}
		  \label{fig:pacman-offline}
		\end{figure}

		For our test case we choose a velocity field inspired by the vortex pair example in \cite{Reiss2015}. Therefore $\vec{u}=\vec{\nabla} \times\omega$ is expressed in terms of the vorticity
		\begin{equation}
			\label{eq:vortex-pair}
			\omega(\vec{x},t) = \omega_0 e^{-t^2 / \tau^2} (e^{-r^2_1(t)/r^2_0}+e^{-r^2_2(t)/r^2_0})\,  \quad r_i(t) = \norm{\vec{x}-\vec{x}_i(t)}_2\,,
		\end{equation}
		which parameterizes a moving vortex pair $\vec{x}_1=L(0.6-ct,0.49)$, $\vec{x}_2=L(0.6-ct,0.51)$, $r_0=5\times10^{-4}$ with an initial amplitude $\omega_0=10^{3}$ decaying slowly in time (decay constant $\tau=3T$).
		The initial distribution of the reactant $q$ is given by:
		\begin{equation}
			q_0(x,y)=
						\begin{cases}
							1 & \sqrt{(x-0.4L)^2+(y-0.5L)^2}>0.2L\,,\\
							0 &\text{else}\,.
						\end{cases}
		\end{equation}
    The velocity field and initial distribution are tuned to mimic a flame kernel interacting with a vortex pair, which is a usual phenomenon in turbulence flame inter-
actions. During the simulation, the synthetic vortex pair \cref{eq:vortex-pair} moves towards burning gas and mixes
unburned ($q = 1$) with burned gas ($q=0$), such that a small island of unburned gas detaches
into the burned area, creating a topology change in the contour line of the front. The time evolution of the FOM is visualized for $t=0.0,0.4,0.8$ in the top row of \cref{fig:pacman-offline}. In the second and third row, the FTR and POD are compared using $r=6$ degrees of freedom. The POD approximation shows the typical staircase behavior as oscillations occur before and after the contour line of the front. The oscillations violate the initial range of values $0\le q \le 1$, which is depicted as red and black areas in \cref{fig:pacman-offline}. Therefore, preservation of physical structure cannot be expected.
Here, the FTR approach gives much better results, restricting the approximation on the initial range of values due to the range of the sigmoid function $f(x)\in[0,1]$.

\newcommand{\q}{\vec{q}}
\newcommand{\RHS}{\vec{F}}
\newcommand{\mapping}{f}
\newcommand{\LinOp}{\matr{L}}
\newcommand{\NonLinOp}{\vec{N}}
\newcommand{\params}{\mathcal{P}}
\newcommand{\mat}[1]{#1}
\renewcommand{\Pa}{\matr{P}_{\vec{a}}}
\newcommand{\Jg}{\matr{J}_{g}}

\section{Galerkin and data-driven Models for Moving Fronts}
\label{sec:non-lin-roms-moving-fronts}

In the previous sections we have addressed the so-called \textit{offline stage} of a model order reduction procedure, in which data is collected and its dimension is reduced. The reduced model generated by the FTR algorithm in \cref{subsec:FTR} is nonlinear, which poses additional challenges for the \textit{online stage}, to predict and interpolate new system states. This section is therefore dedicated to online prediction methods.
In \cref{subsec:data-driven} we use a non-intrusive, i.e.~equation free, approach of \cite{LangeBruntonKutz2021} and introduce an intrusive approach, the hyper-reduced Galerkin method in \cref{subsec:manifoldgalerkin} for 1D and 2D ARD systems.

\subsection{Data-driven Methods}
\label{subsec:data-driven}
With the rise of data-driven methods in model order reduction, non-intrusive prediction methods of the reduced system, e.g. POD-DL-ROM \cite{FrescaManzoni2022}, SINDy \cite{FukamiMuratZhangFukagata2021,QuadeAbelKutzBrunton2018} or Fourier-Koopman forecasting \cite{LangeBruntonKutz2021}, have become prominent. Although the methods make specific assumptions on the system at hand, they can be useful, since they allow rapid evaluation of the reduced variables with good accuracy. This is especially beneficial if the reduced space is a nonlinear manifold, which makes any Galerkin-projection approach more complex and costly, as is shown in the next section.

Following the approach of \cite{LangeBruntonKutz2021}, we can derive new system states and extrapolate in time with help of the Fourier-Koopman framework implemented in \cite{Fourier-Koopman-Github}. The Fourier-Koopman framework imposes the assumption that the reduced state $\vec{a}(t)\in\mathbb{R}^r$ is quasi-periodic in $t$ and can be thus parameterized by:
\begin{equation}
\label{eq:fourier-koopman-parametrization}
    \vec{a}(t) = \matr{A}\vec{\Omega}(t) \quad \text{with} \quad \vec{\Omega}(t)=
    \begin{pmatrix}
        \cos(\vec{\omega}t)\\
        \sin(\vec{\omega}t)
    \end{pmatrix}\,.
\end{equation}
Here, $\matr{A}\in\mathbb{R}^{r\times p}$ and $\vec{\omega}\in \mathbb{R}^{p/2}$ are determined by solving the optimization problem:
	\begin{equation}
	    \label{eq:fourier-koopman-opt}
		\min_{\vec{\omega},\matr{A}} \sum_{n=0}^{N-1} \Vert \vec{a}(t_n) - \matr{A}\vec{\Omega}(t_n)\Vert_2^2\,,
	\end{equation}
in a smart way \cite{LangeBruntonKutz2021}.
Since the dynamical system presented in \cref{subsec:MovingDisc} is quasi-periodic, we can apply the method to the FTR decomposition
\begin{equation}
    \vec{q}(t)\approx\tilde{\vec{q}}(t) = f(\FTRmodes\vec{a}(t))
\end{equation}
using the basis functions $\FTRmodes=[\tilde{\vec{\psi}}_1,\tilde{\vec{\psi}}_2,\tilde{\vec{\psi}}_3]$, shown in \cref{fig:transportfield} together with the amplitudes $\vec{a}(t)=(a_1(t),a_2(t),a_3(t))$ at the sampled time points $\{t_n=n\Delta t\mid n=0,\dots,N-1\}$
\footnote{ The systems dynamics can be further reduced by rewriting $f(\FTRmodes\vec{a}(t)) = f(\matr{\tilde{\FTRmodes}}\vec{\tilde{a}}(t)+\vec{b}), \vec{b}\in\mathbb{R}^M, \vec{\tilde{a}}\in\mathbb{R}^{r-1},\matr{\tilde{\FTRmodes}}\in\mathbb{R}^{M\times (r-1)}$. The offset vector $\vec{b}$ then contains the time independent part of the decomposition shown as constant line in \cref{fig:mov-disc:online}. This can be done similarly for the POD.}
. From the sampled data we compute $\matr{A},\vec{\omega}$. The resulting model $\tilde{\vec{q}}(t)=f(\FTRmodes\matr{A}\vec{\Omega}(t))$ is evaluated at $t_{n+1/2}=(n+1/2)\Delta t$ for $n=0,\dots,2N-1$. Similarly, we can derive an approximation with the POD. Both results are compared in \cref{fig:mov-disc:online}.
%%%%%%%%%%%%%%%%%%%%%%%%%%%%%%%%%%%%%%
\definecolor{color0}{rgb}{0.12156862745098,0.466666666666667,0.705882352941177}
\definecolor{color1}{rgb}{1,0.498039215686275,0.0549019607843137}
\definecolor{color2}{rgb}{0.172549019607843,0.627450980392157,0.172549019607843}
\tikzset{cross/.style={cross out, draw=black, minimum size=2*(#1-\pgflinewidth), inner sep=0pt, outer sep=0pt},
%default radius will be 1pt.
cross/.default={1pt}}
\newcommand{\lineA}{\raisebox{0pt}{\tikz{\draw[color0,solid,line width = 1.0pt](2.mm,0) node[cross,color0]{} (3mm,1.5mm);\draw[-,color0,solid,line width = 1.0pt](0.,0.8mm) -- (4mm,0.8mm)}}}

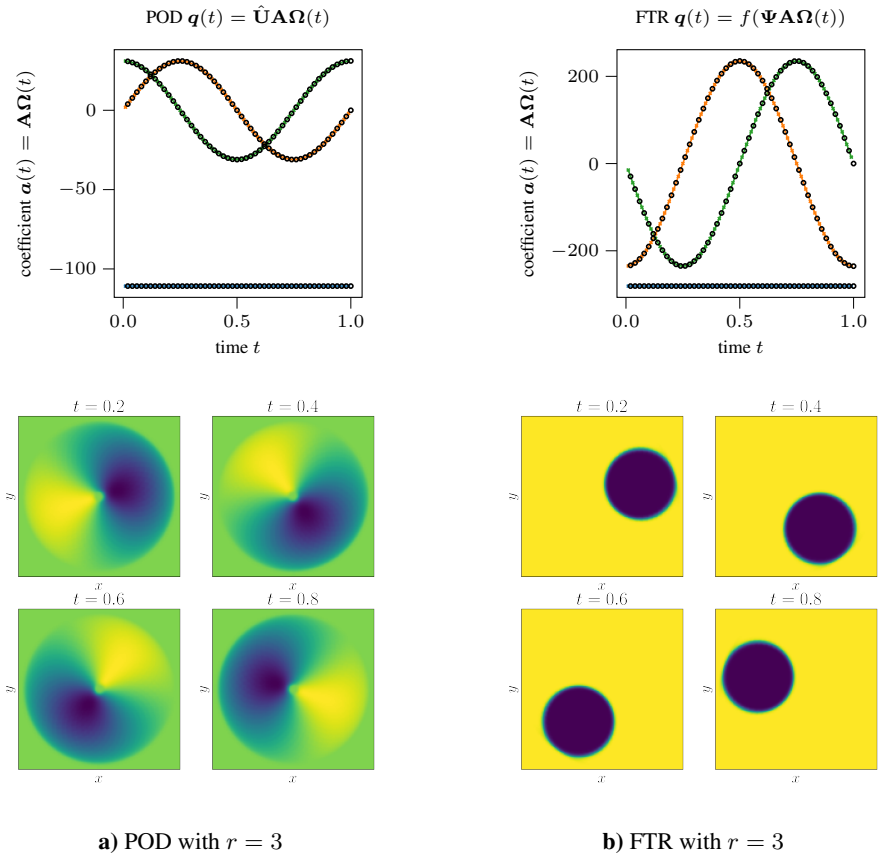
\begin{figure}
    \centering
%\vspace{-0.2cm}
\begin{subfigure}[t]{0.4\textwidth}
	\centering
	\scriptsize
	\begin{center}
    \setlength{\figureheight}{0.8\linewidth}
    \setlength{\figurewidth}{0.8\linewidth}
	% This file was created with tikzplotlib v0.9.15.
\begin{tikzpicture}

\definecolor{color0}{rgb}{0.12156862745098,0.466666666666667,0.705882352941177}
\definecolor{color1}{rgb}{1,0.498039215686275,0.0549019607843137}
\definecolor{color2}{rgb}{0.172549019607843,0.627450980392157,0.172549019607843}

\begin{axis}[
height=\figureheight,
legend cell align={left},
legend style={
  fill opacity=0.8,
  draw opacity=1,
  text opacity=1,
  at={(0.09,0.5)},
  anchor=west,
  draw=white!80!black
},
legend cell align={left},
legend style={fill opacity=0.8, draw opacity=1, text opacity=1, draw=none},
minor xtick={},
minor ytick={},
tick align=outside,
tick pos=left,
title={POD \(\displaystyle \vec{q}(t) = \hat{\mathbf{U}}\mathbf{A}\vec{\Omega}(t)\)},
width=\figurewidth,
x grid style={white!69.0196078431373!black},
xlabel={time \(\displaystyle t\)},
xmin=-0.0395, xmax=1.0495,
xtick style={color=black},
xtick={-0.5,0,0.5,1,1.5},
xticklabels={
  \(\displaystyle {\ensuremath{-}0.5}\),
  \(\displaystyle {0.0}\),
  \(\displaystyle {0.5}\),
  \(\displaystyle {1.0}\),
  \(\displaystyle {1.5}\)
},
y grid style={white!69.0196078431373!black},
ylabel={coefficient \(\displaystyle \vec{a}(t)=\matr{A}\vec{\Omega}(t)\)},
ymin=-117.917031002045, ymax=38.1951663017273,
ytick style={color=black},
ytick={-150,-100,-50,0,50},
yticklabels={
  \(\displaystyle {\ensuremath{-}150}\),
  \(\displaystyle {\ensuremath{-}100}\),
  \(\displaystyle {\ensuremath{-}50}\),
  \(\displaystyle {0}\),
  \(\displaystyle {50}\)
}
]
\addplot [semithick, color0, mark=x, mark size=0.8, mark options={solid}, forget plot]
table {%
0.01 -110.821014404297
0.03 -110.820991516113
0.05 -110.82096862793
0.07 -110.820953369141
0.09 -110.820953369141
0.11 -110.820953369141
0.13 -110.820953369141
0.15 -110.820953369141
0.17 -110.820953369141
0.19 -110.820953369141
0.21 -110.82096862793
0.23 -110.820999145508
0.25 -110.821022033691
0.27 -110.820999145508
0.29 -110.82096862793
0.31 -110.820953369141
0.33 -110.820953369141
0.35 -110.820953369141
0.37 -110.820953369141
0.39 -110.820953369141
0.41 -110.820953369141
0.43 -110.820953369141
0.45 -110.82096862793
0.47 -110.820991516113
0.49 -110.821014404297
0.51 -110.821014404297
0.53 -110.820991516113
0.55 -110.82096862793
0.57 -110.820953369141
0.59 -110.820953369141
0.61 -110.820953369141
0.63 -110.820953369141
0.65 -110.820953369141
0.67 -110.820953369141
0.69 -110.820953369141
0.71 -110.82096862793
0.73 -110.820999145508
0.75 -110.821022033691
0.77 -110.820999145508
0.79 -110.82096862793
0.81 -110.820953369141
0.83 -110.820953369141
0.85 -110.820953369141
0.87 -110.820953369141
0.89 -110.820953369141
0.91 -110.820953369141
0.93 -110.820953369141
0.95 -110.82096862793
0.97 -110.820991516113
0.99 -110.821014404297
};
\label{lineA1}
\addplot [semithick, color1, mark=x, mark size=0.8, mark options={solid}, forget plot]
table {%
0.01 1.95273315906525
0.03 5.82740068435669
0.05 9.61016654968262
0.07 13.2413787841797
0.09 16.6637649536133
0.11 19.823356628418
0.13 22.670316696167
0.15 25.1597518920898
0.17 27.2524013519287
0.19 28.9152698516846
0.21 30.1221199035645
0.23 30.8539352416992
0.25 31.099157333374
0.27 30.8539352416992
0.29 30.1221199035645
0.31 28.9152698516846
0.33 27.2524013519287
0.35 25.1597518920898
0.37 22.670316696167
0.39 19.823356628418
0.41 16.6637649536133
0.43 13.2413787841797
0.45 9.6101655960083
0.47 5.82740020751953
0.49 1.95273268222809
0.51 -1.95273315906525
0.53 -5.82740068435669
0.55 -9.61016654968262
0.57 -13.2413787841797
0.59 -16.6637649536133
0.61 -19.823356628418
0.63 -22.670316696167
0.65 -25.1597518920898
0.67 -27.2524013519287
0.69 -28.9152698516846
0.71 -30.1221199035645
0.73 -30.8539352416992
0.75 -31.099157333374
0.77 -30.8539352416992
0.79 -30.1221199035645
0.81 -28.9152698516846
0.83 -27.2524013519287
0.85 -25.1597518920898
0.87 -22.670316696167
0.89 -19.823356628418
0.91 -16.6637649536133
0.93 -13.2413787841797
0.95 -9.6101655960083
0.97 -5.82740020751953
0.99 -1.95273268222809
};
\label{lineA2}
\addplot [semithick, color2, mark=x, mark size=0.8, mark options={solid}, forget plot]
table {%
0.01 31.0377902984619
0.03 30.5483074188232
0.05 29.5770568847656
0.07 28.1393585205078
0.09 26.2578830718994
0.11 23.9623107910156
0.13 21.2888374328613
0.15 18.2796268463135
0.17 14.9821357727051
0.19 11.4483642578125
0.21 7.734046459198
0.23 3.89775967597961
0.25 -2.6501135153012e-07
0.27 -3.89776015281677
0.29 -7.73404693603516
0.31 -11.4483642578125
0.33 -14.9821357727051
0.35 -18.2796268463135
0.37 -21.2888374328613
0.39 -23.9623107910156
0.41 -26.2578830718994
0.43 -28.1393585205078
0.45 -29.5770568847656
0.47 -30.5483074188232
0.49 -31.0377902984619
0.51 -31.0377902984619
0.53 -30.5483074188232
0.55 -29.5770568847656
0.57 -28.1393585205078
0.59 -26.2578830718994
0.61 -23.9623107910156
0.63 -21.2888374328613
0.65 -18.2796268463135
0.67 -14.9821357727051
0.69 -11.4483642578125
0.71 -7.734046459198
0.73 -3.89775967597961
0.75 2.6501135153012e-07
0.77 3.89776015281677
0.79 7.73404693603516
0.81 11.4483642578125
0.83 14.9821357727051
0.85 18.2796268463135
0.87 21.2888374328613
0.89 23.9623107910156
0.91 26.2578830718994
0.93 28.1393585205078
0.95 29.5770568847656
0.97 30.5483074188232
0.99 31.0377902984619
};
\label{lineA3}
\addplot [semithick, black, mark=o, mark size=0.8, mark options={solid,fill opacity=0}, only marks, forget plot]
table {%
0.02 -110.821003451447
0.04 -110.820980145086
0.06 -110.820958659468
0.08 -110.820952684587
0.1 -110.820953417951
0.12 -110.820953580242
0.14 -110.820952962883
0.16 -110.820954027344
0.18 -110.820952117875
0.2 -110.820958998454
0.22 -110.820982110947
0.24 -110.821015210167
0.26 -110.821015210187
0.28 -110.820982110977
0.3 -110.820958998465
0.32 -110.820952117878
0.34 -110.820954027341
0.36 -110.820952962886
0.38 -110.820953580238
0.4 -110.820953417956
0.42 -110.820952684578
0.44 -110.820958659469
0.46 -110.820980145091
0.48 -110.821003451439
0.5 -110.821019112268
0.52 -110.821003451471
0.54 -110.820980145123
0.56 -110.820958659488
0.58 -110.820952684579
0.6 -110.820953417953
0.62 -110.820953580242
0.64 -110.82095296288
0.66 -110.82095402735
0.68 -110.820952117867
0.7 -110.82095899846
0.72 -110.82098211095
0.74 -110.821015210152
0.76 -110.821015210196
0.78 -110.820982111009
0.8 -110.820958998486
0.82 -110.820952117864
0.84 -110.820954027343
0.86 -110.820952962884
0.88 -110.820953580235
0.9 -110.820953417959
0.92 -110.82095268457
0.94 -110.820958659474
0.96 -110.820980145042
0.98 -110.821003451397
1 -110.821019112267
};
\addplot [semithick, black, mark=o, mark size=0.8, mark options={solid,fill opacity=0}, only marks, forget plot]
table {%
0.02 3.89776108564261
0.04 7.73404285752786
0.06 11.4483643543756
0.08 14.98213383632
0.1 18.279632222528
0.12 21.2888444567072
0.14 23.9623184694187
0.16 26.2578889833673
0.18 28.1393632178771
0.2 29.5770565580279
0.22 30.5483077922127
0.24 31.0377920406145
0.26 31.0377925550203
0.28 30.5483081339092
0.3 29.5770578818193
0.32 28.139364469216
0.34 26.2578904618977
0.36 23.9623191302417
0.38 21.2888479191749
0.4 18.2796344608915
0.42 14.9821416802449
0.44 11.4483653217719
0.46 7.73404172219905
0.48 3.89775993260326
0.5 2.21453878354623e-06
0.52 -3.89775491744911
0.54 -7.73403770469315
0.56 -11.448360853253
0.58 -14.9821371176876
0.6 -18.2796301907333
0.62 -21.2888441885219
0.64 -23.9623158264201
0.66 -26.2578869095879
0.68 -28.1393592281781
0.7 -29.5770566527972
0.72 -30.5483076635419
0.74 -31.0377924051259
0.76 -31.0377918461729
0.78 -30.548310105365
0.8 -29.5770572901704
0.82 -28.1393657837193
0.84 -26.257892555443
0.86 -23.9623227882469
0.88 -21.2888554066124
0.9 -18.2796292132956
0.92 -14.9821392435558
0.94 -11.4483675021979
0.96 -7.73404671585054
0.98 -3.8977720033196
1 1.29618686158826e-05
};
\addplot [semithick, black, mark=o, mark size=0.8, mark options={solid,fill opacity=0}, only marks, forget plot]
table {%
0.02 30.8539325558172
0.04 30.1221212018976
0.06 28.9152658249155
0.08 27.2523973491854
0.1 25.1597434856042
0.12 22.6703098276185
0.14 19.8233497681754
0.16 16.66376410538
0.18 13.2413794705872
0.2 9.61017133014731
0.22 5.82740404725478
0.24 1.95273343979389
0.26 -1.95273057859571
0.28 -5.82740112299587
0.3 -9.6101685515059
0.32 -13.2413759363826
0.34 -16.6637624991301
0.36 -19.8233483209912
0.38 -22.6703071887565
0.4 -25.1597412553621
0.42 -27.2523936618364
0.44 -28.9152646824611
0.46 -30.1221218428822
0.48 -30.8539320598005
0.5 -31.0991569392897
0.52 -30.8539326864647
0.54 -30.1221231892546
0.56 -28.9152665190456
0.58 -27.2523962795275
0.6 -25.1597441664244
0.62 -22.670310962164
0.64 -19.8233519535387
0.66 -16.6637685740145
0.68 -13.2413864378589
0.7 -9.61017312078584
0.72 -5.82740128778728
0.74 -1.95273789387571
0.76 1.95273159756248
0.78 5.8273942222633
0.8 9.61016829012942
0.82 13.2413746541983
0.84 16.6637579993506
0.86 19.8233449087013
0.88 22.6702992748779
0.9 25.1597458633248
0.92 27.252394267123
0.94 28.9152645112307
0.96 30.1221198963896
0.98 30.8539311834869
1 31.0991562958331
};
\label{hwplot1}
\end{axis}

\end{tikzpicture}
	\includegraphics[width=1\linewidth,height=1\linewidth]{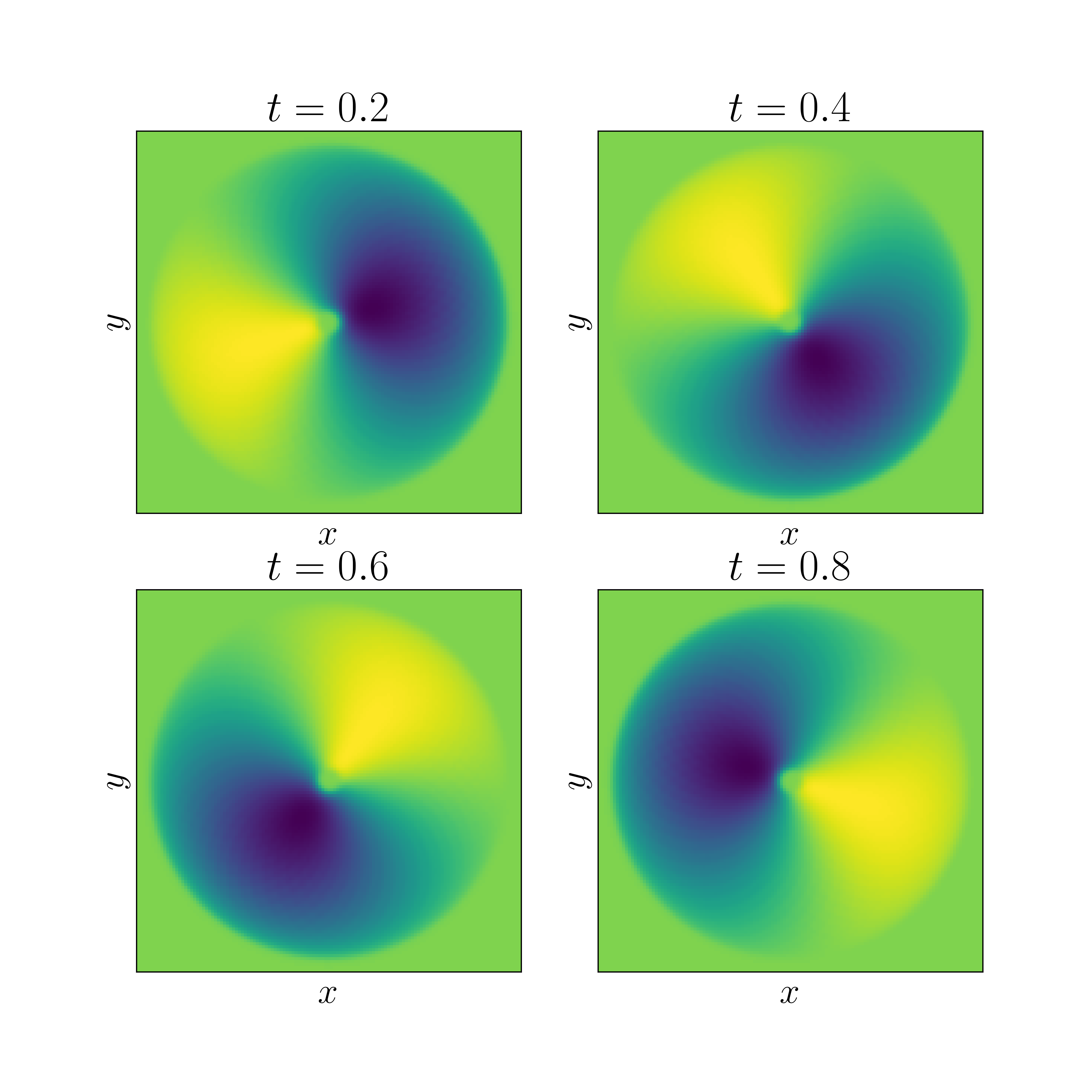}
	\end{center}
    \caption{POD with $r=3$}
    \label{subfig:POD-disk}
\end{subfigure}%
\hspace{0.5cm} %
\begin{subfigure}[t]{0.4\textwidth}
  \centering
		\scriptsize
% 		\textbf{FTR - Ansatz: }
% 		$\vec{q}(t)= f(\FTRmodes\matr{A}\vec{\Omega}(t))$
 		\begin{center}
            \setlength{\figureheight}{0.8\linewidth}
            \setlength{\figurewidth}{0.8\linewidth}
		    % This file was created with tikzplotlib v0.9.15.
\begin{tikzpicture}

\definecolor{color0}{rgb}{0.12156862745098,0.466666666666667,0.705882352941177}
\definecolor{color1}{rgb}{1,0.498039215686275,0.0549019607843137}
\definecolor{color2}{rgb}{0.172549019607843,0.627450980392157,0.172549019607843}

\begin{axis}[
height=\figureheight,
legend cell align={left},
legend style={
  fill opacity=0.8,
  draw opacity=1,
  text opacity=1,
  at={(0.03,0.97)},
  anchor=north west,
  draw=white!80!black
},
minor xtick={},
minor ytick={},
tick align=outside,
tick pos=left,
title={FTR \(\displaystyle \vec{q}(t) = f(\FTRmodes\mathbf{A}\vec{\Omega}(t))\)},
width=\figurewidth,
x grid style={white!69.0196078431373!black},
xlabel={time \(\displaystyle t\)},
xmin=-0.0395, xmax=1.0495,
xtick style={color=black},
xtick={-0.5,0,0.5,1,1.5},
xticklabels={
  \(\displaystyle {\ensuremath{-}0.5}\),
  \(\displaystyle {0.0}\),
  \(\displaystyle {0.5}\),
  \(\displaystyle {1.0}\),
  \(\displaystyle {1.5}\)
},
y grid style={white!69.0196078431373!black},
ylabel={coefficient \(\displaystyle \vec{a}(t)=\matr{A}\vec{\Omega}(t)\)},
ymin=-306.796581368581, ymax=261.121549643845,
ytick style={color=black},
ytick={-400,-200,0,200,400},
yticklabels={
  \(\displaystyle {\ensuremath{-}400}\),
  \(\displaystyle {\ensuremath{-}200}\),
  \(\displaystyle {0}\),
  \(\displaystyle {200}\),
  \(\displaystyle {400}\)
}
]
\addplot [semithick, color0, mark=x, mark size=0.8, mark options={solid}, forget plot]
table {%
0.01 -280.975240315565
0.03 -280.936442874551
0.05 -280.904046561962
0.07 -280.893569302584
0.09 -280.888517074398
0.11 -280.889839938761
0.13 -280.888355194699
0.15 -280.888926389136
0.17 -280.891098201216
0.19 -280.897046990654
0.21 -280.917625310647
0.23 -280.958121847229
0.25 -280.982120868016
0.27 -280.958121847155
0.29 -280.917625310502
0.31 -280.897046990443
0.33 -280.891098200998
0.35 -280.88892638886
0.37 -280.88835519442
0.39 -280.889839938476
0.41 -280.888517074148
0.43 -280.893569302383
0.45 -280.904046561804
0.47 -280.936442874459
0.49 -280.975240315576
0.51 -280.975240315609
0.53 -280.93644287456
0.55 -280.904046561965
0.57 -280.893569302581
0.59 -280.88851707438
0.61 -280.889839938723
0.63 -280.888355194671
0.65 -280.888926389088
0.67 -280.891098201154
0.69 -280.897046990554
0.71 -280.917625310559
0.73 -280.958121847145
0.75 -280.982120867944
0.77 -280.958121847062
0.79 -280.917625310412
0.81 -280.89704699034
0.83 -280.891098200869
0.85 -280.888926388762
0.87 -280.888355194316
0.89 -280.8898399384
0.91 -280.888517074069
0.93 -280.893569302296
0.95 -280.904046561758
0.97 -280.93644287441
0.99 -280.975240315534
};
\addplot [semithick, color1, mark=x, mark size=0.8, mark options={solid}, forget plot]
table {%
0.01 -234.837496474915
0.03 -231.102993375539
0.05 -223.730335180678
0.07 -212.847341777692
0.09 -198.611847908445
0.11 -181.248566285802
0.13 -161.024721991547
0.15 -138.262342230182
0.17 -113.319788599688
0.19 -86.5909946448963
0.21 -58.4992289215098
0.23 -29.4851242186596
0.25 -0.000208005251156441
0.27 29.4847115228499
0.29 58.4988260706268
0.31 86.5906079618368
0.33 113.319424161984
0.35 138.262005780372
0.37 161.024418834074
0.39 181.248301200073
0.41 198.611625077944
0.43 212.847164712877
0.45 223.730206672519
0.47 231.102915445578
0.49 234.837470358358
0.51 234.837496474912
0.53 231.102993375515
0.55 223.730335180659
0.57 212.847341777677
0.59 198.611847908427
0.61 181.248566285778
0.63 161.024721991535
0.65 138.262342230165
0.67 113.319788599671
0.69 86.5909946448683
0.71 58.4992289214875
0.73 29.485124218641
0.75 0.000208005232930895
0.77 -29.4847115228668
0.79 -58.4988260706445
0.81 -86.5906079618492
0.83 -113.319424161984
0.85 -138.262005780386
0.87 -161.024418834081
0.89 -181.248301200092
0.91 -198.611625077956
0.93 -212.847164712878
0.95 -223.730206672545
0.97 -231.102915445594
0.99 -234.837470358371
};
\addplot [semithick, color2, mark=x, mark size=0.8, mark options={solid}, forget plot]
table {%
0.01 -14.7720328555589
0.03 -44.0791068369317
0.05 -72.6875348055381
0.07 -100.152535212683
0.09 -126.038884354048
0.11 -149.939627359829
0.13 -171.474091003592
0.15 -190.305188090131
0.17 -206.135959348975
0.19 -218.718576525058
0.21 -227.863581324871
0.23 -233.432178776181
0.25 -235.307089143228
0.27 -233.432230903863
0.29 -227.863684747876
0.31 -218.718729612748
0.33 -206.136159691739
0.35 -190.305432529961
0.37 -171.474375686099
0.39 -149.939947797017
0.41 -126.039235488624
0.43 -100.152911514892
0.45 -72.6879303483397
0.47 -44.0795154142443
0.49 -14.7724480353247
0.51 14.7720328556044
0.53 44.0791068369746
0.55 72.6875348055839
0.57 100.152535212729
0.59 126.038884354092
0.61 149.939627359868
0.63 171.474091003642
0.65 190.305188090176
0.67 206.135959349014
0.69 218.718576525075
0.71 227.863581324902
0.73 233.43217877622
0.75 235.307089143281
0.77 233.432230903902
0.79 227.863684747918
0.81 218.718729612781
0.83 206.136159691753
0.85 190.305432529996
0.87 171.474375686127
0.89 149.939947797062
0.91 126.039235488666
0.93 100.152911514928
0.95 72.6879303483866
0.97 44.0795154142885
0.99 14.7724480353691
};
\addplot [semithick, black, mark=o, mark size=0.8, mark options={solid,fill opacity=0}, only marks, forget plot]
table {%
0.02 -280.958434387181
0.04 -280.916702870231
0.06 -280.897613074499
0.08 -280.89017663575
0.1 -280.889010698116
0.12 -280.889405170156
0.14 -280.888144096996
0.16 -280.889950222203
0.18 -280.893103643076
0.2 -280.904556378391
0.22 -280.9366219513
0.24 -280.975455189781
0.26 -280.975455207989
0.28 -280.936621992596
0.3 -280.904556391198
0.32 -280.893103652527
0.34 -280.889950218076
0.36 -280.888144101956
0.38 -280.889405162303
0.4 -280.889010706793
0.42 -280.890176618294
0.44 -280.897613077647
0.46 -280.916702866071
0.48 -280.958434385651
0.5 -280.981432160847
0.52 -280.958434440189
0.54 -280.916702912348
0.56 -280.897613087687
0.58 -280.890176629594
0.6 -280.889010699417
0.62 -280.889405171412
0.64 -280.888144093229
0.66 -280.889950225438
0.68 -280.893103621018
0.7 -280.904556388494
0.72 -280.936621955003
0.74 -280.975455174305
0.76 -280.975455214898
0.78 -280.936622037772
0.8 -280.90455641765
0.82 -280.893103640799
0.84 -280.889950227336
0.86 -280.888144097909
0.88 -280.889405156855
0.9 -280.889010715036
0.92 -280.890176606238
0.94 -280.897613093465
0.96 -280.916702783701
0.98 -280.95843429833
1 -280.981432173948
};
\addplot [semithick, black, mark=o, mark size=0.8, mark options={solid,fill opacity=0}, only marks, forget plot]
table {%
0.02 -233.432753879434
0.04 -227.863405616959
0.06 -218.719373453015
0.08 -206.135822116919
0.1 -190.305552216227
0.12 -171.47534942548
0.14 -149.939070446652
0.16 -126.039938043697
0.18 -100.152654963945
0.2 -72.6883185291539
0.22 -44.0795072961844
0.24 -14.7726422928021
0.26 14.7722017495807
0.28 44.0790725374186
0.3 72.6879003313201
0.32 100.152252048757
0.34 126.039575091796
0.36 149.938738324372
0.38 171.475046007278
0.4 190.305288945457
0.42 206.135597074696
0.44 218.719206053901
0.46 227.86331878203
0.48 233.432664461454
0.5 235.306951451273
0.52 233.432721353991
0.54 227.863432471023
0.56 218.719372976716
0.58 206.135817324906
0.6 190.30555527545
0.62 171.475359417735
0.64 149.939086001033
0.66 126.039972535958
0.68 100.152707189187
0.7 72.6883323884669
0.72 44.0794862532034
0.74 14.7726760408198
0.76 -14.7722093883853
0.78 -44.0790205263571
0.8 -72.6878980344215
0.82 -100.15224282806
0.84 -126.039540357993
0.86 -149.938713494785
0.88 -171.47498472182
0.9 -190.305325906491
0.92 -206.135598335318
0.94 -218.719210512292
0.96 -227.863292203118
0.98 -233.432691317348
1 -235.306708975061
};
\addplot [semithick, black, mark=o, mark size=0.8, mark options={solid,fill opacity=0}, only marks, forget plot]
table {%
0.02 -29.4863043001541
0.04 -58.4976467869514
0.06 -86.5912911911149
0.08 -113.3186447984
0.1 -138.262383598563
0.12 -161.02451454002
0.14 -181.247259092884
0.16 -198.612235086647
0.18 -212.846728349536
0.2 -223.73009456954
0.22 -231.102842304059
0.24 -234.837194553878
0.26 -234.837462521849
0.28 -231.102681129767
0.3 -223.730482689962
0.32 -212.846652677975
0.34 -198.612749882116
0.36 -181.247221866044
0.38 -161.025189146691
0.4 -138.262336156613
0.42 -113.319554758477
0.44 -86.591058035099
0.46 -58.4989573702814
0.48 -29.4850801025023
0.5 -0.000227184961614479
0.52 29.4846273143454
0.54 58.4985174201683
0.56 86.590633696997
0.58 113.319156861111
0.6 138.261966834766
0.62 161.02485814065
0.64 181.246931504399
0.66 198.612500403895
0.68 212.846435755271
0.7 223.73034504745
0.72 231.102599455506
0.74 234.837435390818
0.76 234.837219074305
0.78 231.102937954728
0.8 223.730228487241
0.82 212.846925006439
0.84 198.612484731406
0.86 181.247557119616
0.88 161.02490016449
0.9 138.262697848237
0.92 113.319049089049
0.94 86.5917055326656
0.96 58.4980855350337
0.98 29.4868017139986
1 -0.00582821479338784
};
\end{axis}

\end{tikzpicture}
		    \includegraphics[width=1\linewidth,height=1\linewidth]{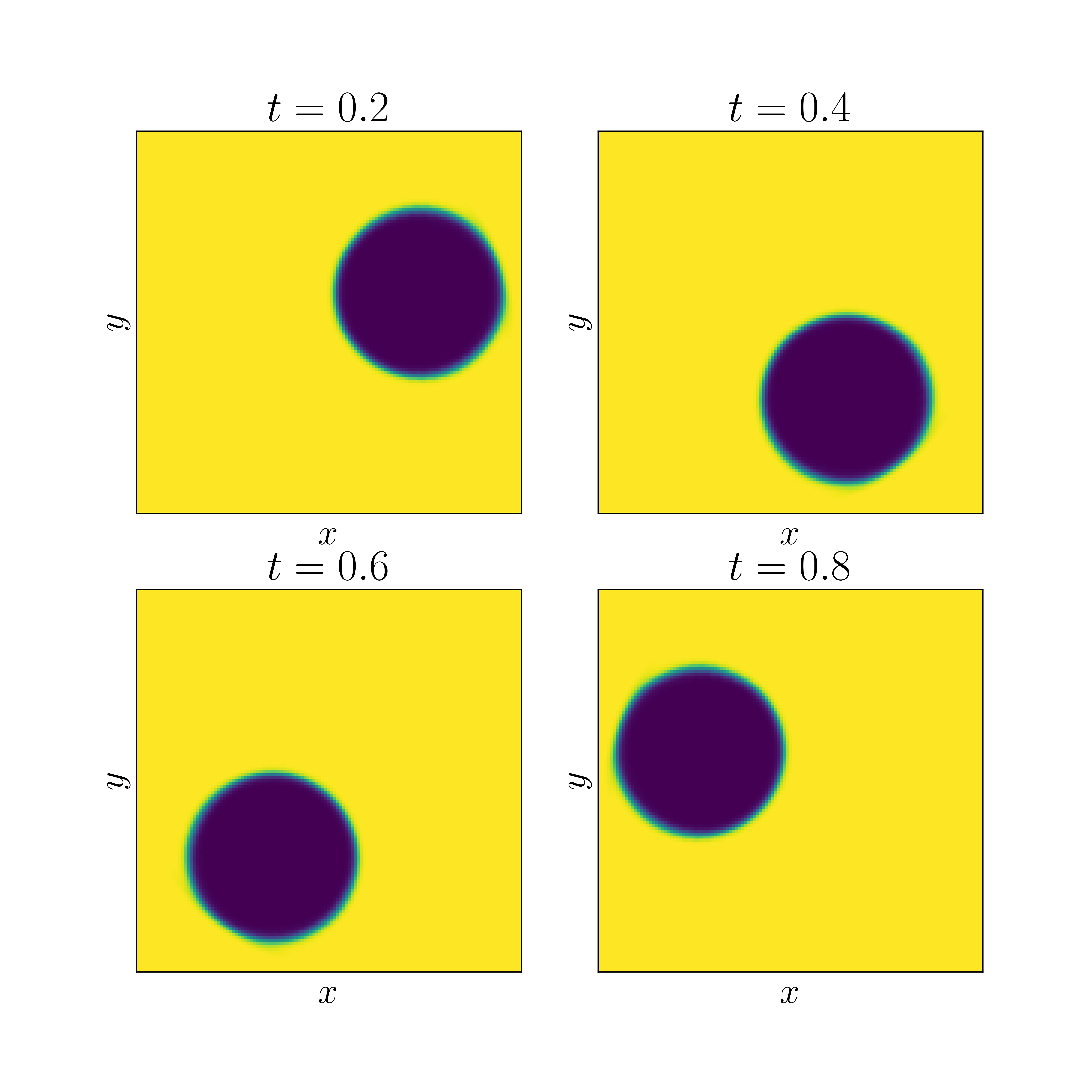}

	\end{center}
    \caption{FTR with $r=3$}
    \label{subfig:FTR-disk}
  \end{subfigure}
  \caption{Predictions using Fourier-Koopman forecasting with three POD modes (a) and three FTR modes (b). The black circles (\ref{hwplot1}) in the upper row indicate the predictions of the amplitudes $\vec{a}(t)=(a_1(t),a_2(t),a_3(t))\widehat{=}$ (\ref{lineA1},\ref{lineA2},\ref{lineA3}) and the colored crosses mark the training samples.
	In the lower row, we show the corresponding snapshots at selected time instances $t=0.2,0.4,0.6,0.8$.}
  \label{fig:mov-disc:online}
\end{figure}
%%%%%%%%%%%%%%%%%%%%%%%%%%%%%%%%%%%%%%
Furthermore, the online-prediction error is stated for $r=2,4,6,8,10,12,15$ in \cref{tabl:online-error-fourier-koopman}.

Note that after solving \cref{eq:fourier-koopman-opt} in the offline stage, the computational effort is reduced to the evaluation of $\tilde{\vec{q}}(t)=f(\FTRmodes\matr{A}\vec{\Omega}(t))$, which only takes milliseconds.

\begin{figure}[htp!]
    \centering
    \begin{subfigure}[t]{0.63\textwidth}
    \centering
    \includegraphics[width=0.48\textwidth,trim={2cm 0cm 1cm 0.5cm},clip]{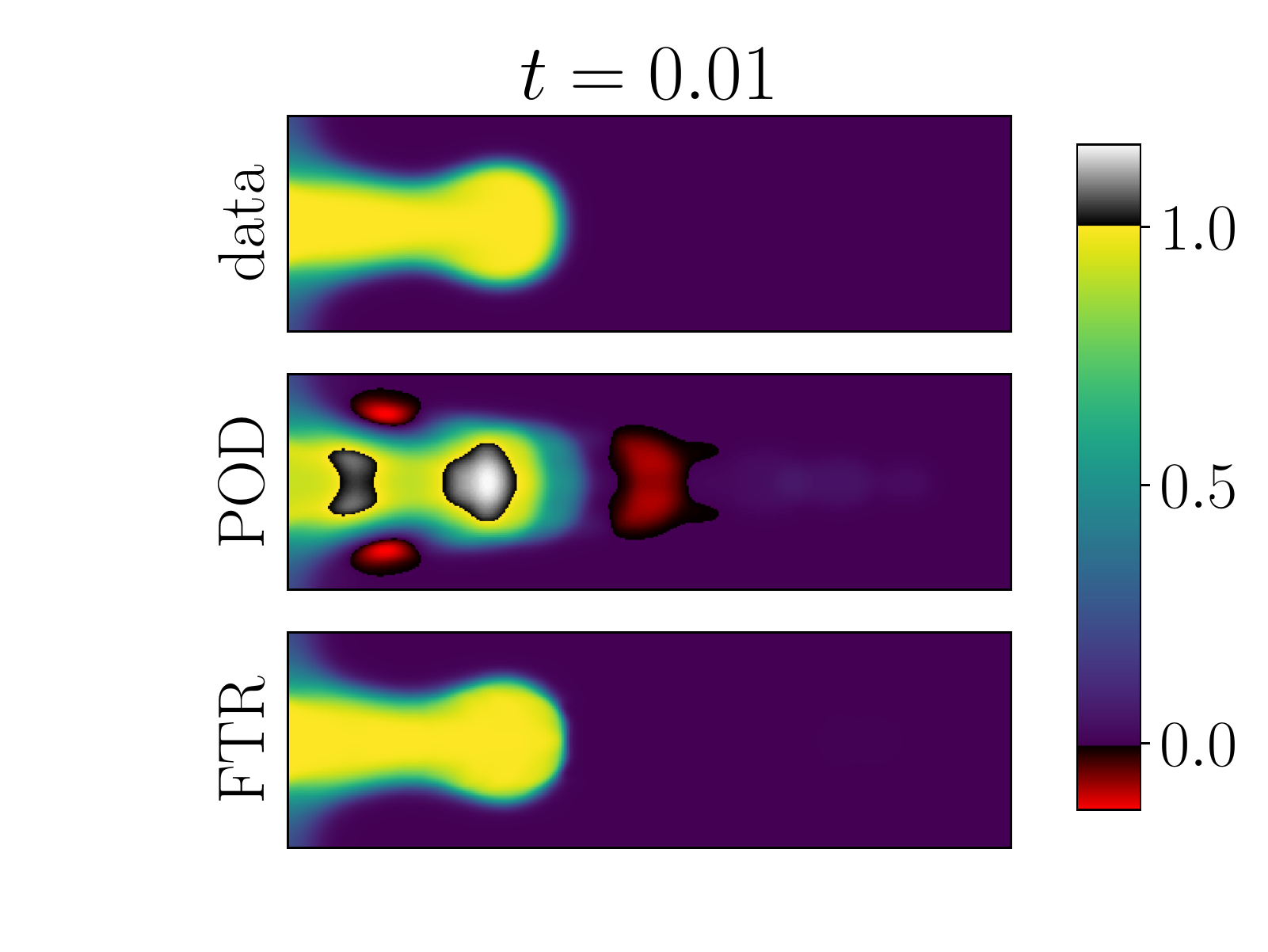}
    \includegraphics[width=0.48\textwidth,trim={2cm 0cm 1cm 0.5cm},clip]{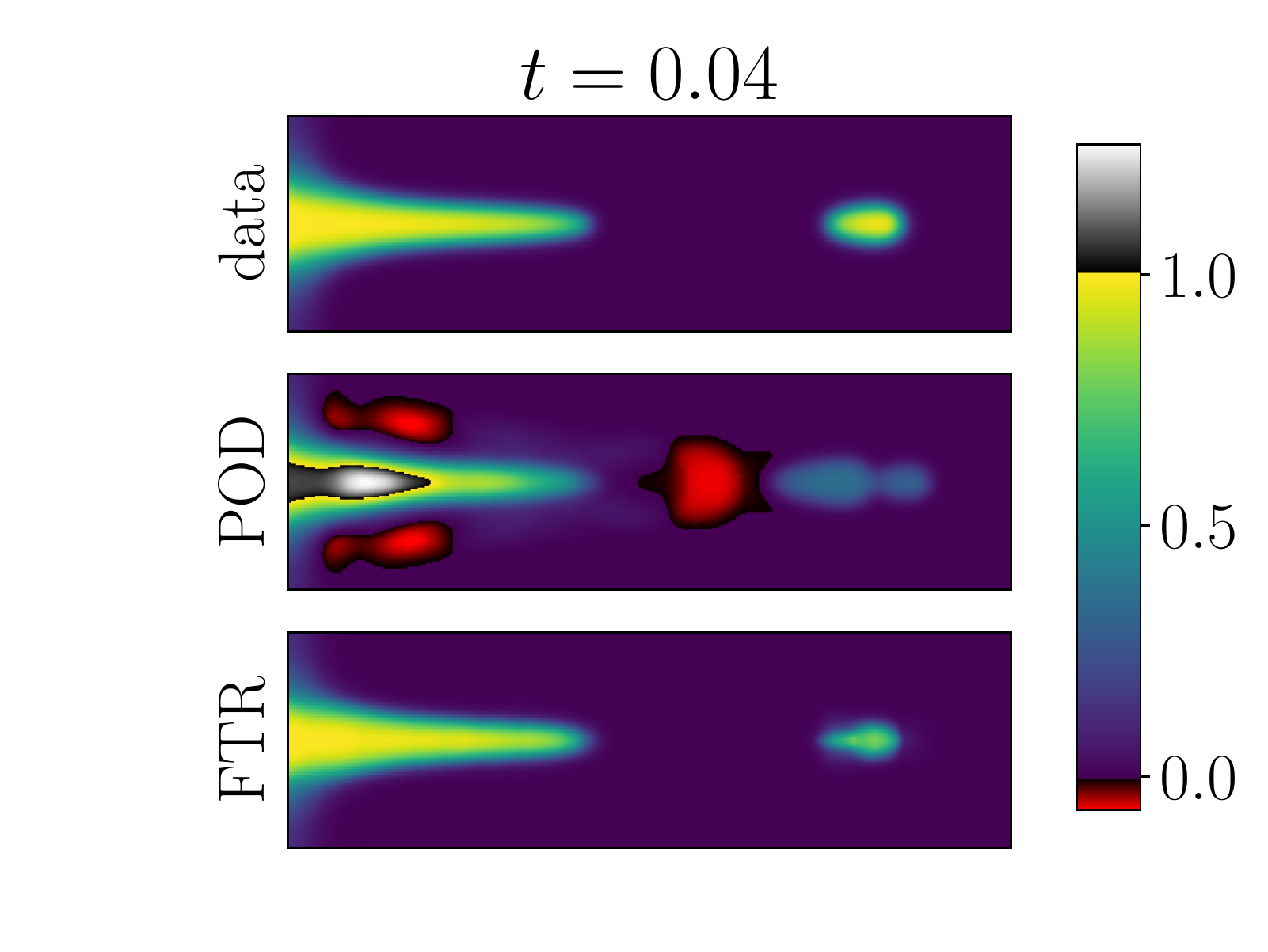}
    \caption{Test snapshots}
    \label{subfig:flame_pinch_off_compare:recon}
    \end{subfigure}%
    \begin{subfigure}[t]{0.38\textwidth}
    \centering
    \setlength{\figureheight}{0.9\linewidth}
    \setlength{\figurewidth}{1\linewidth}
    \input{FTR-online-bunsen_time.tex}
    \caption{FTR-Koopman predictions}
    \label{subfig:flame_pinch_off_compare:time}
    \end{subfigure}
    \caption{Online predictions of the Bunsen flame example. Fig. a) compares the test data in the top row with the FTR-Koopman and POD-Koopman results using $r=8$ degrees of freedom for $t=0.01$ and 0.04. The snapshots show how a burning fuel pocket is detached from the flame at $t=0.04$ causing a change in the topology of the contour line of the front. Fig. b) visualizes the Fourier-Koopman predictions (\ref{tikz:dots_FTR_bunsen}) for $\vec{a}(t)=(a_1(t),a_2(t),a_3(t))\widehat{=}$ (\ref{lineB1},\ref{lineB2},\ref{lineB3}).}
    \label{fig:flame_pinch_off_comparison}
\end{figure}

For a realistic test case, we apply the FTR-Fourier-Koopman procedure to the methane mass fraction $Y_{\mathrm{CH}_4}$ of one flame of a multi-slit Bunsen burner simulation analyzed and studied in \cite{KornilovRookThijeDeGoey2009,JaenschMerkGopalakrishnanBombergEmmertSujithPolifke2017}. The snapshots are generated with a customized, weakly compressible version of \texttt{rhoReactionFOAM} from the \texttt{OpenFOAM} software package (see \cite{WellerTaborJasakFureby1998,JaenschMerkGopalakrishnanBombergEmmertSujithPolifke2017}). In the simulation, a flame is periodically excited by an incoming velocity pulse. The acceleration of the fuel detaches a burning pocket shown in \cref{fig:flame_pinch_off_comparison}. The data set consists of 200 snapshots, with $M = 128\times 430$ grid points, sampled in a time interval $t\in[0.01,0.05]$ in which the Bunsen flame is quasi-periodic. Again, we split the data into train ($t_n=2\Delta t n$) and test samples $t_{n+1/2}=(2 n+ 1)\Delta t$. While we use the train samples to generate the reduced model, the test samples are used to calculate the relative errors stated in \cref{tabl:online-error-fourier-koopman}. The flame pinch-off is not a special case in combustion systems, but it poses challenges to model order reduction methods, as described above.
\Cref{fig:flame_pinch_off_comparison} shows that for the FTR the structure of the solution is well captured and the physical bound $0\le Y_{\mathrm{CH}_4}\le 1$ is preserved.
%%%%%%%%%%%%%%%%%%%%%%%%%%%%%%%
%%%%%%%%%%%%%%%%%%%%%%%%%%%%%%%

 %%%%%%%%%%%%%%%%%%%%%%%%%%%%%%
\begin{table}[htp!]
    \centering
        \begin{tabular}{c cc cc}
    \toprule
     & \multicolumn{2}{c}{Moving Disc} & \multicolumn{2}{c}{Bunsen Flame} \\
    \cmidrule(lr){2-3} \cmidrule(lr){4-5}
    rank $r$     & FTR   & POD     & FTR  & POD  \\
    \midrule
    2   & 2.7e-01   &3.0e-01   &4.2e-01      & 3.1e-01   \\
    4   & 7.4e-03   &2.0e-01   &1.4e-01      & 3.1e-01    \\
    6   & 2.2e-03   &1.5e-01   &1.1e-01       & 2.3e-01    \\
    8   & 1.6e-03   &1.2e-01   &7.6e-02       & 1.8e-01    \\
    10  & 2.2e-03   &1.0e-01   &8.1e-02      & 1.6e-01   \\
    12  & 2.0e-03   &8.8e-02   &7.1e-02      & 1.5e-01    \\
    15  & 1.2e-03   &7.4e-02   &6.9e-02      & 1.4e-01    \\
    \bottomrule
    \end{tabular}
    \caption{Relative error $\sum_{n=0}^{2N-1}\Vert \vec{q}(t_{n+1/2})-\tilde{\vec{q}}(t_{n+1/2})\Vert_2^2/\sum_{n=0}^{2N-1}\Vert \vec{q}(t_{n+1/2})\Vert_2^2$ for the FTR-Fourier-Koopman predictions using the moving disk and Bunsen flame data.
    }
    \label{tabl:online-error-fourier-koopman}
\end{table}
%%%%%%%%%%%%%%%%%%%%%%%%%%%%%%%%%%%%%%

%%%%%%%%%%%%%%%%%%%%%%%%%%%%%%%%%%%%%%%%%%%%%%%%%%%%%%%%%%%%%%%%%%%%%%%%%
% Manifold Galerkin
%%%%%%%%%%%%%%%%%%%%%%%%%%%%%%%%%%%%%%%%%%%%%%%%%%%%%%%%%%%%%%%%%%%%%%%%%
\subsection{Manifold Galerkin Methods}
\label{subsec:manifoldgalerkin}
\begin{figure}[htp!]
    \centering
        \setlength{\figureheight}{0.2\linewidth}
        \setlength{\figurewidth}{0.6\linewidth}
    \begin{center}
        training
        \vspace{-0.3cm}
        \begin{equation*}
            \partial_t q + u(t) \partial_x q = 0
        \end{equation*}
    \end{center}
    \vspace{-0.4cm}
\begin{tikzpicture}

\begin{axis}[%
width=0.461\figurewidth,
height=\figureheight,
at={(0\figurewidth,0\figureheight)},
scale only axis,
point meta min=-0.0000910952,
point meta max=1.0010239820,
axis on top,
xmin=0.5000000000,
xmax=1000.5000000000,
xtick={\empty},
xlabel style={font=\color{white!15!black}},
xlabel={$x$},
ymin=0.5000000000,
ymax=101.5000000000,
ytick={\empty},
ylabel style={font=\color{white!15!black}},
ylabel={$t$},
axis background/.style={fill=white},
title style={font=\bfseries},
title={$u(t)=-2$}
]
\addplot [forget plot] graphics [xmin=0.5000000000, xmax=1000.5000000000, ymin=0.5000000000, ymax=101.5000000000] {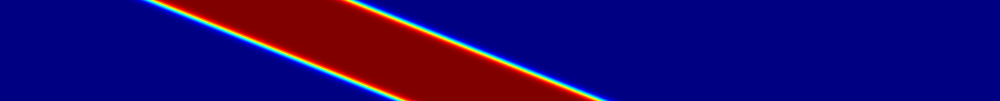};
\end{axis}

\begin{axis}[%
width=0.42\figurewidth,
height=\figureheight,
at={(0.607\figurewidth,0\figureheight)},
scale only axis,
point meta min=0.0000000000,
point meta max=1.0000000000,
axis on top,
xmin=0.5000000000,
xmax=1000.5000000000,
xtick={\empty},
xlabel style={font=\color{white!15!black}},
xlabel={$x$},
ymin=0.5000000000,
ymax=101.5000000000,
ytick={\empty},
ylabel style={font=\color{white!15!black}},
ylabel={$t$},
axis background/.style={fill=white},
title style={font=\bfseries},
title={$u(t) = 2$},
colormap/jet,
colorbar right,
colorbar style={anchor=south west, at={(1.03,0)}, height=1*\pgfkeysvalueof{/pgfplots/parent axis height}}
]
\addplot [forget plot] graphics [xmin=0.5000000000, xmax=1000.5000000000, ymin=0.5000000000, ymax=101.5000000000] {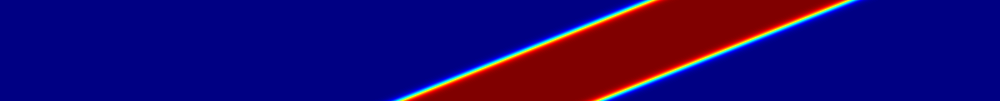};
\end{axis}
\end{tikzpicture}%
    \begin{center}
        testing
        \vspace{-0.3cm}
        \begin{equation*}
            \partial_t q + 5 \sin(2\pi t/T) \partial_x q = 0
        \end{equation*}
    \end{center}
    \vspace{-0.4cm}
\begin{tikzpicture}

\begin{axis}[%
width=0.46\figurewidth,
height=\figureheight,
at={(0\figurewidth,0\figureheight)},
scale only axis,
point meta min=-0.0001422240,
point meta max=1.0015137233,
axis on top,
xmin=0.5000000000,
xmax=1000.5000000000,
xtick={\empty},
xlabel style={font=\color{white!15!black}},
xlabel={$x$},
ymin=0.5000000000,
ymax=101.5000000000,
ytick={\empty},
ylabel style={font=\color{white!15!black}},
ylabel={$t$},
axis background/.style={fill=white},
title style={font=\bfseries},
title={FOM}
]
\addplot [forget plot] graphics [xmin=0.5000000000, xmax=1000.5000000000, ymin=0.5000000000, ymax=101.5000000000] {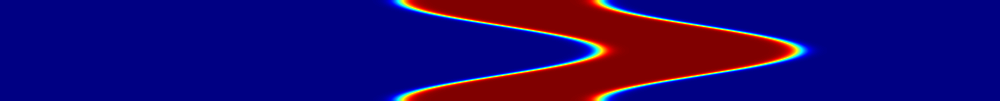};
\end{axis}

\begin{axis}[%
width=0.42\figurewidth,
height=\figureheight,
at={(0.605\figurewidth,0\figureheight)},
scale only axis,
point meta min=0.0000000000,
point meta max=1.0000000000,
axis on top,
xmin=0.5000000000,
xmax=1000.5000000000,
xtick={\empty},
xlabel style={font=\color{white!15!black}},
xlabel={$x$},
ymin=0.5000000000,
ymax=101.5000000000,
ytick={\empty},
ylabel style={font=\color{white!15!black}},
ylabel={$t$},
axis background/.style={fill=white},
title style={font=\bfseries},
title={FTR - ROM},
colormap/jet,
colorbar right,
colorbar style={anchor=south west, at={(1.03,0)}, height=1*\pgfkeysvalueof{/pgfplots/parent axis height}}
]
\addplot [forget plot] graphics [xmin=0.5000000000, xmax=1000.5000000000, ymin=0.5000000000, ymax=101.5000000000] {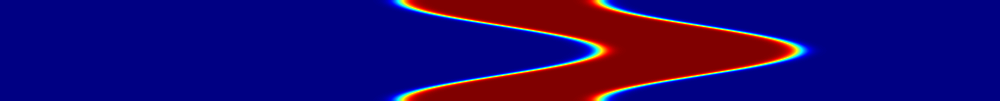};
\end{axis}
\end{tikzpicture}%
    \caption{The 1D advection test case for moving fronts. In the upper row, the two training samples computed with the full order model (FOM) are shown. They are used to build the snapshot matrix $\matr{Q}\in\mathbb{R}^{1000\times 202}$ for the FTR decomposition. In the lower row, the trajectory of the FOM and the FTR-ROM \cref{eq:MLgalerkin-advect} ($r=4$) are compared.}
    \label{fig:1D-advection}
\end{figure}

After discretizing the ARD system \cref{eq:react-diff-advect} in space, we obtain an ODE system of the form
\begin{equation}
    \label{eq:FOM}
\text{(FOM)}\quad
		\left\{\begin{aligned}
		\dot{\q}(t,\mu) &=\RHS(\q,t,\mu)\,\\
				\q(0) &= \q_0\,.
	\end{aligned}\right.
 \end{equation}
Here, the parameters  $\mu\in \params$ contain the velocity field $\vec{u}$, diffusion or reaction constant $\diffconst,\reactconst$.
After discretizing the rescaled system it yields the FOM-RHS
\begin{equation}
\label{eq:ODE-rhs}
 \RHS(\q,t,\mu) = \LinOp(t) \q + \mu\NonLinOp(\q)
 \end{equation}
 with a linear operator $\LinOp\colon [0,T]\to \mathbb{R}^{M\times M}$ and a nonlinear operator $\NonLinOp\colon \mathbb{R}^M \to \mathbb{R}^M$.
 Using a reduced mapping
 \begin{equation}
    g\colon \mathbb{R}^r \to \mathbb{R}^M: \vec{a}\mapsto g(\vec{a}),\quad\text{with Jacobian} \quad \matr{J}_g(\vec{a}) = \left({\frac {\partial g_i}{\partial a_{j}}}(\vec{a})\right)_{\substack{i=1,\ldots,M\\j=1,\ldots,r}}
 \end{equation}
as approximation $\q\approx\tilde{\q}=g(\vec{a})$ of the data and plugging it into \cref{eq:FOM} yields a reduced model:
\begin{empheq}[left=\text{(ROM)}\quad\empheqlbrace]{align}
	 		\dot{\vec{a}}(t,\mu)&=\argmin_{\dot{\vec{a}}\in\mathbb{R}^r}\Vert\matr{J}_g(\vec{a})\dot{\vec{a}}(t,\mu) -\RHS(g(\vec{a}),t,\mu)\Vert_2^2 \label{eq:ROM}\\
 			\vec{a}(0, \mu)&= \argmin_{\vec{a}\in\mathbb{R}^r}\Vert\q_0-g(\vec{a})\Vert_2^2\,. \label{eq:inicond_ROM}
\end{empheq}
 %\end{equation}
 Minimizing the continuous time residual \cref{eq:ROM},
 yields the optimality condition:
 \begin{align}
    0 &= \frac{\d}{\d{\dot{\vec{a}}}}\Vert\matr{J}_g(\vec{a})\dot{\vec{a}} -\RHS(g(\vec{a}),t,\mu)\Vert_2^2\\
      &= 2 \matr{J}_g(\vec{a})^T\matr{J}_g(\vec{a})\dot{\vec{a}} - 2\matr{J}_g(\vec{a})^T  \RHS(g(\vec{a}),t,\mu)\,,\label{eq:MLgalerkin-optimality}
 \end{align}
 which is uniquely solved by
 \begin{equation}
 \label{eq:MLgalerkin}
     \vec{\dot{a}} = \matr{J}_g^+(\vec{a}) \RHS(g(\vec{a}),t,\mu)\,,
 \end{equation}
 if the Jacobin has full column rank \cite{LeeCarlberg2020}.
 Here, $ \matr{J}_g^+$ is the Moore-Penrose pseudo inverse of $\matr{J}_g$.
Note, that for the common POD-Galerkin approach orthogonal mappings $g(\vec{a}) = \matr{U}\vec{a}$, with $\matr{U}^T\matr{U}=\matr{I}$ are used. Therefore, \cref{eq:MLgalerkin} and \cref{eq:MLgalerkin-optimality} are identical. When
neglecting $\vec{N}(\q)$ in \cref{eq:ODE-rhs} for the time being, we obtain a small $r$-dimensional system:
 \begin{equation}
 \label{eq:MLgalerkin-advect0}
     \vec{\dot{a}} = \matr{L}_r(t)\vec{a}\,\quad\text{with}\quad \matr{L}_r(t)=\matr{U}^T\matr{L}(t)\matr{U}\in\mathbb{R}^{r\times r}\,,
 \end{equation}
which can be solved efficiently, when $\matr{L}_r(t)$ is precomputed.
For example in the case of pure advection:
\begin{equation}
\label{eq:pure-advection}
    \dot{q}(\vec{x},t)=\vec{u}\cdot\vec{\nabla}q = \sum_{k=1}^d u_k(t) \partial_{x_k} q(\vec{x},t)\,,
\end{equation}
the spatial derivative has the form $\matr{L}(t) = \sum_{k=1}^d u_k(t) \matr{L}^{(k)}$ and therefore
\begin{equation}
    \label{eq:reduced-advec-operator}
    \matr{L}_r(t) = \sum_{k=1}^d u_k(t) \matr{U}^T \matr{L}^{(k)} \matr{U} \in \mathbb{R}^{r\times r}
\end{equation}
can be precomputed and is much smaller than the operator $\matr{L}(t)\in\mathbb{R}^{M\times M}$, $r\ll M$ of the FOM.
Although POD-Galerkin enables to solve \cref{eq:MLgalerkin-advect0} efficiently, this approach cannot be used for advection dominated systems, because of its slow decaying approximation errors.
Here, nonlinear methods like artificial neural networks can accelerate the convergence of the overall online and offline error. However, any nonlinear reduction method will imply that even linear systems like  \cref{eq:pure-advection} become nonlinear, causing additional effort for evaluating nonlinearities. At least in the special case of an advection system, this can be avoided with the FTR approach.
Due to its special structure $\tilde{q}(\vec{x},t)=f(\phi(\vec{x},t))$, we can rewrite the advection equation $\partial_t q-\vec{u}\cdot\vec{\nabla}q=0$ into the form
\begin{equation}
    f'(\phi)(\partial_t\phi-\vec{u}\cdot\vec{\nabla}\phi)=0\,.
\end{equation}
The prefactor $f'(\phi)$ can be dropped, when assuming that $\phi$ features the same transport then $q$ and thus the nonlinear manifold Galerkin system \eqref{eq:MLgalerkin} can be simplified to a linear Galerkin system for $\vec{\phi}(t)=\FTRmodes \vec{a}(t)$:
 \begin{equation}
 \label{eq:MLgalerkin-advect}
     \vec{\dot{a}} =\FTRmodes\matr{L}(t)\FTRmodes\vec{a}\approx\matr{L}_r(t)\vec{a}\,.
 \end{equation}
 Since the operator $\matr{L}_r$ can be precomputed in the same fashion as for the POD-Galerkin approach, the resulting ROM complexity is reduced and the online/offline error is compensated due to the additional nonlinearity $f$ to retain $\vec{q}=f(\FTRmodes \vec{a})$.
 However, these findings need to be interpreted with caution. One might expect that a pure transport of $q$ implies a pure transport of $\phi$. However, if $f'$ becomes (approximately) zero, $ \phi$  might locally change its value without changing $q$, so that $q$ is transported everywhere, while $\phi$ is not.
 If, however, this cancellation is justified, it can speed up the calculation considerably. The advection of fronts according to a given transport field $\vec{u}(t)$ within milliseconds is impressively shown for a 1D advection example in \cref{fig:1D-advection}. In this example, only two trajectories of constant advection speed ($u(t)=\pm 2$) are used for building the reduced system. Thereafter, almost any parameterization of $u(t)$ can be computed with the ROM.

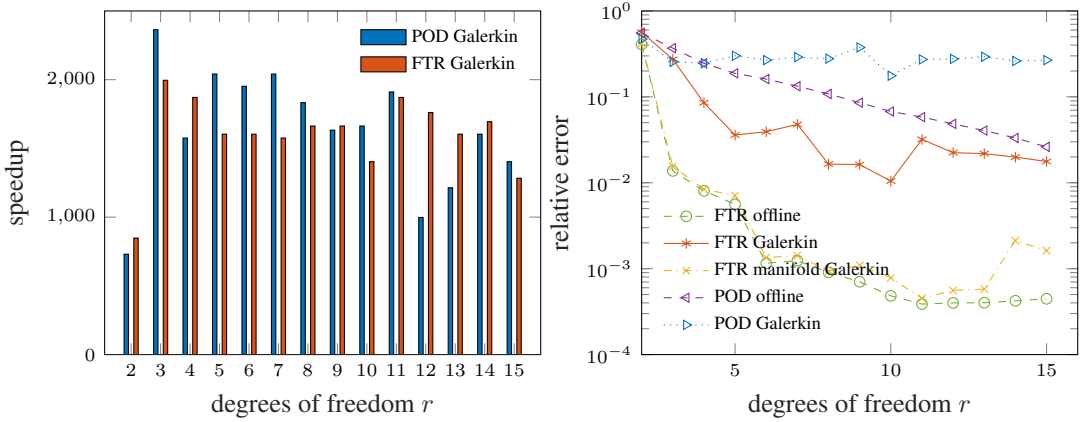
\begin{figure}[htp!]
    \centering
    \setlength{\figureheight}{0.3\linewidth}
    \setlength{\figurewidth}{0.4\linewidth}
    \hspace*{-0.5cm}% This file was created by matlab2tikz.
%
\definecolor{mycolor1}{rgb}{0.00000,0.44700,0.74100}%
\definecolor{mycolor2}{rgb}{0.85000,0.32500,0.09800}%
\begin{tikzpicture}

\begin{axis}[%
width=0.951\figurewidth,
height=\figureheight,
at={(0\figurewidth,0\figureheight)},
scale only axis,
bar shift auto,
xmin=1.1000000000,
xmax=15.9000000000,
xtick={ 2,  3,  4,  5,  6,  7,  8,  9, 10, 11, 12, 13, 14, 15},
xlabel style={font=\color{white!15!black}},
xlabel={degrees of freedom $r$},
ymin=0.0000000000,
ymax=2500.0000000000,
ylabel style={font=\color{white!15!black}},
ylabel={speedup},
axis background/.style={fill=white},
legend style={legend cell align=left, align=left, fill=none, draw=none}
]
\addplot[ybar, bar width=0.15, fill=mycolor1, draw=black, area legend] table[row sep=crcr] {%
2.0000000000	730.0813008132\\
3.0000000000	2363.1578947530\\
4.0000000000	1575.4385964920\\
5.0000000000	2040.9090909151\\
6.0000000000	1952.1739130471\\
7.0000000000	2040.9090908982\\
8.0000000000	1832.6530612184\\
9.0000000000	1632.7272727186\\
10.0000000000	1662.9629629714\\
11.0000000000	1910.6382978675\\
12.0000000000	997.7777777761\\
13.0000000000	1213.5135135108\\
14.0000000000	1603.5714285676\\
15.0000000000	1403.1250000012\\
};
\addplot[forget plot, color=white!15!black] table[row sep=crcr] {%
1.1000000000	0.0000000000\\
15.9000000000	0.0000000000\\
};
\addlegendentry{POD Galerkin}

\addplot[ybar, bar width=0.15, fill=mycolor2, draw=black, area legend] table[row sep=crcr] {%
2.0000000000	847.1698113226\\
3.0000000000	1995.5555555522\\
4.0000000000	1870.8333333349\\
5.0000000000	1603.5714285780\\
6.0000000000	1603.5714285676\\
7.0000000000	1575.4385964920\\
8.0000000000	1662.9629629602\\
9.0000000000	1662.9629629602\\
10.0000000000	1403.1250000012\\
11.0000000000	1870.8333333349\\
12.0000000000	1760.7843137184\\
13.0000000000	1603.5714285676\\
14.0000000000	1694.3396226451\\
15.0000000000	1282.8571428557\\
};
\addplot[forget plot, color=white!15!black] table[row sep=crcr] {%
1.1000000000	0.0000000000\\
15.9000000000	0.0000000000\\
};
\addlegendentry{FTR Galerkin}

\end{axis}
\end{tikzpicture}%
\definecolor{mycolor1}{rgb}{0.00000,0.44700,0.74100}%
\definecolor{mycolor2}{rgb}{0.85000,0.32500,0.09800}%
\definecolor{mycolor3}{rgb}{0.92900,0.69400,0.12500}%
\definecolor{mycolor4}{rgb}{0.49400,0.18400,0.55600}%
\definecolor{mycolor5}{rgb}{0.46600,0.67400,0.18800}%
\begin{tikzpicture}

\begin{axis}[%
width=0.951\figurewidth,
height=\figureheight,
at={(0\figurewidth,0\figureheight)},
scale only axis,
xmin=2.0000000000,
xmax=16.0000000000,
xlabel style={font=\color{white!15!black}},
xlabel={degrees of freedom $r$},
ymode=log,
ymin=0.0001000000,
ymax=1.0000000000,
yminorticks=true,
ylabel style={font=\color{white!15!black}},
ylabel={relative error},
axis background/.style={fill=white},
legend style={at={(0.03,0.03)}, anchor=south west, legend cell align=left, align=left, fill=none, draw=none}
]
\addplot [color=mycolor5, dashed, mark=o, mark options={solid, mycolor5}]
  table[row sep=crcr]{%
2.0000000000	0.4077750088\\
3.0000000000	0.0137609952\\
4.0000000000	0.0080471133\\
5.0000000000	0.0056120108\\
6.0000000000	0.0011583913\\
7.0000000000	0.0012399846\\
8.0000000000	0.0009060126\\
9.0000000000	0.0007048818\\
10.0000000000	0.0004828383\\
11.0000000000	0.0003876773\\
12.0000000000	0.0003994784\\
13.0000000000	0.0004011380\\
14.0000000000	0.0004228859\\
15.0000000000	0.0004459339\\
};
\addlegendentry{FTR offline}

\addplot [color=mycolor2, mark=asterisk, mark options={solid, mycolor2}]
  table[row sep=crcr]{%
2.0000000000	0.5663801476\\
3.0000000000	0.2733827769\\
4.0000000000	0.0853676468\\
5.0000000000	0.0359364892\\
6.0000000000	0.0392849611\\
7.0000000000	0.0478440027\\
8.0000000000	0.0164932287\\
9.0000000000	0.0163082423\\
10.0000000000	0.0104723429\\
11.0000000000	0.0319796529\\
12.0000000000	0.0224060791\\
13.0000000000	0.0218532213\\
14.0000000000	0.0198855897\\
15.0000000000	0.0177006783\\
};
\addlegendentry{FTR Galerkin}

\addplot [color=mycolor3, dashdotted, mark=x, mark options={solid, mycolor3}]
  table[row sep=crcr]{%
2.0000000000	0.4170981555\\
3.0000000000	0.0154749609\\
4.0000000000	0.0084065328\\
5.0000000000	0.0071120892\\
6.0000000000	0.0013461770\\
7.0000000000	0.0014194921\\
8.0000000000	0.0009391045\\
9.0000000000	0.0010914002\\
10.0000000000	0.0007840160\\
11.0000000000	0.0004526283\\
12.0000000000	0.0005579589\\
13.0000000000	0.0005784812\\
14.0000000000	0.0021263481\\
15.0000000000	0.0016261775\\
};
\addlegendentry{FTR manifold Galerkin}

\addplot [color=mycolor4, dashed, mark=triangle, mark options={solid, rotate=90, mycolor4}]
  table[row sep=crcr]{%
2.0000000000	0.5545657688\\
3.0000000000	0.3690850182\\
4.0000000000	0.2455758503\\
5.0000000000	0.1877715655\\
6.0000000000	0.1610633465\\
7.0000000000	0.1330350338\\
8.0000000000	0.1080029181\\
9.0000000000	0.0856458024\\
10.0000000000	0.0675960963\\
11.0000000000	0.0580961401\\
12.0000000000	0.0485128475\\
13.0000000000	0.0405131675\\
14.0000000000	0.0332140686\\
15.0000000000	0.0261845727\\
};
\addlegendentry{POD offline}

\addplot [color=mycolor1, dotted, mark=triangle, mark options={solid, rotate=270, mycolor1}]
  table[row sep=crcr]{%
2.0000000000	0.4704562431\\
3.0000000000	0.2559992674\\
4.0000000000	0.2462001378\\
5.0000000000	0.3008915205\\
6.0000000000	0.2677913994\\
7.0000000000	0.2900615086\\
8.0000000000	0.2785023298\\
9.0000000000	0.3765875214\\
10.0000000000	0.1763147270\\
11.0000000000	0.2727120667\\
12.0000000000	0.2769410366\\
13.0000000000	0.2930084463\\
14.0000000000	0.2618292069\\
15.0000000000	0.2675921467\\
};
\addlegendentry{POD Galerkin}

\end{axis}
\end{tikzpicture}%
    \caption{Relative error vs. speedup for $r=2,\dots,15$ using the FTR- and POD-Galerkin approach. The error and speedup is measured using the testing data shown in \cref{fig:1D-advection}. The FTR/POD offline errors mark the reconstruction error of the training data, which is collected in the offline phase (snapshots are shown in the upper row of \cref{fig:1D-advection}). All FTR Galerkin results are computed using \cref{eq:MLgalerkin-advect} (snapshots are shown in the lower right of \cref{fig:1D-advection}). Accordingly, the POD Galerkin results are computed using \cref{eq:MLgalerkin-advect0}. The FTR manifold Galerkin results are computed using \cref{eq:MLgalerkin} directly.
    }
    \label{fig:1D-advection-speedup}
\end{figure}

 The relative error and speedup for the test trajectory is shown in the lower part of \cref{fig:1D-advection-speedup} and is plotted together with the POD-Galerkin results. We see that the errors are reduced compared to the results of the POD.
 %%%%%%%%%%%%%%%%%%%%%%%%%%%%%%%
%%%%%%%%%%%%%%%%%%%%%%%%%%%%%%%
 Further details about the simulation are reported in \cref{appx:1d-advection-example}. The results apply similarly to higher spatial dimension $d>1$, for example in case of the moving disk (\cref{subsec:MovingDisc}). Note, that the path simulated in the online phase is limited to areas where the transport field is initialized. These are the areas where any front has traveled during the offline phase as shown in \cref{fig:transportfield}.
 This restriction is, however, shared with classical linear methods.
 Dynamics that are not covered in the ansatz space created from the initial set of snapshots are usually not covered by the ROM.

The success of the heuristic approach \cref{eq:MLgalerkin-advect} is somewhat obvious since the level-set function $\phi$ parameterizes the transport (see \cref{subsec:MovingDisc}), which implies it to be a good basis for the advection operator. The idea to use transport capturing level-set functions to accelerate simulations for advection laws is not new. For example it is intensively used by the characteristic mapping method (CMM) \cite{MercierYinNave2020}, which evolves the initial condition $q_0(\vec{x})$ of a PDE along characteristic curves $\vec{X}(\vec{x},t)$, such that $q(\vec{x},t) =q_0(\vec{X}(\vec{x},t))$. Nevertheless, in \cite{MercierYinNave2020} the authors use an invertible mapping $\vec{X}(\vec{x},t)$, which hinders the applicability for systems with topological changes. Intentionally, this is not done here, since we aim for systems, where topological changes are possible. However, it would be interesting to see if snapshots of the characteristic map can be similarly utilized for MOR as the snapshots of the auxiliary field $\matr{\Phi}$ inside the FTR.

Nevertheless, it should be noted that the procedure proposed for the advection equation cannot be generalized to advection-reaction-diffusion equations and therefore special hyper-reduction methods are needed for an efficient ROM.

%%%%%%%%%%%%%%%%%%%%%%%%%%%%%%%%%%%%%%%%%%%%%%%%%%%%%%%%%%%%%%%%%%%%%
\subsection{Hyper-reduction for Moving Fronts}
\label{subsec:hyperreduction}

Apart from the slow decaying POD approximation errors, advection-reaction-diffusion systems pose another difficulty for model order reduction. The dynamics of advection-reaction-diffusion systems take place at a characteristic length scale $l_f$ defined in \cref{eq:characteristic-length}. This characteristic scale is usually much smaller than the size of the domain or the traveling distance of the front. Hence, the FOM-RHS \cref{eq:ODE-rhs} and its gradient posses only few spatial grid points per time step with non-vanishing support. Therefore, the hyper-reduction methods for nonlinear manifolds \cite{KimChoiWidemannZohdi2021,JainTiso2019} cannot be applied. For example, the extended-ECSW scheme proposed by \cite{JainTiso2019}, or the gappy-POD based GNAT procedure \cite{CarlbergFarhatCortial2013} first introduced for nonlinear manifolds in \cite{KimChoiWidemannZohdi2021} cannot be used here, since they preselect a set of sample points, which is fixed for every time step and all $\mu\in\params$.

In contrast, the FTR-hyper-reduction approach can help to identify the locations of the
front to reduce computational complexity, while sustaining an accurate solution. Here, we propose an idea that is similar to the Reduced Integration Domain (RID) method \cite{Ryckelynck2005} for finite elements.
By imposing a threshold criterion on each finite element, RID is choosing a reduced number of elements to describe a balance condition, i.e.~to minimize the residual between internal and external forces. Similar to RID, we choose a selected number of $M_p$ sampled/selected points to minimize the error of the projected right hand side (i.e.~external/internal force):
\begin{align}
\label{eq:rid-min}
    0=\frac{\d}{\d{\dot{\vec{a}}}}\Vert\matr{J}_g(\vec{a})\dot{\vec{a}} -\RHS(g(\vec{a}),t,\mu)\Vert_{\matr{P}^2_{\vec{a}}}^2
      &= 2 \matr{J}_g(\vec{a})^T\matr{P}^T_{\vec{a}}\matr{P}_{\vec{a}}\matr{J}_g(\vec{a})\dot{\vec{a}}\\ &- 2\matr{J}_g(\vec{a})^T \matr{P}^T_{\vec{a}}\matr{P}_{\vec{a}} \RHS(g(\vec{a}),t,\mu)\,.
\end{align}
Each of the $M_p$ selected sample points corresponds to an index $0\le i\le M$, which is represented as the $i$th standard basis vector $\vec{e}_i\in\mathbb{R}^M$ inside the rows of the selection matrix $\Pa\in \mathbb{R}^{M_p\times M}$. Thus, the hyper-reduced Jacobian and right hand side $\Pa \matr{J}_g,\Pa \RHS$ are only computed at $M_p$ sample points. Note that the stencil size of our finite difference scheme requires to compute $f(\vec{\phi})$ on additional supporting mesh points, contained in $\hat{\matr{P}}_{\vec{a}}\in\mathbb{R}^{\hat{M}_p\times M}$. In practice, $\hat{\matr{P}}_{\vec{a}}f(\vec{\phi}),\Pa \matr{J}_g,\Pa \RHS$ are not computed as matrix products, but as pointwise evaluations of $f(\vec{\phi}),\matr{J}_g,\RHS$ at the corresponding sample points.

\begin{figure}[htp!]
%$\vec{F)(\vec{q},t,\mu)$
    \begin{center}
     \includegraphics[width=1\linewidth]{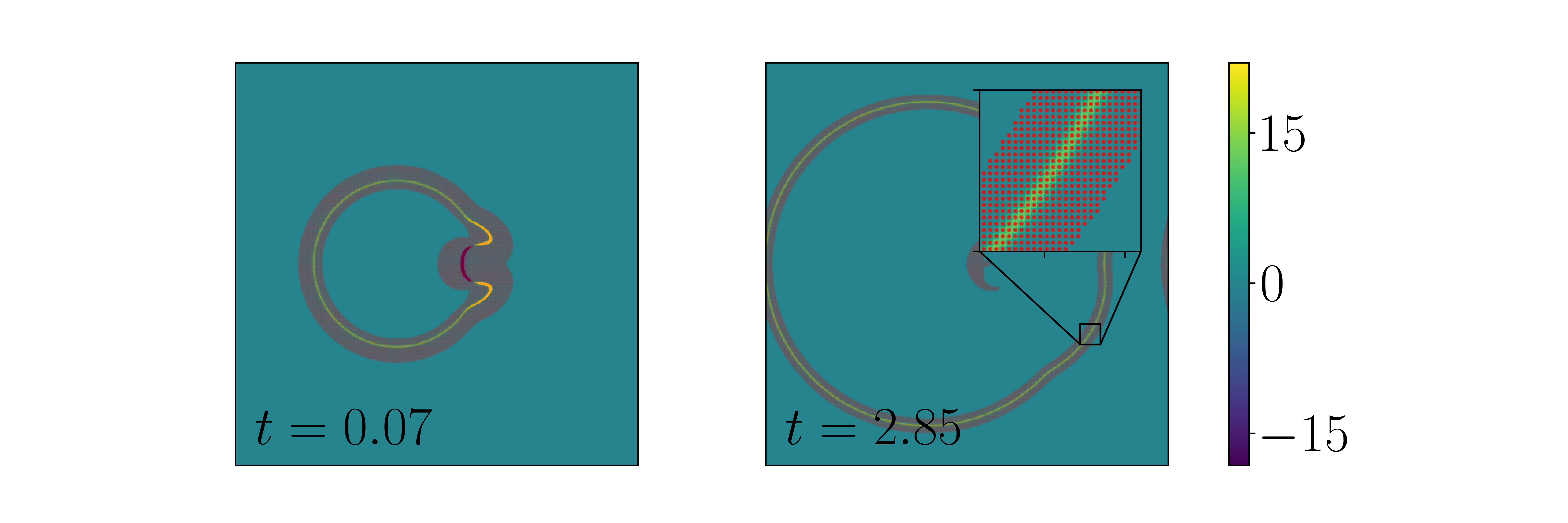}
    \end{center}
    \caption{Color plot of the right hand side $\RHS(\vec{q},t,\mu)$ of the 2D advection-reaction-diffusion system \cref{eq:react-diff-system} for two different time instances $t=0.07$ (left) and $t=2.85$ (right) and the corresponding sample points for a sample fraction of $M_p/M=0.1$. The inset in the right color plot shows a close up of the location of the front.}
    \label{fig:samplemeshpoints}
\end{figure}

In contrast to RID, the selection matrix $\matr{P}_{\vec{a}}\colon \mathbb{R}^r\to\mathbb{R}^{M_p\times M}$ is dependent on the state $\vec{a}(t,\mu)$, which evolves over time (see \cref{fig:samplemeshpoints}).
However, similar to what RID does for external forces, we have to add nodes, i.e.~sample points, at which the right hand side does not vanish.
Since for the FTR $\RHS$ is non-vanishing at the locations of the front, i.e.~at the roots of the level-set function, we can perform a time dependent adaptive thresholding, which defines $\matr{P}_{\vec{a}}$.
The threshold search selects the $M_p$ smallest values of the level-set function $\vec{\phi}=\FTRmodes\vec{a}\in{\mathbb{R}}^M$ at which we evaluate $\hat{\matr{P}}_{\vec{a}}f(\vec{\phi}),\Pa \matr{J}_g,\Pa \RHS$.
For two time instances, the sample points are visualized in \cref{fig:samplemeshpoints} for the 2D ARD-system of \cref{subsec:offline-ARD}.
To reduce the costs of the threshold search, one might recompute the sample points only after a fraction of the characteristic time scale $t_f=l_f/c^*$, where $l_f$ is defined in \cref{eq:characteristic-length}.
Note, that the threshold search is a heuristic in order to perform a cheap  minimization of the residual \cref{eq:rid-min}.

Further, it should be noted that in this work we are using explicit time integration schemes, as they are usually used inside finite difference solvers. Therefore, the aforementioned methods \cite{KimChoiWidemannZohdi2021,JainTiso2019} are not comparable in speedup, since they compare the ROM with implicit time integration schemes used in the offline stage. Nevertheless, applying implicit integration schemes during the online phase may benefit the stability of the resulting ROM.
A promising and efficient method for explicit time integration schemes was proposed by \cite{BlackSchulzeUnger2021} for reaction-diffusion systems in one spatial dimension. Although the framework cannot cope with topological changes, since it relies on a smooth parameterization of the transport, the authors claim speedups of up to a factor of 130.

In the following, we will show some numerical examples utilizing the here presented hyper-reduction approach.

%%%%%%%%%%%%%%%%%%%%%%%%%%%%%%%%%%%%%%%%%%%%%%%%%%%%%%%%%%%%%%%%%%%%%
\subsection{Numerical Examples}
\label{subsec:hyper-ftr-numerical-examples}
In this section, we numerically investigate the applicability of our framework.
Therefore, we define
the offline and online errors:
\begin{equation}
\label{eq:rel-error-test-train}
    \text{offline/online err} = \frac{\norm{\matr{Q}^\text{train/test}-\tilde{\matr{Q}}^\text{train/test}}_\mathrm{F}}{\norm{\matr{Q}^\text{train/test}}_\mathrm{F}}\,.
\end{equation}
Here, $\matr{Q}\in\mathbb{R}^{M\times (N_tN_{\params})}$ is the snapshot matrix containing all snapshots for the $N_t$ time and $N_\params$ parameter instances $\mu \in \params$ in its columns. The superscript "train" ("test") belongs to the snapshots $\mu\in\params^\text{train}$ ($\params^\text{test}$) computed during the offline (online) phase.

The approximation $\tilde{\matr{Q}}^\text{train}$ is therefore either the reconstruction of the training data using the FTR-ansatz (\cref{alg:iterFTR}) or, in case of the POD, the projection onto the first $r$ left singular vectors of $\matr{Q}^\text{train}$ contained in $\matr{U}\in\mathbb{R}^{M\times r}$, i.e.~$\tilde{\matr{Q}}^\text{train}=\matr{U}^T\matr{U}\matr{Q}^\text{train}$.

$\tilde{\matr{Q}}^\text{test}$ refers to the results evaluating the ROM \cref{eq:ROM} for the given time interval and parameters $\mu \in \params^\text{test}$ using the reduced mapping $g:\mathbb{R}^r \to \mathbb{R}^M$.
Specifically, in the case of the POD, the dynamical ROM predictions use $g(\vec{a})=\matr{U}\vec{a}$ as a reduced mapping, whereas $g(\vec{a}) = f(\matr{\Psi}\vec{a})$ for the FTR.

Furthermore, we define the projection error:
\begin{equation}
\label{eq:rel-proj-err}
    \text{proj. err} = \frac{\norm{\matr{Q}^\text{test}-\tilde{\matr{Q}}_*^\text{test}}_\mathrm{F}}{\norm{\matr{Q}^\text{test}}_\mathrm{F}}
\end{equation}
where $\tilde{\matr{Q}}_*^\text{test}$ is the best fit of $\matr{Q}^\text{test}$ with help of our the mapping $g$.

\subsubsection{Reaction-Diffusion System in 1D}
\label{subsubsec:reac-diff-1D}
\begin{figure}[htp!]
    \centering
    \setlength{\figureheight}{0.6\linewidth}
    \setlength{\figurewidth}{0.9\linewidth}
\begin{tikzpicture}

\begin{axis}[%
width=0.282\figurewidth,
height=0.419\figureheight,
at={(0.4\figurewidth,0.581\figureheight)},
scale only axis,
point meta min=0.0000000000,
point meta max=1.0000000000,
axis on top,
xmin=0.5000000000,
xmax=1000.5000000000,
xtick={\empty},
xlabel style={font=\color{white!15!black}},
xlabel={$x$},
ymin=0.5000000000,
ymax=101.5000000000,
ytick={\empty},
ylabel style={font=\color{white!15!black}},
ylabel={$t$},
axis background/.style={fill=white},
title style={font=\bfseries},
title={train data $\delta=1$},
colormap/jet,
colorbar
]
\addplot [forget plot] graphics [xmin=0.5000000000, xmax=1000.5000000000, ymin=0.5000000000, ymax=101.5000000000] {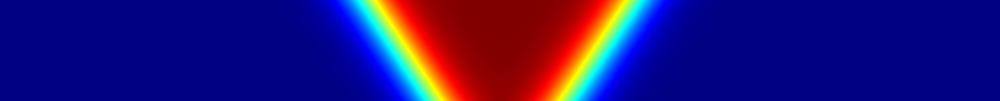};
\end{axis}

\begin{axis}[%
width=0.282\figurewidth,
height=0.419\figureheight,
at={(0.4\figurewidth,0\figureheight)},
scale only axis,
point meta min=0.0000000000,
point meta max=1.0000000000,
axis on top,
xmin=0.5000000000,
xmax=1000.5000000000,
xtick={\empty},
xlabel style={font=\color{white!15!black}},
xlabel={$x$},
ymin=0.5000000000,
ymax=101.5000000000,
ytick={\empty},
ylabel style={font=\color{white!15!black}},
ylabel={$t$},
axis background/.style={fill=white},
title style={font=\bfseries},
title={ROM $\delta_{\mathrm{test}}=0.3$},
colormap/jet,
colorbar
]
\addplot [forget plot] graphics [xmin=0.5000000000, xmax=1000.5000000000, ymin=0.5000000000, ymax=101.5000000000] {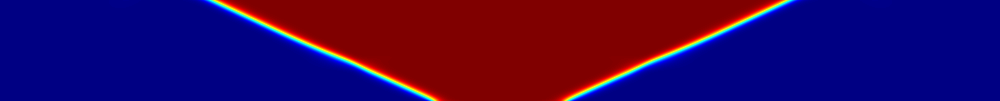};
\end{axis}

\begin{axis}[%
width=0.282\figurewidth,
height=0.419\figureheight,
at={(0\figurewidth,0.581\figureheight)},
scale only axis,
point meta min=-0.0000000000,
point meta max=1.0011460960,
axis on top,
xmin=0.5000000000,
xmax=1000.5000000000,
xtick={\empty},
xlabel style={font=\color{white!15!black}},
xlabel={$x$},
ymin=0.5000000000,
ymax=101.5000000000,
ytick={\empty},
ylabel style={font=\color{white!15!black}},
ylabel={$t$},
axis background/.style={fill=white},
title style={font=\bfseries},
title={train data $\delta=0.2$}
]
\addplot [forget plot] graphics [xmin=0.5000000000, xmax=1000.5000000000, ymin=0.5000000000, ymax=101.5000000000] {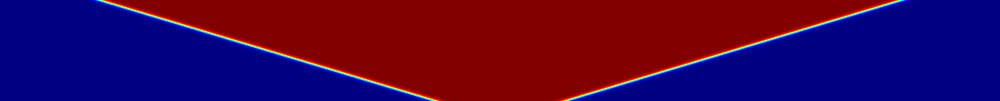};
\end{axis}

\begin{axis}[%
width=0.282\figurewidth,
height=0.419\figureheight,
at={(0\figurewidth,0\figureheight)},
scale only axis,
point meta min=0.0000000000,
point meta max=1.0000000000,
axis on top,
xmin=0.5000000000,
xmax=1000.5000000000,
xtick={\empty},
xlabel style={font=\color{white!15!black}},
xlabel={$x$},
ymin=0.5000000000,
ymax=101.5000000000,
ytick={\empty},
ylabel style={font=\color{white!15!black}},
ylabel={$t$},
axis background/.style={fill=white},
title style={font=\bfseries},
title={FOM $\delta_{\mathrm{test}}=0.3$}
]
\addplot [forget plot] graphics [xmin=0.5000000000, xmax=1000.5000000000, ymin=0.5000000000, ymax=101.5000000000] {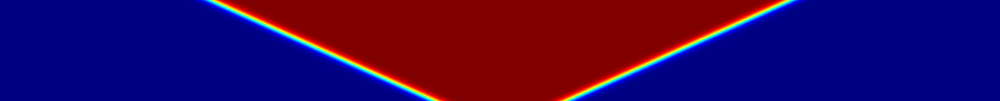};
\end{axis}

\begin{axis}[%
width=1.227\figurewidth,
height=1.227\figureheight,
at={(-0.224\figurewidth,-0.135\figureheight)},
scale only axis,
point meta min=0.0000000000,
point meta max=1.0000000000,
xmin=0.0000000000,
xmax=1.0000000000,
ymin=0.0000000000,
ymax=1.0000000000,
axis line style={draw=none},
ticks=none,
axis x line*=bottom,
axis y line*=left
]
\end{axis}
\end{tikzpicture}%
    \caption{Train and test data of the reaction diffusion system \cref{eq:1Dreact_test} with the FTR-ROM using $r=3$ degrees of freedom.}
    \label{fig:1Dreact_test}
\end{figure}

First, we test our approach on an analytic test case taken and modified from \cite{ZhangZhuLoulaYu2016}.
The test case is based on a one dimensional scalar nonlinear reaction–diffusion equation
\begin{equation}
    \label{eq:1Dreact_test}
    \partial_t q = \partial_{xx} q + \frac{8}{\delta^2} q^2(q-1) \qquad (t,x) \in [0,1]\times[-15,15]
\end{equation}
with corresponding analytical solution
\begin{equation}
    \label{eq:analytical_soluiton_ARD1d}
    q(x,t,\delta)=f\left(\frac{|x|-2t/\delta-2}{\delta}\right),
\end{equation}
given that $f(x)=\sigmoid(2x)$.
We follow the strategy outlined above.
First we compute the numerical solution by discretizing \cref{eq:1Dreact_test} with $M=4000$ grid points and solving it for $\delta\in\params^\text{train}=\{0.2,1\}$ (further details can be found in \cref{appx:1d-advection-example}).
The training data consists of 202 samples, including 101 samples of each training parameter. The training data is visualized as color plot in \cref{fig:1Dreact_test} together with the ROM prediction of the FTR using $r=3$ and $\delta =\delta_\text{test}=0.3$ in \cref{eq:1Dreact_test}.
The FTR algorithm (\cref{alg:iterFTR}) is run for 8000 steps using $\tau = 4$ and different truncation ranks $1<r<10$. After we have computed the reduced mapping $\vec{q}(t,\delta)=f(\FTRmodes\vec{a}(t,\delta))$ from the training set, we can compute the starting values $\vec{a}(0,\delta), \delta\in\params^\text{test}$ to test the ROM \cref{eq:ROM} by minimizing the initial condition of the ROM \cref{eq:inicond_ROM}, using Gauss-Newton iterations \cite{More1978}. As an initial guess for the minimization, we use the set of initial points $\{(\delta,\vec{a}(0,\delta))\mid\delta\in\params^\text{train}\}$ and interpolate them for any given test parameter $\delta\in\params^\text{test}$. Thereafter, the ROM-solution for all test parameters $\delta\in\params^\text{test}$ is compared to the analytical solution \cref{eq:analytical_soluiton_ARD1d}. The results are reported as online errors in \cref{tabl:offline-online-error} together with the offline and projection errors. The online and projection errors are stated for the cumulated snapshots including the time interval $[0,1]$ and all parameters $\delta\in\params^\text{test}=\{0.3,0.4,\dots,0.9\}$. \Cref{tabl:offline-online-error} also compares the results with the POD-Galerkin approach. The starting values for the POD-Galerkin-ROM are simply given by the orthogonal projection of $\vec{q}(0,\delta)$ onto the POD modes. It is remarkable to see that the FTR outperforms the POD by two orders of magnitude.

\begin{figure}
    \centering
    \setlength{\figureheight}{0.4\linewidth}
    \setlength{\figurewidth}{0.9\linewidth}
     % This file was created with tikzplotlib v0.9.17.
\begin{tikzpicture}

\definecolor{color0}{rgb}{0.12156862745098,0.466666666666667,0.705882352941177}
\definecolor{color1}{rgb}{1,0.498039215686275,0.0549019607843137}
\definecolor{color2}{rgb}{0.172549019607843,0.627450980392157,0.172549019607843}
\definecolor{color3}{rgb}{0.83921568627451,0.152941176470588,0.156862745098039}

\begin{axis}[
height=\figureheight,
legend cell align={left},
legend style={
  fill opacity=0.8,
  draw opacity=1,
  text opacity=1,
  at={(0.45,0.97)},
  anchor=north west,
  draw=white!80!black
},
log basis x={10},
log basis y={10},
tick align=outside,
tick pos=left,
width=\figurewidth,
x grid style={white!69.0196078431373!black},
xlabel={cpu-time [s]},
xmin=0.135500564995119, xmax=513.614223389718,
xmode=log,
xtick style={color=black},
xtick={0.01,0.1,1,10,100,1000,10000},
xticklabels={
  \(\displaystyle {10^{-2}}\),
  \(\displaystyle {10^{-1}}\),
  \(\displaystyle {10^{0}}\),
  \(\displaystyle {10^{1}}\),
  \(\displaystyle {10^{2}}\),
  \(\displaystyle {10^{3}}\),
  \(\displaystyle {10^{4}}\)
},
y grid style={white!69.0196078431373!black},
ylabel={rel. error},
ymin=0.00122080924818519, ymax=0.137717326134513,
ymode=log,
ytick style={color=black},
ytick={0.0001,0.001,0.01,0.1,1,10},
yticklabels={
  \(\displaystyle {10^{-4}}\),
  \(\displaystyle {10^{-3}}\),
  \(\displaystyle {10^{-2}}\),
  \(\displaystyle {10^{-1}}\),
  \(\displaystyle {10^{0}}\),
  \(\displaystyle {10^{1}}\)
}
]
\addplot [draw=color0, fill=color0, mark=*, only marks]
table{%
x  y
1.08189429999948 0.014177626700494
3.72935536499972 0.0426528315498698
1.50913509200109 0.00392750832599165
39.6877304900008 0.00621143607809319
151.521196908001 0.112211651195067
225.277721620001 0.0362185548593271
353.157996978 0.0249299974100566
333.361178669 0.0355300396491876
};
\addlegendentry{$M_p/M=0.1$}
\addplot [draw=color1, fill=color1, mark=asterisk, only marks]
table{%
x  y
0.197064820999003 0.0141276044894263
0.468534944000567 0.0213105655103097
0.511958630999288 0.0024368755711217
8.92901356600123 0.00290527502804757
37.6162838189994 0.00224296518656107
48.9977989350009 0.00251210591636464
196.759643026 0.00203800281997098
170.325976438999 0.00186648873793518
};
\addlegendentry{$M_p/M=0.2$}
\addplot [draw=color2, fill=color2, mark=x, only marks]
table{%
x  y
0.22657484299998 0.014156365640191
0.488315833001252 0.0213248807914035
0.443433594999078 0.00268041799505056
8.86156603999916 0.00317332244572266
39.7402081060009 0.00247984590902387
50.8256818450009 0.0025840333575287
189.857672704 0.00208615340628652
195.084876711 0.00190761017751328
};
\addlegendentry{$M_p/M=0.5$}
\addplot [draw=color3, fill=color3, mark=+, only marks]
table{%
x  y
0.327137853999375 0.0141571296707045
0.686877934000222 0.0213244887487675
0.664318199000263 0.0026849601330949
13.4508918069987 0.00317671610638107
53.2804623939992 0.00248035057234887
67.0116572049992 0.00258455304453938
263.330392292 0.00208778213247801
269.059476478 0.00190815963511283
};
\addlegendentry{$M_p/M=1.0$}
\addplot [semithick, black, dashed, forget plot]
table {%
72.7903160679989 0.00152080924818519
72.7903160679989 0.137717326134513
};
\draw (axis cs:1.08189429999948,0.014177626700494) node[
  scale=0.7,
  anchor=north,
  text=black,
  rotate=0.0
]{(2)};
\draw (axis cs:3.72935536499972,0.0426528315498698) node[
  scale=0.7,
  anchor=south,
  text=black,
  rotate=0.0
]{(3)};
\draw (axis cs:1.50913509200109,0.00392750832599165) node[
  scale=0.7,
  anchor=south,
  text=black,
  rotate=0.0
]{(4)};
\draw (axis cs:39.6877304900008,0.00621143607809319) node[
  scale=0.7,
  anchor=south,
  text=black,
  rotate=0.0
]{(5)};
\draw (axis cs:151.521196908001,0.112211651195067) node[
  scale=0.7,
  anchor=north,
  text=black,
  rotate=0.0
]{(6)};
\draw (axis cs:225.277721620001,0.0362185548593271) node[
  scale=0.7,
  anchor=south,
  text=black,
  rotate=0.0
]{(7)};
\draw (axis cs:353.157996978,0.0249299974100566) node[
  scale=0.7,
  anchor=north,
  text=black,
  rotate=0.0
]{(8)};
\draw (axis cs:333.361178669,0.0355300396491876) node[
  scale=0.7,
  anchor=south,
  text=black,
  rotate=0.0
]{(9)};
\draw (axis cs:0.197064820999003,0.0141276044894263) node[
  scale=0.7,
  anchor=north,
  text=black,
  rotate=0.0
]{(2)};
\draw (axis cs:0.468534944000567,0.0213105655103097) node[
  scale=0.7,
  anchor=north,
  text=black,
  rotate=0.0
]{(3)};
\draw (axis cs:0.511958630999288,0.0024368755711217) node[
  scale=0.7,
  anchor=north,
  text=black,
  rotate=0.0
]{(4)};
\draw (axis cs:8.92901356600123,0.00290527502804757) node[
  scale=0.7,
  anchor=north,
  text=black,
  rotate=0.0
]{(5)};
\draw (axis cs:37.6162838189994,0.00224296518656107) node[
  scale=0.7,
  anchor=north,
  text=black,
  rotate=0.0
]{(6)};
\draw (axis cs:48.9977989350009,0.00251210591636464) node[
  scale=0.7,
  anchor=north,
  text=black,
  rotate=0.0
]{(7)};
\draw (axis cs:196.759643026,0.00203800281997098) node[
  scale=0.7,
  anchor=south,
  text=black,
  rotate=0.0
]{(8)};
\draw (axis cs:170.325976438999,0.00186648873793518) node[
  scale=0.7,
  anchor=south,
  text=black,
  rotate=0.0
]{(9)};
\draw (axis cs:0.22657484299998,0.014156365640191) node[
  scale=0.7,
  anchor=south,
  text=black,
  rotate=0.0
]{(2)};
\draw (axis cs:0.488315833001252,0.0213248807914035) node[
  scale=0.7,
  anchor=south,
  text=black,
  rotate=0.0
]{(3)};
\draw (axis cs:0.443433594999078,0.00268041799505056) node[
  scale=0.7,
  anchor=south,
  text=black,
  rotate=0.0
]{(4)};
\draw (axis cs:8.86156603999916,0.00317332244572266) node[
  scale=0.7,
  anchor=south,
  text=black,
  rotate=0.0
]{(5)};
\draw (axis cs:39.7402081060009,0.00247984590902387) node[
  scale=0.7,
  anchor=south,
  text=black,
  rotate=0.0
]{(6)};
\draw (axis cs:50.8256818450009,0.0025840333575287) node[
  scale=0.7,
  anchor=south,
  text=black,
  rotate=0.0
]{(7)};
% \draw (axis cs:189.857672704,0.00208615340628652) node[
%   scale=0.7,
%   anchor=north,
%   text=black,
%   rotate=0.0
% ]{(8)};
\draw (axis cs:195.084876711,0.00190761017751328) node[
  scale=0.7,
  anchor=north,
  text=black,
  rotate=0.0
]{(9)};
\draw (axis cs:0.327137853999375,0.0141571296707045) node[
  scale=0.7,
  anchor=north,
  text=black,
  rotate=0.0
]{(2)};
\draw (axis cs:0.686877934000222,0.0213244887487675) node[
  scale=0.7,
  anchor=north,
  text=black,
  rotate=0.0
]{(3)};
\draw (axis cs:0.664318199000263,0.0026849601330949) node[
  scale=0.7,
  anchor=north,
  text=black,
  rotate=0.0
]{(4)};
\draw (axis cs:13.4508918069987,0.00317671610638107) node[
  scale=0.7,
  anchor=north,
  text=black,
  rotate=0.0
]{(5)};
\draw (axis cs:53.2804623939992,0.00248035057234887) node[
  scale=0.7,
  anchor=north,
  text=black,
  rotate=0.0
]{(6)};
\draw (axis cs:67.0116572049992,0.00258455304453938) node[
  scale=0.7,
  anchor=north,
  text=black,
  rotate=0.0
]{(7)};
\draw (axis cs:263.330392292,0.00208778213247801) node[
  scale=0.7,
  anchor=south,
  text=black,
  rotate=0.0
]{(8)};
\draw (axis cs:269.059476478,0.00190815963511283) node[
  scale=0.7,
  anchor=north,
  text=black,
  rotate=0.0
]{(9)};
\draw (axis cs:72.7903160679989,0.117728909317924) node[
  scale=0.8,
  anchor=north east,
  text=black,
  rotate=0.0
]{FOM};
\end{axis}

\end{tikzpicture}
    \caption{Error vs. CPU-time for the accumulated parameter range $\delta\in\params^\text{test}$. Different ranks $r$ are indicated as $(r)$ above or below the markers. The dashed line indicates the CPU-time needed for solving the FOM. The sampled fraction $M_p/M$ in the hyper-reduced ROM is given in terms of the size $M$ of the FOM.}
    \label{fig:cpu_time_vs_error}
\end{figure}
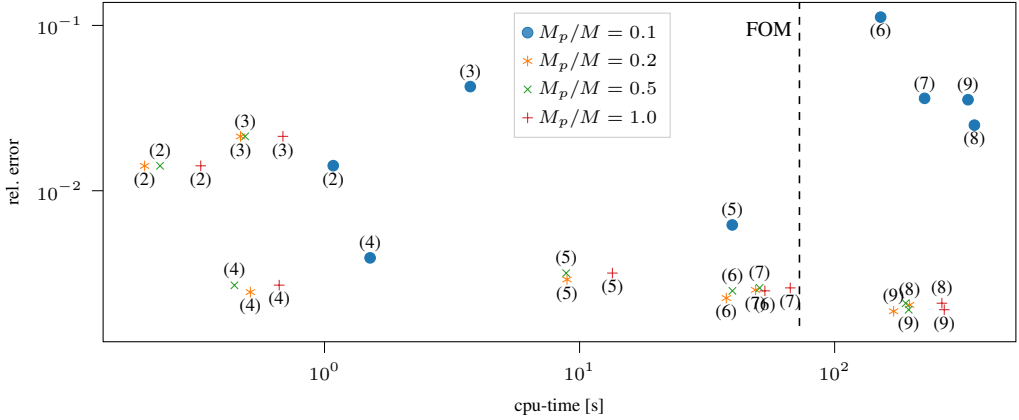

Next, we are interested in whether the gain in precision can be translated to speedups. Therefore, we study the performance of the hyper-FTR explained in \cref{subsec:hyperreduction}. \Cref{fig:cpu_time_vs_error} compares CPU-time and error for $M_p=0.1M,0.2M,0.5M$ and $M$ number of grid points, where $M$ is the dimension of the FOM.
The figure indicates that even without hyper-reduction, speedups can be achieved compared to the FOM, due to larger step sizes in the reduced coordinates. Comparing the hyper-FTR with a sample fraction of $M_p/M=0.2$ to 1 we see another speedup in CPU-time. For a reduction below $0.1M$ grid points, the solution is unstable and can lead to additional time steps, making the overall simulation slower.

%The hyper-FTR is compared with the classical POD-DEIM hyper-reduction approach. Here we make use of the discrete empirical interpolation method (DEIM) \cite{ChaturantabutSorensen2010} implemented as Q-DEIM \cite{ZlatkoGugercin2016}.

\begin{table}[htp!]
    \centering
 \begin{tabular}{rrrrrr}
\toprule
&\multicolumn{3}{c}{FTR} & \multicolumn{2}{c}{POD} \\
    \cmidrule(lr){2-4} \cmidrule(lr){5-6}
   rank &   offine error &   online error &   proj. error &   online error &   proj. error \\
\midrule
      2 &            8.2e-03 &            1.4e-02 &                3.0e-03 &            3.6e-01 &                2.7e-01 \\
      3 &            2.6e-03 &            2.1e-02 &                6.6e-03 &            2.8e-01 &                2.0e-01 \\
      4 &            6.2e-04 &            2.7e-03 &                5.3e-04 &            2.4e-01 &                1.4e-01 \\
      5 &            5.3e-04 &            3.2e-03 &                7.2e-04 &            2.3e-01 &                1.1e-01 \\
      6 &            5.4e-04 &            2.5e-03 &                3.7e-04 &            2.5e-01 &                9.0e-02 \\
      7 &            5.0e-04 &            2.6e-03 &                2.7e-04 &            3.4e-01 &                7.3e-02 \\
      8 &            4.4e-04 &            2.1e-03 &                1.6e-04 &            2.7e-01 &                6.0e-02 \\
      9 &            2.0e-04 &            1.9e-03 &                2.2e-04 &            2.0e-01 &                5.0e-02 \\
\bottomrule
\end{tabular}

    \caption{Offline, online and projection errors for FTR and POD. The errors are reported for the cumulated snapshot data of the train and test parameters used in \cref{subsubsec:reac-diff-1D}.}
    \label{tabl:offline-online-error}
\end{table}

\subsubsection{Advection-Reaction-Diffusion System in 2D}
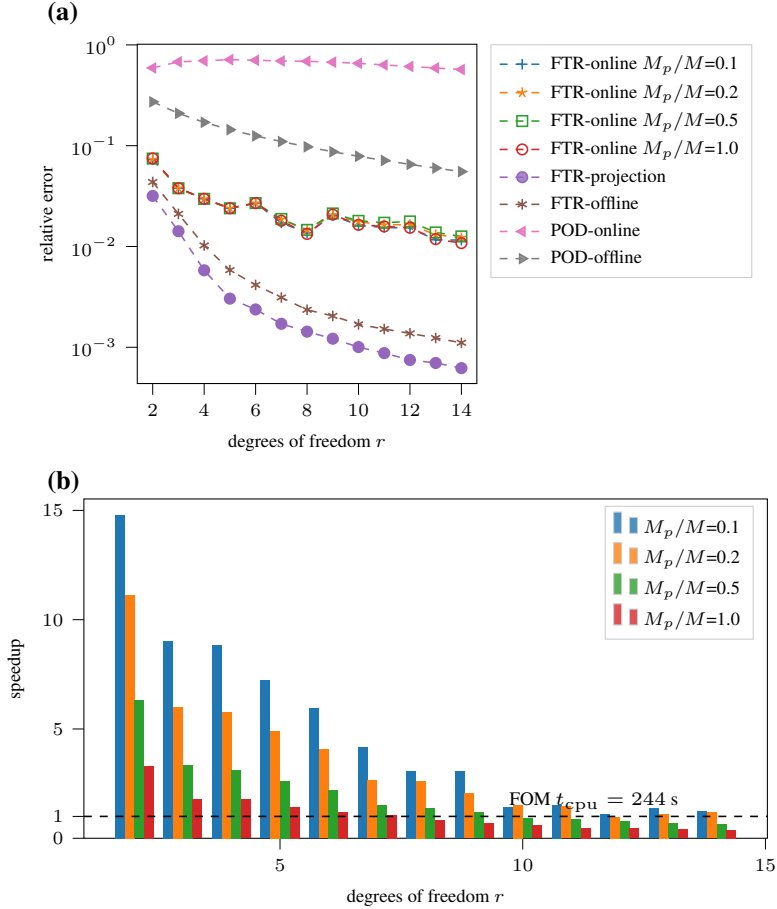
\begin{figure}

    \setlength{\figureheight}{0.4\linewidth}
    \setlength{\figurewidth}{0.4\linewidth}
    \begin{flushleft}
    \hspace*{3cm}\textbf{(a)}
    \end{flushleft}
    \centering
    % This file was created with tikzplotlib v0.9.17.
\begin{tikzpicture}

\definecolor{color0}{rgb}{0.12156862745098,0.466666666666667,0.705882352941177}
\definecolor{color1}{rgb}{1,0.498039215686275,0.0549019607843137}
\definecolor{color2}{rgb}{0.172549019607843,0.627450980392157,0.172549019607843}
\definecolor{color3}{rgb}{0.83921568627451,0.152941176470588,0.156862745098039}
\definecolor{color4}{rgb}{0.580392156862745,0.403921568627451,0.741176470588235}
\definecolor{color5}{rgb}{0.549019607843137,0.337254901960784,0.294117647058824}
\definecolor{color6}{rgb}{0.890196078431372,0.466666666666667,0.76078431372549}

\begin{axis}[
height=\figureheight,
legend cell align={left},
legend style={
  fill opacity=0.8,
  draw opacity=1,
  text opacity=1,
  at={(1.04,1)},
  anchor=north west,
  draw=white!80!black
},
log basis y={10},
tick align=outside,
tick pos=left,
width=\figurewidth,
x grid style={white!69.0196078431373!black},
xlabel={degrees of freedom \(\displaystyle r\)},
xmin=1.4, xmax=14.6,
xtick style={color=black},
xtick={0,2,4,6,8,10,12,14,16},
xticklabels={
  \(\displaystyle {0}\),
  \(\displaystyle {2}\),
  \(\displaystyle {4}\),
  \(\displaystyle {6}\),
  \(\displaystyle {8}\),
  \(\displaystyle {10}\),
  \(\displaystyle {12}\),
  \(\displaystyle {14}\),
  \(\displaystyle {16}\)
},
y grid style={white!69.0196078431373!black},
ylabel={relative error},
ymin=0.000436609659130152, ymax=1.01511401957085,
ymode=log,
ytick style={color=black},
ytick={1e-05,0.0001,0.001,0.01,0.1,1,10,100},
yticklabels={
  \(\displaystyle {10^{-5}}\),
  \(\displaystyle {10^{-4}}\),
  \(\displaystyle {10^{-3}}\),
  \(\displaystyle {10^{-2}}\),
  \(\displaystyle {10^{-1}}\),
  \(\displaystyle {10^{0}}\),
  \(\displaystyle {10^{1}}\),
  \(\displaystyle {10^{2}}\)
}
]
\addplot [semithick, color0, dashed, mark=+, mark options={solid,fill opacity=0}]
table {%
2 0.0725723064519199
3 0.0377989509532122
4 0.0297196637864113
5 0.0239363034542761
6 0.0269895837200707
7 0.0171509087901726
8 0.0136855586843327
9 0.0203302180960742
10 0.0162189681271827
11 0.0154393272037413
12 0.0152418868069443
13 0.0116484098961696
14 0.0115521363313882
};
\addlegendentry{FTR-online $M_p/M$=0.1}
\addplot [semithick, color1, dashed, mark=star, mark options={solid,fill opacity=0}]
table {%
2 0.0730251137895063
3 0.0378182768332661
4 0.0297568903850871
5 0.0239222963868737
6 0.0269696074885071
7 0.0185770614632445
8 0.014502622327415
9 0.0205200062311621
10 0.0175684071582196
11 0.0164013559145669
12 0.0164785591146193
13 0.0130977260936946
14 0.012205733136347
};
\addlegendentry{FTR-online $M_p/M$=0.2}
\addplot [semithick, color2, dashed, mark=square, mark options={solid,fill opacity=0}]
table {%
2 0.0745352662898091
3 0.037839762157362
4 0.0296986320018274
5 0.0239366008543418
6 0.0269658356668252
7 0.0187525922977997
8 0.0146101335805107
9 0.0212614170326111
10 0.017949883070751
11 0.0172068465030236
12 0.017717451109012
13 0.0137779037429554
14 0.0125841591002773
};
\addlegendentry{FTR-online $M_p/M$=0.5}
\addplot [semithick, color3, dashed, mark=o, mark options={solid,fill opacity=0}]
table {%
2 0.0747672510786933
3 0.0378530163849239
4 0.0297195703471032
5 0.02395117726898
6 0.0269442269210057
7 0.0180121533451234
8 0.0133141927345495
9 0.0207306450272601
10 0.0163857594696006
11 0.0158085190777073
12 0.0153292960190455
13 0.0118168761979525
14 0.0107967735223357
};
\addlegendentry{FTR-online $M_p/M$=1.0}
\addplot [semithick, color4, dashed, mark=*, mark options={solid}]
table {%
2 0.0317102022337074
3 0.0141696744954932
4 0.00580286346319023
5 0.00303717627507393
6 0.00237103976058244
7 0.00171075999599629
8 0.00142838521346802
9 0.0012161472583887
10 0.00100416408046534
11 0.000875368399894198
12 0.000748158333982055
13 0.000698577007109573
14 0.000621029882432832
};
\addlegendentry{FTR-projection}
\addplot [semithick, color5, dashed, mark=asterisk, mark options={solid}]
table {%
2 0.0435949708459934
3 0.0210550984415973
4 0.0102376039583448
5 0.00582464704341502
6 0.00415523729012852
7 0.00311898742820736
8 0.00235326006712577
9 0.00204555793850509
10 0.00168827461000907
11 0.00151660779563982
12 0.00137457004978561
13 0.00123618280643639
14 0.00110737249465049
};
\addlegendentry{FTR-offline}
\addplot [semithick, color6, dashed, mark=triangle*, mark options={solid,rotate=90}]
table {%
2 0.590897185811785
3 0.67896241713495
4 0.695576043797145
5 0.713667085272665
6 0.703881691924275
7 0.692460810430824
8 0.687740188403607
9 0.673991869779725
10 0.656550913410807
11 0.632335771622905
12 0.61071364496498
13 0.590256084277138
14 0.56954953273585
};
\addlegendentry{POD-online}
\addplot [semithick, white!49.8039215686275!black, dashed, mark=triangle*, mark options={solid,rotate=270}]
table {%
2 0.272257361505081
3 0.209411111819684
4 0.170471178998648
5 0.144166018047243
6 0.124858556800937
7 0.10989168803553
8 0.0974784509522021
9 0.0871045999061421
10 0.0784831646561915
11 0.071289428967054
12 0.0652432063956812
13 0.0599406048183806
14 0.0553150855478573
};
\addlegendentry{POD-offline}
\end{axis}

\end{tikzpicture}
    %%%)
    \setlength{\figurewidth}{0.7\linewidth}
    \begin{flushleft}
    \hspace*{3cm}\textbf{(b)}
    \end{flushleft}
    % This file was created with tikzplotlib v0.9.15.
\begin{tikzpicture}

\definecolor{color0}{rgb}{0.12156862745098,0.466666666666667,0.705882352941177}
\definecolor{color1}{rgb}{1,0.498039215686275,0.0549019607843137}
\definecolor{color2}{rgb}{0.172549019607843,0.627450980392157,0.172549019607843}
\definecolor{color3}{rgb}{0.83921568627451,0.152941176470588,0.156862745098039}

\begin{axis}[
height=\figureheight,
legend cell align={left},
legend style={fill opacity=0.8, draw opacity=1, text opacity=1, draw=white!80!black},
tick align=outside,
tick pos=left,
width=\figurewidth,
x grid style={white!69.0196078431373!black},
xlabel={degrees of freedom \(\displaystyle r\)},
xmin=0.96, xmax=15.04,
xtick style={color=black},
xtick={0,5,10,15,20},
xticklabels={
  \(\displaystyle {0}\),
  \(\displaystyle {5}\),
  \(\displaystyle {10}\),
  \(\displaystyle {15}\),
  \(\displaystyle {20}\)
},
y grid style={white!69.0196078431373!black},
ylabel={speedup},
ymin=0, ymax=15.5224993011298,
ytick style={color=black},
ytick={0,1,5,10,15},
yticklabels={
  \(\displaystyle {0}\),
  \(\displaystyle {1}\),
  \(\displaystyle {5}\),
  \(\displaystyle {10}\),
  \(\displaystyle {15}\)
}
]
\draw[draw=none,fill=color0] (axis cs:1.6,0) rectangle (axis cs:1.8,14.7833326677427);
\addlegendimage{ybar,ybar legend,draw=none,fill=color0}
\addlegendentry{$M_p/M$=0.1}

\draw[draw=none,fill=color0] (axis cs:2.6,0) rectangle (axis cs:2.8,9.01808226038621);
\draw[draw=none,fill=color0] (axis cs:3.6,0) rectangle (axis cs:3.8,8.83752713011157);
\draw[draw=none,fill=color0] (axis cs:4.6,0) rectangle (axis cs:4.8,7.21303510189897);
\draw[draw=none,fill=color0] (axis cs:5.6,0) rectangle (axis cs:5.8,5.94681171786119);
\draw[draw=none,fill=color0] (axis cs:6.6,0) rectangle (axis cs:6.8,4.14463024177576);
\draw[draw=none,fill=color0] (axis cs:7.6,0) rectangle (axis cs:7.8,3.07445630101301);
\draw[draw=none,fill=color0] (axis cs:8.6,0) rectangle (axis cs:8.8,3.05974577199434);
\draw[draw=none,fill=color0] (axis cs:9.6,0) rectangle (axis cs:9.8,1.41506446247534);
\draw[draw=none,fill=color0] (axis cs:10.6,0) rectangle (axis cs:10.8,1.51209373015759);
\draw[draw=none,fill=color0] (axis cs:11.6,0) rectangle (axis cs:11.8,1.07980889800429);
\draw[draw=none,fill=color0] (axis cs:12.6,0) rectangle (axis cs:12.8,1.35830998024488);
\draw[draw=none,fill=color0] (axis cs:13.6,0) rectangle (axis cs:13.8,1.22990172013161);
\draw[draw=none,fill=color1] (axis cs:1.8,0) rectangle (axis cs:2,11.1266985236303);
\addlegendimage{ybar,ybar legend,draw=none,fill=color1}
\addlegendentry{$M_p/M$=0.2}

\draw[draw=none,fill=color1] (axis cs:2.8,0) rectangle (axis cs:3,6.00688552584982);
\draw[draw=none,fill=color1] (axis cs:3.8,0) rectangle (axis cs:4,5.78435367766417);
\draw[draw=none,fill=color1] (axis cs:4.8,0) rectangle (axis cs:5,4.91071727909166);
\draw[draw=none,fill=color1] (axis cs:5.8,0) rectangle (axis cs:6,4.0708740612491);
\draw[draw=none,fill=color1] (axis cs:6.8,0) rectangle (axis cs:7,2.65131973540786);
\draw[draw=none,fill=color1] (axis cs:7.8,0) rectangle (axis cs:8,2.61324531541825);
\draw[draw=none,fill=color1] (axis cs:8.8,0) rectangle (axis cs:9,2.05102690189419);
\draw[draw=none,fill=color1] (axis cs:9.8,0) rectangle (axis cs:10,1.53083406547589);
\draw[draw=none,fill=color1] (axis cs:10.8,0) rectangle (axis cs:11,1.44839396598739);
\draw[draw=none,fill=color1] (axis cs:11.8,0) rectangle (axis cs:12,0.974062182725757);
\draw[draw=none,fill=color1] (axis cs:12.8,0) rectangle (axis cs:13,1.10869477190115);
\draw[draw=none,fill=color1] (axis cs:13.8,0) rectangle (axis cs:14,1.18342706031466);
\draw[draw=none,fill=color2] (axis cs:2,0) rectangle (axis cs:2.2,6.30654312237885);
\addlegendimage{ybar,ybar legend,draw=none,fill=color2}
\addlegendentry{$M_p/M$=0.5}

\draw[draw=none,fill=color2] (axis cs:3,0) rectangle (axis cs:3.2,3.35140055498753);
\draw[draw=none,fill=color2] (axis cs:4,0) rectangle (axis cs:4.2,3.13576160459798);
\draw[draw=none,fill=color2] (axis cs:5,0) rectangle (axis cs:5.2,2.60894305414811);
\draw[draw=none,fill=color2] (axis cs:6,0) rectangle (axis cs:6.2,2.20350834230044);
\draw[draw=none,fill=color2] (axis cs:7,0) rectangle (axis cs:7.2,1.50611807716919);
\draw[draw=none,fill=color2] (axis cs:8,0) rectangle (axis cs:8.2,1.38130034642053);
\draw[draw=none,fill=color2] (axis cs:9,0) rectangle (axis cs:9.2,1.17569530323028);
\draw[draw=none,fill=color2] (axis cs:10,0) rectangle (axis cs:10.2,0.901583011873028);
\draw[draw=none,fill=color2] (axis cs:11,0) rectangle (axis cs:11.2,0.852645502509643);
\draw[draw=none,fill=color2] (axis cs:12,0) rectangle (axis cs:12.2,0.797810821428296);
\draw[draw=none,fill=color2] (axis cs:13,0) rectangle (axis cs:13.2,0.679764884833134);
\draw[draw=none,fill=color2] (axis cs:14,0) rectangle (axis cs:14.2,0.658605853253986);
\draw[draw=none,fill=color3] (axis cs:2.2,0) rectangle (axis cs:2.4,3.29467101844845);
\addlegendimage{ybar,ybar legend,draw=none,fill=color3}
\addlegendentry{$M_p/M$=1.0}

\draw[draw=none,fill=color3] (axis cs:3.2,0) rectangle (axis cs:3.4,1.77004452139332);
\draw[draw=none,fill=color3] (axis cs:4.2,0) rectangle (axis cs:4.4,1.7637574957223);
\draw[draw=none,fill=color3] (axis cs:5.2,0) rectangle (axis cs:5.4,1.43151128161804);
\draw[draw=none,fill=color3] (axis cs:6.2,0) rectangle (axis cs:6.4,1.1976310340329);
\draw[draw=none,fill=color3] (axis cs:7.2,0) rectangle (axis cs:7.4,1.06124557567399);
\draw[draw=none,fill=color3] (axis cs:8.2,0) rectangle (axis cs:8.4,0.811142798144947);
\draw[draw=none,fill=color3] (axis cs:9.2,0) rectangle (axis cs:9.4,0.684040017909166);
\draw[draw=none,fill=color3] (axis cs:10.2,0) rectangle (axis cs:10.4,0.580233454081893);
\draw[draw=none,fill=color3] (axis cs:11.2,0) rectangle (axis cs:11.4,0.473223154445031);
\draw[draw=none,fill=color3] (axis cs:12.2,0) rectangle (axis cs:12.4,0.454714577206314);
\draw[draw=none,fill=color3] (axis cs:13.2,0) rectangle (axis cs:13.4,0.396805563685674);
\draw[draw=none,fill=color3] (axis cs:14.2,0) rectangle (axis cs:14.4,0.376583279187956);
\addplot [semithick, black, dashed, forget plot]
table {%
0.96 1
15.04 1
};
\draw (axis cs:9.5,1.5) node[
  scale=1.2,
  anchor=base west,
  text=black,
  rotate=0.0
]{\tiny FOM $t_{\mathrm{cpu}}=244$~s};
\end{axis}

\end{tikzpicture}
    %%%
    \caption{Hyper-FTR results: (a) relative errors defined in \cref{eq:rel-error-test-train,eq:rel-proj-err} for $\params^\text{test/train}$ using POD and FTR decomposition.  (b) Speedup vs. degrees of freedom for the cumulated parameter range $\gamma\in\params^\text{test}$. The speedups are compared for different numbers of sampled grid points $M_p\le M$. The full order model (FOM) using $M=512^2$ grid points is marked with a dashed line.}
    \label{fig:pacman-results}
\end{figure}

Finally, we test the online performance of the hyper-FTR on the advection-reaction-diffusion example of \cref{eq:react-diff-system}, introduced in \cref{subsec:offline-ARD}. We build the training/testing data $\matr{Q}^\text{train}$ from 101 equally spaced snapshots (visualized in \cref{fig:pacman-offline}) with $t\in[0,3]$, $\gamma\in \params^\text{train}=\{10,30,50,70,100\}$ and respectively $\gamma\in \params^\text{test}=\{20,40,60,80,90\}$.
The online, offline, and projection errors of the test case are shown in \cref{fig:pacman-results} together with the speedup generated by the hyper-reduction scheme. It is visible that the FTR outperforms the POD with respect to the offline and online errors. Furthermore, the utilized hyper-reduction strategy results in speedups with moderate online errors. Note that reducing the integration domain to about $10\%$ of its original size (see sample points in \cref{fig:samplemeshpoints}) does not affect the online error, as can be seen from \cref{fig:pacman-results} (a). \Cref{fig:pacman-results} (b) shows, that for small $r$, the additional costs ($\mathcal{O}(rM)$) for the matrix multiplication $\vec{\phi}(t,\mu)=\FTRmodes\vec{a}(t,\mu)$ are negligible, compared to the evaluation of $\RHS$. However, as soon as $r$ becomes large, the speedups of the hyper-reduction scheme are compensated by the computation of $\vec{\phi}$ inside the threshold search. The balance point at which the additional costs compensate the costs of the RHS is problem dependent, but computing the sample points from $\vec{\phi}$ is a bottleneck of this method.
Nevertheless, when aiming for more complex examples like combustion systems or 3D ARD systems, the outlined hyper-reduction approach will benefit from a more computationally complex RHS, which will shift the balance point towards a higher number of modes.

\section{Discussion and Conclusion}
\label{sec:conclusion_and_outlook}

In this work, we have introduced the \textit{front transport reduction} (FTR) method to decompose and simulate transports of complex moving fronts. The decomposition parameterizes moving fronts with the help of a transport-dependent auxiliary field $\vec{\phi}$ and a function $f$ to approximate the front profile. 
Two different decomposition algorithms have been proposed based on singular value thresholding (\cref{alg:iterFTR}) and artificial neural networks (\cref{subsec:ML}). These methods are purely data-driven since they only require a set of snapshots $\q(t,\mu)\in\mathbb{R}^M$ of the FOM as input. 
The resulting approximation $\q(t,\mu)\approx f(\vec{\phi}(t,\mu))$ is well suited for model order reduction of reacting fronts, since $\vec{\phi}(t,\mu)=\matr{U}\vec{a}(t,\mu)\in\mathbb{R}^M$ can be represented by a few $r\ll M$ spatial modes collected in $\matr{U}\in\mathbb{R}^{M\times r}$. 

We emphasize that the utilized front-structure is inherent for advection-reaction-diffusion (ARD) systems (see for example \cite{HadelerRothe1975,BerestyckiHamelNadirashvili2004,BerestyckiHamelRoques2004,Fischer1937}). Making explicit use of the physical structure has advantages over other linear and nonlinear dimension reduction methods, for reasons we discuss in the following:
It was shown, that for various ARD systems the FTR requires fewer modes to decompose the input snapshots compared with the proper orthogonal decomposition (POD), i.e.~it has a better compression quality.
Regarding artificial autoencoder networks, the FTR is similar in the sense that it uses a linear layer activated by a problem dependent nonlinear front function as a decoder. Here, other authors \cite{KimChoiWidemannZohdi2021} use multiple nonlinear activated layers $\vec{q}\approx f(f\dots f(\vec{a}))$, resulting in costly evaluations of the network itself. This can limit the overall performance of the ROM when evaluating the additional nonlinearities. Furthermore, the autoencoder networks are often difficult to tune and require training on GPUs. 
Similar to artificial autoencoder networks, the FTR can approximate topological changes in the evolution of the contour line of the front since it does not make explicit assumptions on the mapping. Here, methods like the shifted POD, previously applied to similar problems in \cite{BlackSchulzeUnger2021}, cannot be used, since they assume one-to-one mappings to align the front.

The ability of the FTR to predict new system states has been demonstrated for non-intrusive (\cref{subsec:data-driven}) and intrusive (\cref{subsec:manifoldgalerkin}) ROMs.
Since the FTR gives additional insights into the underlying structure (transport field $\vec{\phi}$), it allows us to use this information when predicting new system states.
As an example, we heuristically reduced the integration domain during the online evaluation of the Galerkin projected ODE system, using the knowledge of $\vec{\phi}$. This can be seen as an adaptive version of the reduced integration domain method \cite{Ryckelynck2005}. Other nonlinear hyper-reduction methods preselect a set of sample points, on which the dynamics are evaluated. Since for the studied systems, only sample points close to the front are relevant for the dynamics, such hyper-reduction methods may fail.
Although the outlined hyper-reduction procedure yields speedups in CPU-time, it needs a substantially larger number of sample points $M_p$ than required by the dimensions of the ROM $r\ll M_p\ll M$. Therefore, the construction of more efficient hyper-reduction schemes is left open for future research.

To apply our findings to more complex advection-reaction-diffusion systems such as combustion systems in fluid mechanics with multiple reacting species, the decomposition has to be extended to allow arbitrary traveling front shapes. Here, the FTR method would benefit from a generalization or combination with the shifted POD \cite{Reiss2021}, as this would allow to decompose multiple traveling wave systems with topological changes. Furthermore, it would be interesting to see if our approach can be applied to multi-phase flows, as they inherit a similar front structure separating the fluids. Here, the similarity with level-set-based methods like the characteristic mapping method \cite{MercierYinNave2020} should be exploited.

\addcontentsline{toc}{section}{Code Availability}
\section*{Code Availability}
For some selected examples we provide MATLAB and Python code at:
\vspace{0.1cm}
\begin{center}
\faGithub \quad \url{https://github.com/Philipp137/FrontTransportReduction}
\end{center}

\addcontentsline{toc}{section}{Acknowledgement}
\section*{Acknowledgement}
Philipp Krah gratefully acknowledges the support of the Deutsche Forschungsgemeinschaft (DFG) as part of GRK2433 DAEDALUS. Furthermore, we thank Iris Wohnsiedler, Kai Schneider, Mathias Lemke and Volker Mehrmann for useful comments and discussions.

\addcontentsline{toc}{section}{Declaration of competing interest}
\section*{Declaration of competing interest}
The authors declare that they have no known competing financial interests or personal relationships that could have appeared to influence the work reported in this paper.

\addcontentsline{toc}{section}{Author Contribution Statement (CRediT)}
\section*{Author Contribution Statement (CRediT)}
\begin{tabular}{ll}
\textbf{Philipp Krah:} & concept, method, implementation, writing original draft \\
\textbf{Steffen Büchholz:} & implementation of neural networks, reviewing \& editing\\
\textbf{Matthias Häringer:} & computations of the Bunsen flame, reviewing \& editing\\
\textbf{Julius Reiss:} & proposed initial idea and application to ARD-systems,\\ & supervision, reviewing \& editing
\end{tabular}

%\section*{Code and Data Availability}
% Software can be found under \cite{WABBIT_github,WABBIT_pythongithub}. All scripts to reproduce the results are available at \cite{WABBIT_convergence}.
\bibliographystyle{unsrtnat}
\bibliography{references.bib}   % name your BibTeX data base

\newpage
\appendix
% !TeX root = ../main.tex
\newcommand{\bxi}{\boldsymbol{\xi}}
\newcommand{\bPsi}{\boldsymbol{\Psi}}
\newcommand{\f}{\vec{f}}
\newcommand{\residual}{\vec{r}}
\renewcommand{\d}{\mathrm{d}}
\newcommand{\Lagrange}{\mathcal{L}}
\newcommand{\shrink}{\mathcal{S}}

%%%%%%%%%%%%%%%%%%%%%%%%%%%%%%%%%%%%%%%%%%%%%%%%%%%%%%%%%%%%%%%%%%%%%%%%%%%%%%%%%%
\section{Details on the Autoencoders Network Architecture and training hyperparameters}
\label{appx:AutoencoderArchitecture}

In this section we provide detailed information on the architecture and training hyperparameters for the autoencoder networks. For both autoencoder variants (NN and FTR-NN), the encoder architecture $g_\text{enc}:\mathbb{R}^M\to\mathbb{R}^r$ is the same. Its task is to encode the spatial field $\vec{q}\in\mathbb{R}^M$ into a latent space $\vec{a}\in\mathbb{R}^r$. It consists of four convolutional layers, each followed by a ELU activation and a batch normalization layer. After flattening the output, two fully connected layers follow, with another ELU activation and batch normalization layer in between. The output of the second fully connected layer represents the latent space with $r$ degrees of freedom and is not activated. A summary of the encoder architecture is listed in \cref{tabl-appx:Enc-details}.
The decoder, $g_\text{dec}:\mathbb{R}^r\to\mathbb{R}^M$ maps the latent representation back to the spatial domain. 

There are two different decoders used in this paper labeled NN and FTR-NN autoencoder.
The NN decoder mirrors the encoder architecture, using transposed convolutional layers instead of convolutional layers. The FTR-NN decoder consists of only a single fully connected layer with no bias with $M$ (number of grid points) output channels. It applies a simple Matrix multiplication $\matr{W}\vec{a}$, where $\vec{a}$ is the vector with the latent representations and $\matr{W}$ is the learnable weight matrix of the layer. Afterwards the resulting output $\vec{\phi}=\matr{W}\vec{a}$ is reshaped into the spatial domain. In analogy to the FTR ansatz $\vec{q}\approx\tilde{\vec{q}}=f(\vec{\phi})$, both networks are activated with the physics dependent front function $f$ in the output layer. The layer details for both decoder networks are listed in \cref{tabl-appx:Dec-details}.

After splitting the data by taking every other time step into a training set and a test set, each network was trained using the ADAM optimizer with a learning rate of $0.0025$ for up to $2 \cdot 10^4$ iterations, using all training samples as input batch. Every 500 iterations, the performance is tested on the test set. The network parameters that yield the best test results are saved.

\begin{figure}[htp!]
%$\vec{F)(\vec{q},t,\mu)$
    \begin{center}
    \setlength{\figureheight}{0.35\linewidth}
    \setlength{\figurewidth}{0.35\linewidth}
    % This file was created with tikzplotlib v0.9.17.
\begin{tikzpicture}

\begin{groupplot}[group style={group size=3 by 2}]
\nextgroupplot[
height=\figureheight,
tick pos=left, ytick style={draw=none}, xtick style={draw=none}, yticklabels={,,}, xticklabels={,,},
title={\(\displaystyle \lambda_{\mathrm{smooth}}=10^{-3}\)},
width=\figurewidth,
xmin=0, xmax=129,
ylabel={FTR-NN},
ymin=0, ymax=129
]
\addplot graphics [includegraphics cmd=\pgfimage,xmin=0, xmax=129, ymin=0, ymax=129] {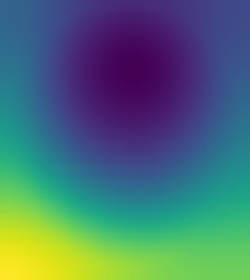};

\nextgroupplot[
height=\figureheight,
tick pos=left,
ytick style={draw=none},
xtick style={draw=none},
yticklabels={,,},
xticklabels={,,},
title={\(\displaystyle \lambda_{\mathrm{smooth}}=10^{-6}\)},
width=\figurewidth,
xmin=0, xmax=129,
ymin=0, ymax=129
]
\addplot graphics [includegraphics cmd=\pgfimage,xmin=0, xmax=129, ymin=0, ymax=129] {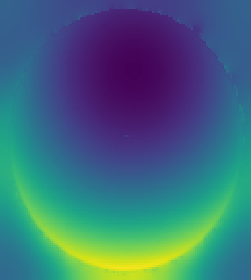};

\nextgroupplot[
colorbar,
colorbar style={ytick={-100,0,100,200,300},yticklabels={
  \(\displaystyle {\ensuremath{-}100}\),
  \(\displaystyle {0}\),
  \(\displaystyle {100}\),
  \(\displaystyle {200}\),
  \(\displaystyle {300}\)
},ylabel={}},
colormap/viridis,
height=\figureheight,
point meta max=205.735382080078,
point meta min=-22.6586303710938,
tick pos=left, ytick style={draw=none}, xtick style={draw=none}, yticklabels={,,}, xticklabels={,,},
title={\(\displaystyle \lambda_{\mathrm{smooth}}=10^{-9}\)},
width=\figurewidth,
xmin=0, xmax=129,
ymin=0, ymax=129
]
\addplot graphics [includegraphics cmd=\pgfimage,xmin=0, xmax=129, ymin=0, ymax=129] {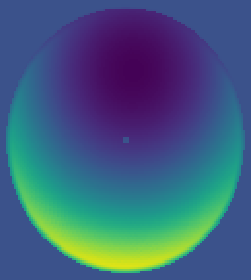};

\nextgroupplot[
height=\figureheight,
tick pos=left, ytick style={draw=none}, xtick style={draw=none}, yticklabels={,,}, xticklabels={,,},
width=\figurewidth,
xmin=0, xmax=129,
ylabel={NN},
ymin=0, ymax=129
]
\addplot graphics [includegraphics cmd=\pgfimage,xmin=0, xmax=129, ymin=0, ymax=129] {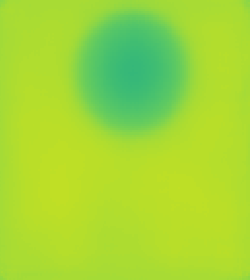};

\nextgroupplot[
height=\figureheight,
tick pos=left, ytick style={draw=none}, xtick style={draw=none}, yticklabels={,,}, xticklabels={,,},
width=\figurewidth,
xmin=0, xmax=129,
ymin=0, ymax=129
]
\addplot graphics [includegraphics cmd=\pgfimage,xmin=0, xmax=129, ymin=0, ymax=129] {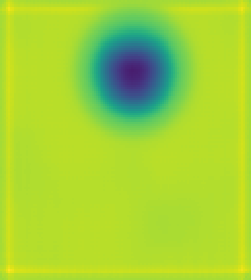};

\nextgroupplot[
colorbar,
colorbar style={ytick={-100,-50,0,50},yticklabels={
  \(\displaystyle {\ensuremath{-}100}\),
  \(\displaystyle {\ensuremath{-}50}\),
  \(\displaystyle {0}\),
  \(\displaystyle {50}\)
},ylabel={}},
colormap/viridis,
height=\figureheight,
point meta max=21.4449024200439,
point meta min=-91.2352828979492,
tick pos=left, ytick style={draw=none}, xtick style={draw=none}, yticklabels={,,}, xticklabels={,,},
width=\figurewidth,
xmin=0, xmax=129,
ymin=0, ymax=129
]
\addplot graphics [includegraphics cmd=\pgfimage,xmin=0, xmax=129, ymin=0, ymax=129] {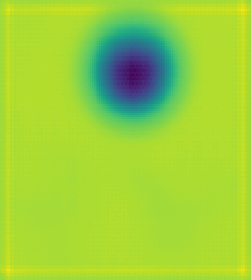};
\end{groupplot}

% \begin{axis}[
% axis y line=right,
% height=\figureheight,
% tick align=outside,
% width=\figurewidth,
% xmin=0, xmax=1,
% xtick pos=left,
% y grid style={white!69.0196078431373!black},
% ymin=-22.6586303710938, ymax=205.735382080078,
% ytick pos=right,
% ytick style={color=black},
% ytick={-100,0,100,200,300},
% yticklabel style={anchor=west},
% yticklabels={
%   \(\displaystyle {\ensuremath{-}100}\),
%   \(\displaystyle {0}\),
%   \(\displaystyle {100}\),
%   \(\displaystyle {200}\),
%   \(\displaystyle {300}\)
% }
% ]
% \path [draw=white, fill=white, line width=0.004pt]
% ;
% \addplot graphics [includegraphics cmd=\pgfimage,xmin=0, xmax=1, ymin=-22.6586303710938, ymax=205.735382080078] {disc_smoothness-006.png};
% \end{axis}

% \begin{axis}[
% axis y line=right,
% height=\figureheight,
% tick align=outside,
% width=\figurewidth,
% xmin=0, xmax=1,
% xtick pos=left,
% y grid style={white!69.0196078431373!black},
% ymin=-91.2352828979492, ymax=21.4449024200439,
% ytick pos=right,
% ytick style={color=black},
% ytick={-100,-50,0,50},
% yticklabel style={anchor=west},
% yticklabels={
%   \(\displaystyle {\ensuremath{-}100}\),
%   \(\displaystyle {\ensuremath{-}50}\),
%   \(\displaystyle {0}\),
%   \(\displaystyle {50}\)
% }
% ]
% \path [draw=white, fill=white, line width=0.004pt]
% ;
% \addplot graphics [includegraphics cmd=\pgfimage,xmin=0, xmax=1, ymin=-91.2352828979492, ymax=21.4449024200439] {disc_smoothness-007.png};
% \end{axis}

\end{tikzpicture}
    \end{center}
    \caption{Color plot of one snapshot of the  FTR-NN and NN levelset field $\phi$ using three degrees of freedom and different smoothness strength $\smoothness$. The smoothness parameter $\smoothness$ controls the strength of the smoothness constraint \cref{eq:NN-reg}.} 
    \label{fig-apx:disc-smoothness}
\end{figure}

\begin{table}[htp!]
    \centering
    \begin{tabular}{l c c c c}
    \toprule
    Layer & Details  &  &\\ 
    \toprule
     & Input  & Output & Kernel\\ 
     & channels & channels & Size & Stride\\ 
    \cmidrule(lr){2-5}
    Input of $q$ ($M$ grid points) & 1 \\
    2D Convolution & 1 & 8 & 5 & 1 \\
    ELU + 2D BatchNorm \\
    2D Convolution  & 8 & 16 & 5 & 2 \\
    ELU + 2D BatchNorm \\
    2D Convolution  & 16 & 32 & 5 & 2 \\
    ELU + 2D BatchNorm \\
    2D Convolution  & 32 &  16 & 5 & 2 \\
    ELU + 2D BatchNorm \\
    Flatten Spatialy \\
    \cmidrule(lr){1-3}
    Fully Connected & $16 \cdot \tilde{M}$ & 512 \\
    ELU + 1D BatchNorm \\
    Fully Connected & 512 & $r$ \\
    Output of latent representation $a$ & & $r$ \\
    \bottomrule
    \end{tabular}
			\caption{Encoder network details. $\tilde{M}$ describes the number of remaining spatial grid points after all convolutional layers are applied. Each convolutional layer reduces the spatial resolution in each spatial direction by $N_\mathrm{out} = \left(N_\mathrm{in} - \text{kernel size} \right) / \text{stride} + 1$}.
			\label{tabl-appx:Enc-details}
\end{table}

\begin{table}[htp!]
    \centering
    \begin{tabular}{l c c c c}
    \toprule
    Layer & Details  &  &\\ 
    \toprule
     & Input  & Output & Kernel\\ 
     & channels & channels & Size & Stride\\ 
    \cmidrule(lr){2-5}
    Input of latent representation $a$ & $r$ \\
    Fully Connected & $r$ & 512 \\
    ELU + 1D BatchNorm \\
    Fully Connected & 512 &$16 \cdot \tilde{M}$ \\
    ELU \\
    Unflatten Spatialy & & 16\\
    \cmidrule(lr){1-3}
    2D BatchNorm \\
    2D Transposed Convolution & 16 & 32 & 5 & 2 \\
    ELU + 2D BatchNorm \\
    2D Transposed Convolution & 32 & 16 & 5 & 2 \\
    ELU + 2D BatchNorm \\
    2D Transposed Convolution & 16 & 8 & 5 & 2 \\
    ELU + 2D BatchNorm \\
    2D Transposed Convolution & 8 & 1 & 5 & 1 \\
    Output of $\phi$ ($M$ grid points) \\
    \bottomrule
    \end{tabular}
			\caption{NN decoder network details}
			\label{tabl-appx:Dec-details}
\end{table}

%%%%%%%%%%%%%%%%%%%%%%%%%%%%%%%%%%%%%%%%%%%%%%%%%%%%%%%%%%%%%%%%%%%%%%%%%%%%%%%%%%
\section{Simulation details of the 1D advection PDE and reaction diffusion PDE example}
\label{appx:1d-advection-example}
In this section we give additional details on the two PDE-examples with analytical solution. Namely, the PDE-example for
\begin{equation}
		    \label{eq-appx:advec-system}
			\text{advection}\quad\quad \begin{cases}
				0 &=\partial_t q - u(t)\partial_x q\\
				q(x,t) &= f(\abs{x-u(t)}-2)
			\end{cases}\,,
		\end{equation}

shown in \cref{fig:1D-advection} and
\begin{equation}
    \label{eq-appx:1Dreact_test}
    \text{reaction-diffusion} \quad
    \begin{cases}
    0 &=\partial_t q - \partial_{xx} q + \frac{8}{\delta^2} q^2(q-1)\\
    q(x,t) &= f(\frac{\abs{x}-2-t/\delta}{\delta})
    \end{cases}\,.
\end{equation}
In both examples we use central finite difference of 6th order with periodic boundary conditions and an explicit Runge-Kutta integration method of 5th(4th) order for adaptive time stepping of the FOM and ROM ODE-system \cite{DormandPrince1980}. 
The numerical parameters for the FTR-decomposition and discretization are stated in \cref{tabl-appx:PDE-details}.

%%%%%%%%%%%%%%%%%%%%%%%%%%%%%%
\begin{table}[hbp!]
    \centering
        \begin{tabular}{l c c}
    \toprule
    property                    & advection  & reaction-diffusion  \\ 
    \midrule
    %------------------------
    \textbf{FOM - parameters}\\
    %------------------------
    Simulation time $T$             & $2.5$  & $ 1$ \\
    Domain     $\domain$        & $[-20,20]$ & $[-15,15]$ \\
    Grid resolution $M$             & $1000$ & $ 4000$ \\
    \midrule
    %------------------------
    \textbf{ROM - parameters}\\
    %------------------------
    Number of snapshots          & 202 & 202 \\
    FTR iterations               & 3000 & 8000 \\
    FTR step width $\tau$        & 1    & 4 \\
    front function $f(x)$              & $ 0.5(1-\tanh(2.5 x))$ & $ 0.5(1-\tanh(x))$\\
% 				Front function & $f(x) = \sigmoid (x\lambda)$\\
% 				Front thickness & $\lambda=100$\\
    \bottomrule
    \end{tabular}
	\caption{Parameters of the 1D advection and reaction-diffusion simulations and the decomposition procedure (\cref{alg:iterFTR})}
	\label{tabl-appx:PDE-details}
\end{table}

\end{document}